\documentclass[11pt,twoside,english,bibliography=totoc]{article}
\usepackage[T1]{fontenc}
\usepackage[utf8]{inputenc}
\usepackage{geometry}
\usepackage{color}
\usepackage{babel}
\usepackage{amsmath}
\usepackage{amsthm}
\usepackage{amssymb}
\usepackage{stackrel}
\usepackage[all]{xy}
\PassOptionsToPackage{normalem}{ulem}
\usepackage{ulem}
\usepackage[unicode=true,pdfusetitle,
 bookmarks=true,bookmarksnumbered=false,bookmarksopen=false,
 breaklinks=false,pdfborder={0 0 1},backref=page,colorlinks=true]
 {hyperref}
\hypersetup{
 linkcolor=blue,citecolor=blue}

\makeatletter
\numberwithin{equation}{section}
\numberwithin{figure}{section}
\theoremstyle{plain}
\newtheorem{thm}{\protect\theoremname}[section]
\theoremstyle{plain}
\newtheorem{conjecture}[thm]{\protect\conjecturename}
\theoremstyle{plain}
\newtheorem{cor}[thm]{\protect\corollaryname}
\theoremstyle{definition}
\newtheorem{defn}[thm]{\protect\definitionname}
\theoremstyle{plain}
\newtheorem{prop}[thm]{\protect\propositionname}
\theoremstyle{remark}
\newtheorem{rem}[thm]{\protect\remarkname}
\theoremstyle{plain}
\newtheorem{lem}[thm]{\protect\lemmaname}
\theoremstyle{plain}
\newtheorem*{thm*}{\protect\theoremname}

\usepackage{babel}

\usepackage{ytableau}

\usepackage{multirow}

\usepackage{appendix}

\usepackage{bbm}

\usepackage{graphicx}
\newcommand{\Conv}{\mathop{\scalebox{1.5}{\raisebox{-0.2ex}{$\ast$}}}}

\makeatother

\providecommand{\conjecturename}{Conjecture}
\providecommand{\corollaryname}{Corollary}
\providecommand{\definitionname}{Definition}
\providecommand{\lemmaname}{Lemma}
\providecommand{\propositionname}{Proposition}
\providecommand{\remarkname}{Remark}
\providecommand{\theoremname}{Theorem}

\begin{document}
\global\long\def\F{\mathbb{\mathbb{\mathbf{F}}}}%
 
\global\long\def\rk{\mathbb{\mathrm{rk}}}%
 
\global\long\def\crit{\mathbb{\mathrm{Crit}}}%
 
\global\long\def\Hom{\mathrm{Hom}}%
 
\global\long\def\defi{\stackrel{\mathrm{def}}{=}}%
 
\global\long\def\tr{{\cal T}r }%
 
\global\long\def\id{\mathrm{id}}%
 
\global\long\def\Aut{\mathrm{Aut}}%
 
\global\long\def\wl{w_{1},\ldots,w_{\ell}}%
 
\global\long\def\alg{\mathrm{alg}}%
 
\global\long\def\ff{\stackrel{*}{\le}}%
 
\global\long\def\mobius{M\dacute{o}bius}%
 
\global\long\def\chimax{\chi^{\mathrm{max}}}%
 
\global\long\def\fq{\mathbb{F}_{q}}%
 
\global\long\def\gl{\mathrm{GL}}%
 
\global\long\def\glm{\mathrm{GL}_{m}\left(\k\right)}%
 
\global\long\def\gln{\mathrm{GL}_{N}\left(q\right)}%
 
\global\long\def\k{K}%
 
\global\long\def\fix{\mathrm{fix}}%
 
\global\long\def\cl{\mathrm{cl}}%
 
\global\long\def\scl{\mathrm{scl}}%
 
\global\long\def\sqlh{\mathrm{sql^{*}}}%
 
\global\long\def\ssql{\mathrm{ssql}}%
 
\global\long\def\sp{\mathrm{s\pi}}%
 
\global\long\def\spm{\mathrm{s\pi^{\left(m\right)}}}%
 
\global\long\def\spq{\mathrm{s\pi_{q}}}%
 
\global\long\def\spf{\mathrm{s\pi^{\phi}}}%
 
\global\long\def\decomp{{\cal D}\mathrm{ecomp}_{\twoheadrightarrow}^{2}}%
 
\global\long\def\decompt{{\cal D}\mathrm{ecomp}_{\twoheadrightarrow}^{3}}%
 
\global\long\def\algdecomp{{\cal D}\mathrm{ecomp}_{\mathrm{alg}}^{2}}%
 
\global\long\def\algdecompt{{\cal D}\mathrm{ecomp}_{\mathrm{alg}}^{3}}%
 
\global\long\def\LB{L^{\twoheadrightarrow} }%
\global\long\def\CB{C^{\twoheadrightarrow} }%
\global\long\def\RB{R^{\twoheadrightarrow} }%
 
\global\long\def\irr{\mathrm{Irr} }%
 
\global\long\def\wh{\mathrm{Wh} }%
 
\global\long\def\cod{\mathrm{codom} }%
 
\global\long\def\Sp{\mathrm{Sp}}%
 
\global\long\def\s{\mathrm{S}}%
 
\global\long\def\U{\mathrm{U}}%
 
\global\long\def\O{\mathrm{O}}%
 
\global\long\def\arrm{\mathrm{\overrightarrow{\mu}}}%
 
\global\long\def\arrn{\mathrm{\overrightarrow{\nu}}}%
 
\global\long\def\I{{\cal I}}%
 
\global\long\def\irr{\mathrm{Irr}}%
 
\global\long\def\stairr{\mathrm{StabIrr}}%
 
\global\long\def\triv{\mathrm{triv}}%
 
\global\long\def\chialg{\chi_{\mathrm{prop.alg.}}}%
 
\global\long\def\supp{\mathrm{Supp}}%
 
\global\long\def\std{\mathrm{std}}%

\title{Stable Invariants of Words from Random Matrices II:\\
Formulas and Extensions}
\author{Doron Puder~~~~and~~~~Yotam Shomroni}
\maketitle
\begin{abstract}
Let $w$ be a word in a free group. As was revealed by Magee and Puder
in \cite{MPunitary}, the stable commutator length ($\scl$) of $w$,
a well-known topological invariant, can also be defined in terms of
certain stable Fourier coefficients of $w$-random unitary matrices.
In the first part of the current work \cite{PSh23}, we demonstrated
how this phenomenon is much broader: we proved more instances of such
results and conjectured others. These new results and conjectures
involved other topological invariants (relatives of $\scl$) and different
families of groups. 

In the current paper we further extend and support this theory. We
provide another instance of the theory and prove that the stable primitivity
rank, too, can be expressed in terms of stable Fourier coefficients
of $w$-random elements of groups. We introduce concrete formulas
for stable Fourier coefficients of $w$-random elements in the symmetric
group $S_{N}$ and its generalizations in the form of the wreath products
$G\wr S_{N}$ where $G$ is an arbitrary compact group. We also define
new stable invariants related to these groups, and prove they give
bounds to many of the stable Fourier coefficients. 

As an aside, we generalize to tuples of words a result of Puder and
Parzanchevski \cite{PP15} about the expected number of fixed points
of $w$-random permutations.
\end{abstract}
\tableofcontents{}

\section{Introduction\label{sec:Introduction}}

In \cite{MPunitary}, Magee and the first named author revealed a
surprising connection between random matrices and stable commutator
length ($\scl$)\footnote{The exact definition of $\scl$ is not very important at the moment.
It is enough to know it is a well-known and well-studied invariant:
there is even a whole book dedicated to it \cite{calegari2009scl}. }: for any word $w$ in a free group $\F$,\footnote{Throughout the paper, we assume without loss of generality that $\F$
has finite rank $r\in\mathbb{Z}_{\ge1}$.} $\scl(w)$ can be expressed in terms of stable Fourier coefficients
of $w$-random unitary matrices. Here, a $w$-random element of a
compact group $G$ is obtained by substituting the letters of $w$
with independent, Haar-random elements of $G$. For example, if $w=abab^{-2}$,
a $w$-random element of $G$ is $ghgh^{-2}$ where $g$ and $h$
are independent, Haar-random elements of $G$.\footnote{Whenever we spell out a word, unless stated otherwise, we use letters
which are assumed to be distinct elements of a basis of the ambient
free group.} By a Fourier coefficient we mean the expected value of an irreducible
character, namely, the expected value of the trace of a continuous
irreducible unitary representation. \emph{Stable} Fourier coefficients
correspond to \emph{stable} irreducible characters, in the sense of
\cite{church2015fi}.

In the prequel \cite{PSh23} of the current paper, we suggest vast
generalizations of this phenomenon. The main point is that the random-matrix
side of the story can be tweaked by considering groups other than
$U(N)$ or even just different families of Fourier coefficients, and
\cite{PSh23} presents results and conjectures as to what relatives
of $\scl$ may replace it on the other side of the equations. That
paper deals with various families of groups, including, inter alia,
unitary groups, orthogonal groups, and matrix groups over finite fields.
In the current paper we focus on the symmetric groups $\s_{\bullet}\defi\{\s_{N}\}_{N\ge1}$
(also discussed in \cite{PSh23}) and, more generally, the wreath
products $G\wr\s_{\bullet}\defi\{G\wr S_{N}\}_{N\ge1}$ where $G$
is an arbitrary compact group (\cite{PSh23} discusses only the case
where $G$ is cyclic or $S^{1}$). Our main results are a proof that
the stable primitivity rank, a relative of $\scl$ introduced by Wilton
\cite[Def.~10.6]{wilton2024rational}, can be defined, too, in terms
of Fourier coefficients of $w$-random elements of groups (Theorem
\ref{thm:spi profinite}), and explicit formulas for stable Fourier
coefficients of $S_{\bullet}$ (Theorem \ref{thm:formula for stable irreps of S_N})
and of $G\wr S_{\bullet}$ (Theorem \ref{thm:formula for stable irreps of G wr S_bullet}).
We also introduce additional relatives of $\scl$ and give bounds
on many of the stable Fourier coefficients. Some of these results
are used in \cite{PSh23}.

The study of $w$-random elements of groups has its roots in works
by free probabilists and combinatorialists in the 1980's and 1990's.
These works usually consider the limit values of some Fourier coefficients,
and do not distinguish between different words as long as both are
non-powers. A more delicate study of the subject was initiated more
recently in \cite{puder2014primitive,PP15}. For $N\ge2$, denote
by $\mathrm{std}_{N}$ the character of the standard, $(N-1)$-dimensional
representation of $S_{N}$, so $\mathrm{std_{N}(\sigma)+1}$ is the
number of fixed points of $\sigma\in S_{N}$. Denote by $\mathbb{E}_{w}[\mathrm{std}_{N}]$
the expected value of $\mathrm{std}_{N}$ evaluated on a $w$-random
permutation in $S_{N}$. It was shown already in \cite{broder1987second}
that if $w\ne1$ is not a proper power then $\mathbb{E}_{w}[\mathrm{std}_{N}]\stackrel{N\to\infty}{\longrightarrow}0$
(and that this limit is a positive integer for proper powers), and
in \cite{nica1994number} that $\mathbb{E}_{w}[\mathrm{std}_{N}]$
coincides with some rational function $f_{w}(N)\in\mathbb{Q}(N)$
for large enough $N$. The precise degree of $f_{w}$ was studied
in \cite{puder2014primitive,PP15} and shown to be determined by $\pi(w)$,
the primitivity rank of $w\in\F$, a notion introduced in \cite{puder2014primitive}.
It is defined by
\begin{equation}
\pi\left(w\right)\defi\min\left\{ \rk H\,\middle|\,w\in H\le\F~\mathrm{and}~w~\mathrm{non\textnormal{-}primitive~in}~H\right\} .\label{eq:def of pi}
\end{equation}
By \cite[Thm.~1.8]{PP15}, 

\begin{equation}
\mathbb{E}_{w}\left[\mathrm{std}_{N}\right]=N^{1-\pi(w)}\left(c_{w}+O\left(\frac{1}{N}\right)\right),\label{eq:PP15}
\end{equation}
where $c_{w}\in\mathbb{Z}_{\ge0}$ is the number of subgroups $H\le\F$
of rank $\pi(w)$ in which $w$ is not primitive.\footnote{To clarify \eqref{eq:def of pi} and \eqref{eq:PP15}, recall that
by Nielsen-Schreier theorem, every subgroup of a free group is free.
An element of a free group is called primitive if it is contained
in some basis. The minimum over the empty set is $\infty$, and by
\cite[Lem.~4.1]{puder2014primitive}, $\pi\left(w\right)=\infty$
if and only if $w$ is primitive in $\F$. The number $c_{w}$ is
finite by \cite[\S4]{puder2014primitive}.} 

Denote by $\mathrm{std}$ the sequence of irreducible characters $\{\mathrm{std}_{N}\}_{N\ge2}$.
This is only one example of a \emph{stable irreducible character}
of the symmetric groups $S_{\bullet}$, corresponding to the Young
diagrams $(N-1,1$). More generally, every fixed partition $\mu$
gives rise to a stable irreducible character of $S_{\bullet}$ corresponding
to the family of Young diagrams where the structure outside the first
row is fixed and given by $\mu$. (In the case of $\mathrm{std}$,
this partition is $\mu=(1)\vdash1$). We elaborate more in $\S\S$\ref{subsec:reps of G wr S_N}.
We denote this stable irreducible character by $\chi^{\mu[\bullet]}=\{\chi^{\mu[N]}\}_{N\ge|\mu|+\mu_{1}}$.
As explained in \cite[\S\S4.2]{PSh23}, for every non-empty partition
$\mu$, the average $\mathbb{E}_{w}[\chi^{\mu[N]}]$ coincides with
a rational function in $\mathbb{Q}(N)$ for every large enough $N$,
and we define $\beta(w,\chi^{\mu[\bullet]})\in\mathbb{Q}\cup\{\infty\}$
to be the number such that 
\begin{equation}
\mathbb{E}_{w}\left[\chi^{\mu[N]}\right]=\left(\dim\chi^{\mu\left[N\right]}\right)^{-\beta\left(w,\chi^{\mu[\bullet]}\right)}\left(c+O\left(N^{-1}\right)\right),\label{eq:def of beta for S_bullet}
\end{equation}
for some constant $c\ne0$. (In particular, $\beta(w,\chi^{\mu[\bullet]})=\infty$
if and only if $\mathbb{E}_{w}[\chi^{\mu[\bullet]}]=0$ for every
large enough $N$). In this language, \eqref{eq:PP15} yields that
$\beta(w,\chi^{(1)[\bullet]})=\pi(w)-1$.

The above-mentioned result from \cite{MPunitary} states that the
infimum of $\beta(w,\chi)$ over all non-trivial stable polynomial
irreducible characters $\chi$ of $\U(\bullet)=\{\U(N)\}_{N}$ coincides
with $2\scl(w)$. A similar result was obtained in \cite[Thm.~1.5]{PSh23}
for generalized symmetric groups. We conjecture that in the case of
$S_{\bullet}$, the infimum of $\beta(w,\chi)$ over all non-trivial
stable irreducible characters is equal to \emph{$\sp(w)$ -- the
stable primitivity rank }of $w$, see Definition \ref{def:spi}. Roughly,
$\sp(w)$ is equal to the minimal quantity $\frac{\rk H-1}{d}$, where
$H\le\F$ contains in a non-trivial fashion different conjugates of
powers of $w$, of total exponent $d$. The precise conjecture is
the following:
\begin{conjecture}[{\cite[Conj.~1.2]{PSh23}}]
\label{conj:spi as statistical invariant} Let ${\cal I}=\left\{ \chi^{\mu[\bullet]}\right\} _{\mu}$
be the set of stable irreducible characters of $S_{\bullet}$. Then
for every $w\in\F$,
\begin{equation}
\inf_{\triv\ne\chi\in{\cal I}}\beta(w,\chi)=\sp(w).\label{eq:spi-as-statistical-invariant-conj}
\end{equation}
Moreover, the infimum on the left hand side is attained.
\end{conjecture}

While we are still not able to prove this conjecture, we do prove
here, nonetheless, that $\sp(w)$ is equal to the infimum over $\beta(w,\chi)$'s,
only for a different set of stable irreducible characters $\chi$.
These characters come from groups of the form $G\wr S_{\bullet}$,
where $G$ is an arbitrary compact group.  The group $G\wr S_{N}$
can be realized as the group of $N\times N$ matrices with entries
in $G\sqcup\{0\}$ and exactly one non-zero entry in every block and
every row -- see $\S\S$\ref{subsec:reps of G wr S_N}. In \cite{PSh23},
these ideas were developed in the case where $G$ was $C_{m}\defi\{z\in S^{1}\,\mid\,z^{m}=1\}$
for some $1\ne m\in\mathbb{Z}_{\ge0}$ (so $G$ is either finite non-trivial
cyclic or $C_{0}=S^{1})$. 

As we elaborate in $\S\S$\ref{subsec:reps of G wr S_N}, stable irreducible
characters of $G\wr S_{\bullet}$ are classified by finitely-supported
partition-valued functions $\arrm\colon\irr(G)\to{\cal P}$, where
$\irr(G)$ is the set of irreducible characters of $G$, ${\cal P}$
is the set of all partitions of non-negative integers, and the support
$\supp(\arrm)$ of $\arrm$ is the set of irreducible characters which
are mapped to non-empty partitions. We denote by $\chi^{\arrm[\bullet]}\defi\{\chi^{\arrm[N]}\}_{N\ge N_{0}}$
the stable irreducible character corresponding to $\arrm$. For every
$w\in\F$ and every $\arrm$, similarly to \eqref{eq:def of beta for S_bullet},
we may define the number $\beta(w,\chi^{\arrm[\bullet]})$ which controls
the rate of decay of the expectation $\mathbb{E}_{w}[\chi^{\arrm[N]}]$.
\begin{thm}
\label{thm:spi profinite}For $m\ge2$, let ${\cal J}_{m}$ denote
the set of stable non-trivial irreducible characters of $S_{m}\wr S_{\bullet}$
corresponding to $\arrm\colon\irr(S_{m})\to{\cal P}$ with $\triv\notin\supp(\arrm)$.\footnote{When $m=1$ there are no such stable characters.}
Let ${\cal J}=\bigsqcup_{m\ge2}{\cal J}_{m}$. Then,
\[
\sp(w)=\inf_{\chi\in{\cal J}}\beta\left(w,\chi\right).
\]
Moreover, the infimum on the right hand side is attained. 
\end{thm}

Call an invariant $f$ of words \emph{profinite} if $f(w_{1})=f(w_{2})$
whenever $w_{1}$ and $w_{2}$ induce the same measures on all finite
groups \cite[Def.~1.3]{PSh23}. 
\begin{cor}
\label{cor:spi profinite}If $w_{1}$ and $w_{2}$ induce the same
measure on $S_{m}\wr S_{N}$ for all $m$ and $N$, then $\sp(w_{1})=\sp(w_{2})$.
Hence, $\sp$ is a profinite invariant.
\end{cor}

We introduced some further evidence towards Conjecture \ref{conj:spi as statistical invariant}
in \cite[\S\S4.2]{PSh23}. One significant piece of evidence is the
following formula for $\mathbb{E}_{w}\left[\chi^{\mu[N]}\right]$.
This formula, stated as Theorem 4.13 in \cite{PSh23} and proven in
the current paper, may very well be useful for other applications.
The unexplained terms and notions used in the statement of Theorem
\ref{thm:formula for stable irreps of S_N} are defined throughout
the paper\footnote{The notation $\chi^{\mu[N]}$ is defined in $\S\S$\ref{subsec:reps of G wr S_N};
$\algdecompt$ and $C_{\eta_{2}}^{\alg}$ in $\S\S$\ref{subsec:M=0000F6bius-inversions};
and the term \emph{efficient} in Definition \ref{def:efficient}.
The remaining condition are explained in $\S$\ref{sec:Formula-for-stable-chars-of-S_N}.}. 
\begin{thm}
\label{thm:formula for stable irreps of S_N}Let $d\in\mathbb{Z}_{\ge0}$,
let $\mu\vdash d$ be a partition and let $1\ne w\in\F$ be a non-power\footnote{The formula as stated is false for proper powers, but can be modified
a bit to apply to proper powers too --- see Remark \ref{rem:about proper powers}.}. For every $N\ge d+\mu_{1}$ we have 
\begin{equation}
\mathbb{E}_{w}\left[\chi^{\mu\left[N\right]}\right]=\frac{1}{d!}\sum_{\sigma\in S_{d}}\chi^{\mu}\left(\sigma\right)\sum_{\substack{\left(\eta_{1},\eta_{2},\eta_{3}\right)\in\algdecompt\left(w^{\sigma}\to\F\right)\colon\\
\eta_{1}~\mathrm{is~efficient,~and}\\
\cod\left(\eta_{2}\right)~\mathrm{has~no~cycles}
}
}C_{\eta_{2}}^{\mathrm{alg}}\left(N\right).\label{eq:formula for stable Fouriers of S_N}
\end{equation}
\end{thm}

Roughly, the summands $\left(\eta_{1},\eta_{2},\eta_{3}\right)$ in
\eqref{eq:formula for stable Fouriers of S_N} correspond to certain
decompositions of the map of labeled graphs from the graph $\Gamma_{w^{\sigma}}$
of disjoint cycles corresponding to powers of $w$ of total exponent
$d$ (based on the cycle structure of $\sigma$, see Definition \ref{def:Gamma_w}),
to the bouquet $\Omega$ representing $\F$ (see $\S\S$\ref{subsec:Core-graphs}).
The invariant $\sp(w)$ clearly shows up in this formula, as most
summands correspond to diagrams as in Definition \ref{def:spi} of
$\sp(w)$ --- see \cite[\S\S4.2]{PSh23} for more details. We compare
the formula \eqref{eq:formula for stable Fouriers of S_N} to other
formulas for $\mathbb{E}_{w}[\chi^{\mu[N]}]$ in $\S\S$\ref{subsec:about the formula we prove}.
This formula is further generalized to all stable irreducible characters
of $G\wr S_{\bullet}$ for arbitrary compact $G$ in Theorem \ref{thm:formula for stable irreps of G wr S_bullet}.
\medskip{}

Still considering the groups $G\wr S_{\bullet}$, when $\arrm$ is
supported on the trivial character $\triv$ of $G$ and $\arrm(\triv)=\mu$,
the character $\chi^{\arrm[N]}$ is obtained by composing the quotient
$G\wr S_{N}\twoheadrightarrow S_{N}$ with the character $\chi^{\mu[N]}$
of $S_{N}$. In particular, if Conjecture \ref{conj:spi as statistical invariant}
holds, the infimum of $\beta(w,\chi^{\arrm[\bullet]})$ over all possible
$\arrm$ is at most $\sp(w)$.\footnote{In fact, we believe that for any compact group $G$, the infimum $\inf_{\triv\ne\chi}\beta(w,\chi)$
over all non-trivial stable irreducible characters of $G\wr S_{\bullet}$
is equal to $\sp(w)$. When $G=1$, this is Conjecture \ref{conj:spi as statistical invariant}.
For the general case, this follows from \cite[Conj.~4.14]{PSh23}
combined with Theorem \ref{thm:formula for stable irreps of G wr S_bullet}.} On the other extreme, we have the partitions where $\arrm(\triv)=\emptyset$,
such as those from Theorem \ref{thm:spi profinite}. One part of Theorem
\ref{thm:spi profinite} holds for every stable irreducible character
of $G\wr S_{\bullet}$ for arbitrary $G$, as long as $\arrm(\triv)=\emptyset$.
\begin{thm}
\label{thm:spi bounds beta from below when triv outside support}Let
$G$ be a compact group and $\arrm\colon\irr(G)\to{\cal P}$ satisfy
$1\le|\arrm|<\infty$ and $\arrm(\triv)=\emptyset$. Then for all
$w\in\F$,
\[
\sp(w)\le\beta\left(w,\chi^{\arrm[\bullet]}\right).
\]
\end{thm}

In the special case where $G=C_{m}$ for $1\ne m\in\mathbb{Z}_{\ge0}$,
Theorem \ref{thm:spi bounds beta from below when triv outside support}
is stated (without a proof) as Theorem 6.17 in \cite{PSh23}. It is
already known that for arbitrary $\arrm\colon\irr(G)\to{\cal P}$
with $1\le|\arrm|<\infty$ (possibly with $\triv\in\supp(\arrm)$),
$\mathbb{E}_{w}[\chi^{\arrm[N]}]=O(N^{1-\pi(w)})$, and even $\mathbb{E}_{w}[\chi^{\arrm[N]}]=O(N^{-\pi(w)})$
if $|\arrm|\ge2$ and $w$ non-power \cite[Thm.~1.6, stated for finite G]{shomroni2023wreathII}.
Theorem \ref{thm:spi bounds beta from below when triv outside support}
improves on these bounds for most $\arrm$ as long as $\triv\notin\supp(\arrm)$
(or for all such $\arrm$ if Wilton's conjecture \cite[Conj.~4.7]{PSh23}
that $\sp(w)=\pi(w)-1$ holds).

A more specialized family of stable irreducible characters of $G\wr S_{\bullet}$
is the family where $\arrm$ is supported on a single non-trivial
irreducible character $\phi\in\irr(G)$. In this case, we can bound
$\beta(w,\chi^{\arrm[\bullet]})$ from below by the $\phi$-primitivity
rank of $w$, denoted $\spf(w)$, which is fully defined in Definition
\ref{def:spf} below. This new invariant of words is yet another new
relative of $\scl$. It generalizes $\sp$ (which can be thought of
as $\sp^{\triv}$), as well as $\spm$ for $1\ne m\in\mathbb{Z}_{\ge0}$
(\cite[Def.~6.1]{PSh23} or Definition \ref{def:spm} below), which
is the special case of $\spf$ where $G=C_{m}$ and $\phi$ is any
faithful linear character. 
\begin{thm}
\label{thm:spf bounds beta from below when supported on phi}Let $G$
be a compact group and $\arrm\colon\irr(G)\to{\cal P}$ satisfy $\supp(\arrm)=\{\phi\}$
for some non-trivial character $\triv\ne\phi\in\irr(G)$. Then for
all $w\in\F$,
\[
\spf(w)\le\beta\left(w,\chi^{\arrm[\bullet]}\right).
\]
\end{thm}

The lower bound in Theorem \ref{thm:spf bounds beta from below when supported on phi}
improves on the bound in Theorem \ref{thm:spi bounds beta from below when triv outside support}
as $\sp\le\spf$ (Corollary \ref{cor:properties of spf}\eqref{enu:sp <=00003D spf}),
often with strict inequality. 

We should stress that in general, our understanding of $\spf$ is
not as deep as that of $\spm$. In \cite{PSh23} we prove that similarly
to $\scl$, $\spm$ obtains only rational values (or $\infty$), and
admits a gap in $\left(0,1\right)$. (In the proof of the latter fact
in \cite[\S\S6.4]{PSh23} we rely on Theorem \ref{thm:spi bounds beta from below when triv outside support}
when specialized to $C_{m}\wr S_{\bullet}$.) We also show there that
for a faithful $\phi\in\irr(C_{m})$, the infimum of $\beta(w,\chi^{\arrm[\bullet]})$
over all $\arrm$ with support $\left\{ \phi\right\} $ is precisely
$\spm(w)$ \cite[Thm.~1.5]{PSh23}. Unfortunately, we do not know
if these properties hold for $\spf$ for a general $\phi$. In $\S$\ref{sec:spf}
we elaborate some properties that $\spf$ does satisfy for an arbitrary
$\phi$. 

One ingredient for all the above results is two formulas for the expectation
under the $w$-measure of the induction of certain characters in $G\wr S_{\bullet}$
-- Proposition \ref{prop:basic formula for induction f boxtensor triv}
and Corollary \ref{cor:basic formula for induction, with algebraic morphisms}.
These characters include the stable irreducible character $\mathbb{E}_{w}[\chi^{\arrm[N]}]$
whenever $\arrm(\triv)=\emptyset$ (Lemma \ref{lem:if no triv, then stable is induction}). 

An additional ingredient in the proof of Theorem \ref{thm:spi profinite}
is a generalization of the main result of \cite{PP15}, which was
stated in \eqref{eq:PP15} above. Already \cite{PP15} contains a
generalization of \eqref{eq:PP15}, where one considers f.g.~(finitely
generated) subgroups of $\F$ rather than single words \cite[Thm.~1.8]{PP15}.
Here we introduce a different generalization. There is a natural extension
of the primitivity rank $\pi$, which is defined on single words,
to a function defined on finite multisets of words $\{w_{1},\ldots,w_{k}\}\subset\F$.
This function returns the maximal Euler characteristic of a ``proper
algebraic extension'' of $\{w_{1},\ldots,w_{k}\}$ -- see Definition
\ref{def:proper algebraic extension} below. We denote it by $\chialg(w_{1},\ldots,w_{k}$).
For non-trivial single words it satisfies $\chialg(w)=1-\pi(w)$.
In the following theorem, $\mathbb{E}_{S_{N}}[\mathrm{std}(w_{1})\cdots\mathrm{std}(w_{k})]$
is the average, over all $(\sigma_{1},\ldots,\sigma_{r})\in(S_{N})^{r}$,
of the product $\mathrm{std}_{N}(w_{1}(\sigma_{1},\ldots,\sigma_{r}))\cdots\mathrm{std}_{N}(w_{k}(\sigma_{1},\ldots,\sigma_{r}))$.
\begin{thm}
\label{thm:generalized PP15}Let $\{w_{1},\ldots,w_{k}\}\in\F$ be
a multiset of non-trivial words. Then
\[
\mathbb{E}_{S_{N}}\left[\prod_{i=1}^{k}\left(\fix\left(w_{i}\right)-1\right)\right]=\mathbb{E}_{S_{N}}\left[\prod_{i=1}^{k}\mathrm{std}\left(w_{i}\right)\right]=N^{\chialg\left(w_{1},\ldots,w_{k}\right)}\left(c_{w_{1},\ldots,w_{k}}+O\left(\frac{1}{N}\right)\right),
\]
where $c_{w_{1},\ldots,w_{k}}$ is the (finite) number of ``proper
algebraic extensions'' of $\{w_{1},\ldots,w_{k}\}$ of Euler characteristic
$\chialg(w_{1},\ldots,w_{k})$.
\end{thm}

\subsection{Comparison with other formulas for stable Fourier coefficients of
$S_{\bullet}$\label{subsec:about the formula we prove}}

We now compare the formula from Theorem \ref{thm:formula for stable irreps of S_N}
to other known formulas for stable Fourier coefficients of word measures
on $S_{\bullet}$: one appearing in \cite{hanany2020word} and another
in \cite{cassidy2024random}. First, the three formulas originate
from three different ways to express stable irreducible characters
of $S_{\bullet}$. 
\begin{itemize}
\item The paper \cite{hanany2020word} studies the functions $f_{N}^{\nu}\colon S_{N}\to\mathbb{Z}$
where $\nu=(\nu_{1},\ldots,\nu_{\ell})\in{\cal P}$ is a partition,
which are defined by $f_{N}(\sigma)=\prod_{i}(\fix(\sigma^{\nu_{i}}))$.
It shows \cite[Prop.~B.2]{hanany2020word} that the functions $f_{N}^{\nu}$
form a linear basis for the algebra of stable class functions on $S_{\bullet}$,
hence $\chi^{\mu[N]}$ is a linear combination of such functions (with
somewhat complicated coefficients). 
\item Cassidy's work uses the fact that $\chi^{\mu[N]}$ has a non-zero
multiplicity in the decomposition of a tensor power $V^{\otimes k}$
of the defining $N$-dimensional representation $V$ of $S_{N}$,
where $k=|\mu|$ (the character of $V^{\otimes k}$ is $\sigma\mapsto\fix(\sigma)^{k}$),
and relies on a projection formula from $V^{\otimes k}$ to the isotypic
component of $\chi^{\mu[N]}$ that Cassidy derives \cite[Thm.~1.1]{cassidy2023projection}.
\item We use a different approach, relying on a different basis for the
ring of stable class functions of $S_{\bullet}$. This is the basis
\[
\left\{ \mathrm{Ind}_{S_{d}\times S_{N-d}}^{S_{N}}\left(\chi^{\nu}\boxtimes\mathrm{triv}\right)\right\} {}_{d\in\mathbb{Z}_{\ge0},\nu\vdash d}.
\]
The linear combinations expressing stable irreducible characters rely
on Pieri's rule (Proposition \ref{prop:pieri} below). The coefficients
in these linear combinations are much simpler compared with the ones
in \cite[Prop.~B.2]{hanany2020word}.
\end{itemize}
Recall Conjecture \ref{conj:spi as statistical invariant} that $\mathbb{E}_{w}[\chi^{\mu[N]}]=O(\dim(\chi^{\mu[N]}))^{-\sp(w)}$
(with $O$ replaced by $\Theta$ for some $\mu$). Definition \ref{def:spi}
of $\sp$ requires the diagrams to be efficient (see Definition \ref{def:efficient}),
algebraic (Definition \ref{def:algebraic morphisms}), and avoid isomorphisms
on connected components of the codomain (the 'proper' part in being
proper algebraic -- Definition \ref{def:proper algebraic extension}). 
\begin{itemize}
\item In \cite[\S\S6.2]{hanany2020word} a formula is introduced for $\mathbb{E}_{w}[f_{N}^{\nu}]$
using only algebraic morphisms, which together with the above-mentioned
formula \cite[Prop.~B.2]{hanany2020word} yields a formula for $\mathbb{E}_{w}[\chi^{\mu[N]}]$
using only algebraic morphisms. However, these morphisms are not necessarily
efficient nor \emph{proper} algebraic, and $\mathbb{E}_{w}[f_{N}^{\nu}]$
can be of high order of magnitude. Hence, delicate cancellations need
to be established in order to prove Conjecture \ref{conj:spi as statistical invariant}
using this formula.
\item Cassidy derives from the aforementioned \cite[Thm.~1.1]{cassidy2023projection}
another formula for $\mathbb{E}_{w}[\chi^{\mu[N]}]$ \cite[Thm.~4.7]{cassidy2024random}.
In that formula, each summand can be thought of as corresponding to
an efficient morphism, but proper algebraic is replaced by the strictly
weaker property that, roughly, every edge in the codomain is covered
at least twice by the morphism (this is similar to the property assumed
in the ``$w$-cycles conjecture'' established in \cite{helfer2016counting,louder2017stackings}).
Cassidy shows that each summand is of order $O(\dim(\chi^{\mu[N]})^{-1})$,
hence proving that $\mathbb{E}_{w}[\chi^{\mu[N]}]=O(\dim(\chi^{\mu[N]}){}^{-1})$.
However, the formula does not restrict to algebraic morphisms, and
using it to prove Conjecture \ref{conj:spi as statistical invariant}
probably requires, too, delicate cancellations. 
\item In contrast, in the formula from Theorem \ref{thm:formula for stable irreps of S_N}
all decompositions $(\eta_{1},\eta_{2},\eta_{3})$ satisfy that the
morphism $\eta_{1}\colon\Gamma_{w^{\sigma}}\to\Sigma$ is both efficient
and algebraic. Granted, it may not be \emph{proper} algebraic, and
only if it is do we know that the summand $C_{\eta_{2}}^{\alg}(N)$
is of order $O(\dim(\chi^{\mu[N]}))^{-\sp(w)}$. Yet, computer simulations
show that \emph{all} summands in this formula are of order $O(\dim(\chi^{\mu[N]}))^{-\sp(w)}$
-- this is formally conjectured in \cite[Conj.~4.14]{PSh23}, to
which Conjecture \ref{conj:spi as statistical invariant} can be reduced.
Importantly, all summands do not depend on the stable irreducible
character $\chi^{\mu[\bullet]}$, but only on morphisms of graphs
with certain properties. So it seems that this formula may possibly
lead to Conjecture \ref{conj:spi as statistical invariant} with no
additional cancellations. To illustrate, for a primitive word $w$
in which every letter appears at least twice, the sum in the formula
from Theorem \ref{thm:formula for stable irreps of S_N} is empty,
in accordance with the a priori fact that $\mathbb{E}_{w}[\chi^{\mu[N]}]\equiv0$
in this case. The other two formulas mentioned above give non-empty
sums which completely cancel in some non-trivial fashion.\\
As a result, our formula here suggests (together with the conjecture
\cite[Conj.~4.14]{PSh23}) the exact value of the coefficient of the
term of order $N^{-|\mu|\cdot\sp(w)}$ in $\mathbb{E}_{w}[\chi^{\arrm[N]}]$:
this coefficient seems to count $\sp$-extremal diagrams of $w$ of
degree $|\arrm|$ -- see \cite[Prop.~4.15]{PSh23} for the precise
statement.
\end{itemize}

\subsection{Outline and notation}

We begin with some preliminaries in $\S$\ref{sec:Preliminaries}:
in $\S\S$\ref{subsec:Core-graphs} we review the notions of core
graphs, and their free and algebraic morphisms; in $\S\S$\ref{subsec:Word measures of S - known analysis}
we mention some crucial points from previous works on the analysis
of word measures in $S_{\bullet}$, and mostly the different Möbius
inversions and their properties, which we use throughout the paper;
in $\S\S$\ref{subsec:reps of G wr S_N} we summarize the definition
of stable irreducible characters in $G\wr S_{\bullet}$; and in $\S\S$\ref{subsec:The-expectations-and-beta}
the general definition of $\beta(w,\chi)$. Section \ref{sec:generalized PP15}
proves Theorem \ref{thm:generalized PP15}, which is needed later
in the proof of Theorem \ref{thm:spi profinite}. After recalling
the definition of $\sp$, Section \ref{sec:spf} defines $\spf$ and
proves some of its basic properties. In $\S$\ref{sec:Lower-bounds-for-stable-chars}
we first find a formula (Proposition \ref{prop:basic formula for induction f boxtensor triv}
and Corollary \ref{cor:basic formula for induction, with algebraic morphisms})
for $\mathbb{E}_{w}[\mathrm{Ind}_{S_{d}\times S_{N-d}}^{S_{N}}(\chi\boxtimes\triv_{N-d})]$
for any character $\chi\in\irr(S_{d})$, as well as for similar induction
of characters in the more general case of $G\wr S_{\bullet}$. Then,
we use these formulas to prove Theorems \ref{thm:spi bounds beta from below when triv outside support}
and \ref{thm:spf bounds beta from below when supported on phi} giving
upper bounds for certain stable Fourier coefficients. We continue
this analysis in $\S$\ref{sec:sp is profinite}, where we prove Theorem
\ref{thm:spi profinite}. Finally, the formula for induction of characters
is used again in $\S$\ref{sec:Formula-for-stable-chars-of-S_N},
where we first prove Theorem \ref{thm:formula for stable irreps of S_N}
in $\S\S$\ref{subsec:Sketch-of-proof}-\ref{subsec:Finishing-the-proof of S_n formula},
and then the formula for stable irreducible characters of $G\wr S_{\bullet}$
for arbitrary compact group $G$ in $\S\S$\ref{subsec:general formula for G wr S}.

\subsubsection*{Notation}

Throughout the paper, $\F$ is a fixed finite-rank free group. The
bouquet $\Omega$ has wedge point $o$ and $\rk\F$ petals, where
each petal is oriented and identified with an element of some fixed
basis of $\F$, so that $\pi_{1}\left(\Omega,o\right)=\F$. The graphs
$\Gamma_{w}$, $\Gamma_{w^{\nu}}$ and $\Gamma_{w^{\sigma}}$ and
the morphisms $\eta_{w},\eta_{w^{\nu}}$ and $\eta_{w^{\sigma}}$
are introduced in Definition \ref{def:Gamma_w} (here $\nu$ is an
integer partition and $\sigma$ a permutation). Given a graph $\Gamma$,
we denote by $V(\Gamma)$ and $E(\Gamma)$ the vertex-set and edge-set
of $\Gamma$, respectively. A Serre graph is defined in $\S\S$\ref{subsec:Core-graphs}.

We denote $\left[d\right]=\left\{ 1,\ldots,d\right\} $. For a finite
set $S$ and $d_{1},\ldots,d_{k}\in\mathbb{Z}_{\ge0}$ with $\sum d_{i}=|S|$,
we denote by $\binom{S}{d_{1}}$ the set of subsets of $S$ of cardinality
$d_{1}$, and by $\binom{S}{d_{1}\,\ldots\,d_{k}}$ the set of all
ordered set-partitions of $S$ into an ordered $k$-tuple of subsets
$(B_{1},\ldots,B_{k}$) with $B_{1}\sqcup\ldots\sqcup B_{k}=S$ and
$\left|B_{i}\right|=d_{i}$. 

A partition of $n\in\mathbb{Z}_{\ge0}$ is $\mu=(\mu_{1},\ldots,\mu_{\ell})$
where $\mu_{1}\ge\mu_{2}\ge\ldots\ge\mu_{\ell}\ge1$ and $\sum\mu_{i}=n$.
We denote the fact that $\mu$ is a partition of $n$ by $\mu\vdash n$,
and given $\mu$, we let $\left|\mu\right|$ denote $n=\sum\mu_{i}$
and $\ell(\mu)=\ell$ denote the number of parts in $\mu$. Let ${\cal P}_{n}$
denote the set of all partitions of $n$, and ${\cal P}=\bigsqcup_{n\in\mathbb{Z}_{\ge0}}{\cal P}_{n}$.
When $n=0$, we denote by $\emptyset$ the empty partition, and if
$\mu=\emptyset$ we set $\mu_{1}=0$ (so that $\mu_{1}$ it defined
for every $\mu\in{\cal P}$). For $N\ge|\mu|+\mu_{1}$, denote by
$\mu[N]$ the partition of $N$ corresponding to a Young diagram with
$\mu$ outside the first row -- see \eqref{eq:mu=00005BN=00005D}.
The stable irreducible character of $S_{\bullet}$ corresponding to
$\mu$ is $\chi^{\mu[\bullet]}\defi(\chi^{\mu[N]})_{N\ge|\mu|+\mu_{1}}$.

Throughout the paper, $G$ denotes an arbitrary compact group. We
denote by $\irr(G)$ the set of irreducible complex characters of
$G$ corresponding to irreducible continuous complex representations.
We sometimes let $G_{N}$ denote the wreath product $G\wr S_{N}$.
As explained in $\S\S$\ref{subsec:reps of G wr S_N}, for any compact
group $G$ the irreducible representations of $G\wr S_{N}$ correspond
to $\arrm\colon\irr(G)\to{\cal P}$ with $|\arrm|\defi\sum_{\phi\in\irr(G)}\arrm(\phi)=N$.
The corresponding character is denoted $\chi^{\arrm}$. For $N\ge|\arrm|+\arrm(\triv)_{1}$,
denote by $\arrm[N]\colon\irr(G)\to{\cal P}$ the map with $|\arrm[N]|=N$
which is identical to $\arrm$ outside the trivial character and such
that $\arrm[N](\triv)$ had $\arrm(\triv)$ outside the first row
-- see $\S\S$\ref{subsec:reps of G wr S_N}. The stable irreducible
character of $G\wr S_{\bullet}$ corresponding to $\arrm$ is $\chi^{\arrm[\bullet]}\defi(\chi^{\arrm[N]})_{N\ge|\arrm|+\arrm(\triv)_{1}}$.

For two functions $f\colon\mathbb{Z}_{\ge N_{0}}\to\mathbb{C}$ and
$g\colon\mathbb{Z}_{\ge N_{0}}\to\mathbb{R}$, we write $f=O(g)$
if there are a constant $c>0$ and $N\in\mathbb{Z}$ such that for
every $n\ge N$, $\left|f(n)\right|\le c|g(n)|$. We write $f=\Theta(g)$
if $f=O(g)$ and $g=O(f)$.

\subsection*{Acknowledgments}

This work was supported by the European Research Council (ERC) under
the European Union’s Horizon 2020 research and innovation programme
(grant agreement No 850956), by the Israel Science Foundation, ISF
grants 1140/23, by the National Science Foundation under Grant No.
DMS-1926686, as well as by the Kovner Member Fund at the IAS, Princeton.

\section{Preliminaries\label{sec:Preliminaries}}

\subsection{Serre graphs, core graphs and algebraic and free morphisms\label{subsec:Core-graphs}}

It will be occasionally useful for us to think of graphs in Serre's
terminology. A \textbf{Serre graph} $\Gamma$ consists of a vertex
set $V\left(\Gamma\right)$ and an edge set $E\left(\Gamma\right)$,
which contains two oriented edges for every geometric edge: every
$e\in E\left(\Gamma\right)$ comes with its inverse edge $\overline{e}\in E\left(\Gamma\right)$,
so $e\mapsto\overline{e}$ is a fixed-point-free involution. The map
$\iota=\iota_{\Gamma}\colon E\left(\Gamma\right)\to V\left(\Gamma\right)$
maps every oriented edge to its starting point, so the edge $e$ connects
$\iota\left(e\right)$ with $\iota\left(\overline{e}\right)$.

The term core graphs is used throughout this paper in the following
sense.
\begin{defn}
\label{def:core graphs}A \textbf{core graph} is a finite simplicial
graph with no leaves and no isolated vertices, namely, with all vertices
of degree $\ge2$. A core graph is not necessarily connected, and
may be empty.
\end{defn}

Consider the bouquet \marginpar{$\Omega$}$\Omega=\bigwedge^{r}\mathbb{S}^{1}$
with $r$ oriented petals identified with the elements of a basis
of $\F$, so $\pi_{1}(\Omega)=\F$. When a core graph $\Gamma$ is
connected and comes equipped with a graph-morphism $\eta\colon\Gamma\to\Omega$
which is an immersion (namely, a locally injective map) into $\Omega$,
it is called a Stallings core graph, after \cite{stallings1983topology}.
In this case, it represents $\eta_{*}(\pi_{1}(\Gamma))$ -- a conjugacy
class of a f.g.~(finitely generated) subgroup of $\F$. A not-necessarily-connected
core graph equipped with an immersion to $\Omega$ is called a \emph{multi
core graph }in\emph{ }\cite{hanany2020word}, and represents a multiset
of conjugacy classes of non-trivial f.g.~subgroups. 

Conversely, if $H\le\F$ is a non-trivial f.g.~subgroup, there is
a unique Stallings core graph (with an immersion to $\Omega$) corresponding
to the conjugacy class $H^{\F}$. Throughout the paper, we will need
multi core graphs corresponding to a multiset of powers of $w$:
\begin{defn}
\label{def:Gamma_w}Let $w\in\F$. If $w\ne1$, denote by \marginpar{$\eta_{w},\Gamma_{w}$}$\eta_{w}\colon\Gamma_{w}\to\Omega$
the Stallings core graph corresponding to $\left\langle w\right\rangle $.
For any partition $\nu=\left(\nu_{1},\ldots,\nu_{\ell}\right)\vdash d$
($d\in\mathbb{Z}_{\ge0}$), denote by \marginpar{$\eta_{w^{\nu}},\Gamma_{w^{\nu}}$}$\eta_{w^{\nu}}\colon\Gamma_{w^{\nu}}\to\Omega$
the multi core graph corresponding to the multiset $\left\{ \left\langle w^{\nu_{1}}\right\rangle ^{\F},\ldots,\left\langle w^{\nu_{\ell}}\right\rangle ^{\F}\right\} $.
If $\nu=\emptyset\vdash0$, then $\Gamma_{w^{\nu}}$ is the empty
graph. If $w=1$, then $\Gamma_{w^{\nu}}$ is the empty graph for
all $\nu$. For $\sigma\in S_{d}$, we denote by \marginpar{$\eta_{w^{\sigma}},\Gamma_{w^{\sigma}}$}$\eta_{w^{\sigma}}\colon\Gamma_{w^{\sigma}}\to\Omega$
the multi core graph corresponding to $w^{\nu}$, where $\nu$ is
the partition of $d$ corresponding to the cycle structure of $\sigma$.
\end{defn}

Topologically, if $w\ne1$ then $\Gamma_{w^{\nu}}$ is a disjoint
union of $\ell(\nu)$ cycles. Note that the number of edges in $\Gamma_{w}$
is the length of the cyclic reduction of $w$. See \cite[Fig.~3.1]{hanany2020word}
for an illustration of $\Gamma_{w^{\left(3,2,1,1\right)}}$.

Let $H\le\F$ be free groups. We say that $\F$ is an \emph{algebraic}
\emph{extension} of $H$ if there are no intermediate proper free
factors of $\F$. This concept goes back to Takahasi \cite{takahasi1951note}
and was coined in \cite[Def.~11.1]{KM02}. It is a simple observation
that if $w\ne1$, then $\pi\left(w\right)$ from \eqref{eq:def of pi}
is the smallest rank of a proper algebraic extension of $\left\langle w\right\rangle $.
In \cite[\S4]{hanany2020word} this notion was extended to multiple
subgroups and to morphisms of graphs. The precise definition involves
the notion of a free morphism of core graphs.
\begin{defn}[free morphisms]
\label{def:free morphism}\cite[Def.~4.2]{hanany2020word} Let $\eta\colon\Gamma\to\Delta$
be a morphism of core graphs. Assume first that $\Delta$ is connected
and let $\Gamma_{1},\ldots,\Gamma_{\ell}$ be the connected components
of $\Gamma$. We say that $\eta$ is a \emph{free} morphism if there
is a choice of a subgroup in the conjugacy class $\eta_{*}(\pi_{1}(\Gamma_{i}))$
in $\pi_{1}(\Delta)$ for all $i$, there exist a subgroup $K\le\pi_{1}(\Delta)$
with 
\[
\pi_{1}(\Delta)=\left(\Conv_{i=1}^{\ell}\eta_{*}\left(\pi_{1}\text{\ensuremath{\left(\Gamma_{i}\right)}}\right)\right)\ast K.
\]
If $\Delta$ is not connected, we call $\eta$ free if it is free
in every connected component of $\Delta$, namely, if $\eta\mid_{\eta^{-1}(\Delta')}\colon\eta^{-1}(\Delta')\to\Delta'$
is free for every connected component $\Delta'$ of $\Delta$.
\end{defn}

\begin{defn}[algebraic morphisms]
\label{def:algebraic morphisms}\cite[Def.~4.6]{hanany2020word}
A morphism $\eta\colon\Gamma\to\Delta$ of core graphs is called \emph{algebraic}
if whenever $\Gamma\stackrel{\eta_{1}}{\to}\Sigma\stackrel{\eta_{2}}{\to}\Delta$
is a decomposition of $\eta$ with $\eta_{2}$ free, we have that
$\eta_{2}$ is an isomorphism.
\end{defn}

Algebraic morphisms can be roughly described also as follows. First,
the morphism $\eta\colon\Gamma\to\Delta$ is algebraic if and only
if it is algebraic when restricted to any connected component of $\Delta$.
If $\Delta$ is connected and $\Gamma_{1},\ldots,\Gamma_{\ell}$ are
the connected components of $\Gamma$, then $\eta\colon\Gamma\to\Delta$
is algebraic if there is no free splitting of $\pi_{1}(\Delta)$ relative
to $\eta_{*}(\pi_{1}(\Gamma_{1})),\ldots,\eta_{*}(\pi_{1}(\Gamma_{1}))$,
namely, if there is no non-trivial free decomposition $\pi_{1}(\Delta)=J_{1}*J_{2}$
so that for each $i=1,\ldots,\ell$, the group $\eta_{*}(\pi_{1}(\Gamma_{i}))$
can be conjugated into $J_{1}$ or $J_{2}$. The latter description
is not accurate in some degenerate cases: for example, the morphism
from the empty core graph to a circle graph is free and therefore
\emph{not} algebraic.

An important fact we use below is that every immersion of core graphs
$\eta\colon\Gamma\to\Delta$ admits a unique decomposition, up to
isomorphism, of the form 
\begin{equation}
\Gamma\stackrel{\eta_{\text{alg}}}{\longrightarrow}\Sigma\stackrel{\eta_{\text{free}}}{\longrightarrow}\Delta\label{eq:algebraic-free decomposition}
\end{equation}
where $\eta_{\text{alg}}$ is algebraic, $\eta_{\text{free}}$ is
free and $\eta=\eta_{\text{free}}\circ\eta_{\text{alg}}$ \cite[Thm.~4.9]{hanany2020word}.

\subsection{Previous analysis of word measures on $S_{\bullet}$\label{subsec:Word measures of S - known analysis}}

The current subsection gathers some definitions and results, mostly
from \cite{puder2014primitive,PP15,hanany2020word}, that we need
in order to compile and prove the formulas in the paper, including
the one in Theorem \ref{thm:formula for stable irreps of S_N}. 

\subsubsection{A geometric interpretation of $\mathbb{E}_{S_{N}}[\protect\fix(w)]$
and a generalization}

For every $k\in\mathbb{Z}_{\ge1}$ and $N\ge1$ define $\zeta_{k}\colon S_{N}\to\mathbb{Z}_{\ge0}$
by $\zeta_{k}\left(\sigma\right)=\#\mathrm{fix}(\sigma^{k})$ (the
number of fixed points of $\sigma^{k}$). For every partition $\nu=\left(\nu_{1},\ldots,\nu_{\ell}\right)$
define $\zeta_{\nu}\defi\zeta_{\nu_{1}}\zeta_{\nu_{2}}\cdots\zeta_{\nu_{\ell}}$.
First, $\mathbb{E}_{w}\left[\zeta_{\nu}\right]$ is given a geometric
interpretation (\cite[Lem.~6.2]{PP15} and \cite[Prop.~3.7]{hanany2020word}):
Let $\eta_{w}\colon\Gamma_{w}\to\Omega$ be as in Definition \ref{def:Gamma_w}.
The expected value of $\zeta_{1}$ under the $w$-measure in $S_{N}$
is equal to the \emph{average number of lifts of $\eta_{w}$ to a
random degree-$N$ cover of $\Omega$}.\footnote{\label{fn:random N cover}Here, a random degree-N cover of a finite
graph $\Delta$ is in the sense of \cite[P.~9244]{hanany2020word}.
If $\Delta$ is connected, the $N$-covers are in one to one correspondence
with $\Hom(\pi_{1}(\Delta),S_{N})$, and a random $N$-cover corresponds
to a uniformly random homomorphism in this set.} Similarly, as a function on $S_{N}$, $\mathbb{E}_{w}\left[\zeta_{\nu}\right]$
is equal to the average number of lifts of $\eta_{w^{\nu}}\colon\Gamma_{w^{\nu}}\to\Omega$
to a random $N$-cover of $\Omega$.

Based on this observation, we may generalize $\mathbb{E}_{w}\left[\zeta_{\nu}\right]$
and study a similar quantity defined for \emph{any} morphism of finite
graphs. Given such a morphism $\eta\colon\Gamma\to\Delta$, define
$\Phi_{\eta}\left(N\right)$ to be the average number of lifts of
$\eta$ to $\widehat{\Delta_{N}}$ -- a random $N$-cover of $\Delta$.
\[
\Phi_{\eta}\left(N\right)\defi\mathbb{E}\left[\begin{gathered}\xymatrix{ & \widehat{\Delta_{N}}\ar@{->>}[d]^{\rho~~}\\
\Gamma\ar[r]_{\eta}\ar@{-->}[ur]^{\#} & \Delta
}
\end{gathered}
\right]
\]
The quantity $\Phi_{\eta}\left(N\right)$ can also be defined in group-theoretic
terms with respect to the group $\pi_{1}\left(\Delta\right)$ and
the subgroup $\eta_{*}\left(\pi_{1}\left(\Gamma\right)\right)$ \cite[\S2,3]{hanany2020word}.
The analysis continues with immersions of core graphs only.

\subsubsection{\label{subsec:M=0000F6bius-inversions}Möbius inversions}

It turns out to be very useful to decompose $\Phi_{\eta}$ as a sum
of other functions using Möbius inversions. We use two types of decompositions
of $\eta$: one based on \emph{surjective} immersions of core graphs,
and the other on algebraic immersions. 
\begin{defn}
\label{def:decompositions}\cite[Def.~6.4 and 6.13]{hanany2020word}
Let $\eta\colon\Gamma\to\Delta$ be an immersion of core graphs. There
are finitely many decompositions of $\eta$ into triples $\left(\eta_{1},\eta_{2},\eta_{3}\right)$
of immersions such that $\eta=\eta_{3}\circ\eta_{2}\circ\eta_{1}$,
and $\eta_{1}$ and $\eta_{2}$ are surjective (algebraic, respectively).
We denote the set of these decompositions by $\decompt\left(\eta\right)$
($\algdecompt\left(\eta\right)$, respectively).\footnote{In \cite{hanany2020word}, the set $\decompt$ is denoted ${\cal D}\mathrm{ecomp}_{B}^{3}$,
where $B$ is a fixed basis of $\F$. This was supposed to stress
that when the graphs are Stallings (multi) core graphs and represent
(multisets of) subgroups of $\F$, the corresponding morphism may
be surjective with respect to one basis and non-surjective with respect
to another. In contrast, the set of decompositions into algebraic
morphisms does not depend on the basis -- they always represent the
same group-theoretic morphisms in the level of fundamental groups.
In the current paper, we stick to geometric language and the notation
$\decompt$ seems more appropriate.} Similarly, define $\decomp(\eta)$ and $\algdecomp(\eta)$ as the
set of decompositions $\eta=\eta_{2}\circ\eta_{1}$ with $\eta_{1}$
surjective (algebraic, respectively). Two decompositions are considered
equivalent if there are isomorphisms between the intermediate domains
which commute with the decompositions (this is explained in detail
in \cite[Def.~6.4]{hanany2020word}).
\end{defn}

If $\eta$ itself is surjective, then for every $\left(\eta_{1},\eta_{2},\eta_{3}\right)\in\decompt\left(\eta\right)$,
the morphism $\eta_{3}$ is also surjective. If $\eta$ is not surjective,
there is a unique decomposition $\eta=\eta_{\mathrm{inj}}\circ\eta_{\mathrm{sur}}$
with $\eta_{\mathrm{sur}}$ surjective and $\eta_{\mathrm{inj}}$
injective, $\decompt\left(\eta\right)$ is in bijection with $\decompt\left(\eta_{\mathrm{sur}}\right)$
and $\Phi_{\eta}=\Phi_{\eta_{\mathrm{sur}}}$ \cite[Prop.~4.3, Items (1) and (4)]{hanany2020word}.
Likewise, if $\eta$ itself is algebraic, then for every $\left(\eta_{1},\eta_{2},\eta_{3}\right)\in\algdecompt\left(\eta\right)$,
the morphism $\eta_{3}$ is also algebraic. If $\eta$ is not algebraic,
recalling the unique decomposition $\eta=\eta_{\mathrm{free}}\circ\eta_{\mathrm{alg}}$
from \eqref{eq:algebraic-free decomposition}, we have $\Phi_{\eta}=\Phi_{\eta_{\mathrm{alg}}}$
\cite[Prop.~4.3(4)]{hanany2020word}.

The fact that $\decomp\left(\eta\right)$ and $\decompt\left(\eta\right)$
are finite follows from the easy observation that there are only finitely
many surjective morphisms from a given finite graph. Because every
algebraic immersion of core graphs is surjective \cite[Thm.~4.7(1)]{hanany2020word},
$\algdecompt\left(\eta\right)\subseteq\decompt\left(\eta\right)$
for every algebraic $\eta$ and, in particular, $\algdecompt\left(\eta\right)$
is finite too (similarly with $\algdecomp(\eta)$). 

These two types of decompositions allow us to define two families
of Möbius inversions of $\Phi$. We think of $\Phi$ as a map from
morphisms of core graphs to functions $\{\mathbb{Z}_{\ge1}\to\mathbb{Q}\}$
with $\Phi_{\eta}$ denoting the image of $\eta$, and the Möbius
inversions will be maps with similar domains and codomains. First,
there are unique functions, $\LB$ ($L$ for left inversion) defined
on any immersion and $\CB$ ($C$ for central inversion) and $\RB$
($R$ for right inversion), defined on surjective immersions, which
satisfy

\begin{equation}
\Phi_{\eta}=\sum_{\left(\eta_{1},\eta_{2}\right)\in\decomp(\eta)}\LB_{\eta_{2}}~~~~~~~~~~~~~~~~~~\Phi_{\eta}=\sum_{\left(\eta_{1},\eta_{2}\right)\in\decomp\left(\eta\right)}\RB_{\eta_{1}}\label{eq:L-surj and R_surj}
\end{equation}
\begin{equation}
\Phi_{\eta}=\sum_{\left(\eta_{1},\eta_{2},\eta_{3}\right)\in\decompt\left(\eta\right)}\CB_{\eta_{2}},\label{eq:C-surj}
\end{equation}
for every immersion $\eta$ of core graphs. The existence and uniqueness
of $\LB$, $\CB$ and $\RB$ can be seen easily by induction on $\left|\decomp\left(\eta\right)\right|$
or $\left|\decompt\left(\eta\right)\right|$ -- see \cite[\S5]{PP15}
and \cite[\S6]{hanany2020word}.\footnote{There is one subtle point here: two morphisms are considered equivalent
if they differ by post-composition with an isomorphism of the codomain.
In particular, in each of the decompositions of $\Phi_{\eta}$ as
a sum as in \eqref{eq:L-surj and R_surj} or \eqref{eq:C-surj}, there
is a single summand of the form $\LB_{\eta}$, $\RB_{\eta}$ or $\CB_{\eta}$,
and the remaining summands are associated with morphisms $\eta'$
with strictly smaller sets $\decompt(\eta')$ or $\decomp(\eta')$.} Similarly, there are unique functions -- $L^{\alg}$ defined on
any immersions and $C^{\mathrm{alg}}$ and $R^{\alg}$ defined on
algebraic immersions, which satisfy
\begin{equation}
\Phi_{\eta}=\sum_{\left(\eta_{1},\eta_{2}\right)\in\algdecomp(\eta)}L_{\eta_{2}}^{\alg}~~~~~~~~~~~~~~~~~~\Phi_{\eta}=\sum_{\left(\eta_{1},\eta_{2}\right)\in\algdecomp\left(\eta\right)}R_{\eta_{1}}^{\alg}\label{eq:L-alg and R-alg}
\end{equation}
\begin{equation}
\Phi_{\eta}=\sum_{\left(\eta_{1},\eta_{2},\eta_{3}\right)\in\algdecompt\left(\eta\right)}C_{\eta_{2}}^{\alg},\label{eq:C-alg}
\end{equation}
for every immersion $\eta$ of core graphs.

\subsubsection{Useful properties of the Möbius inversions\label{subsec:Useful-properties-of-Mobius}}

The following properties of the Möbius inversions defined above will
be useful in subsequent sections. The most involved part of the analysis
of $\mathbb{E}_{w}\left[\zeta_{\nu}\right]=\Phi_{\eta_{w^{\nu}}}$
in \cite{PP15,hanany2020word} deals with the functions $\CB$ and
$C^{\mathrm{alg}}$. We state here the main result only for $C^{\mathrm{alg}}$:
\begin{prop}
\label{prop:bound on C^alg}\cite[Cor.~6.16]{hanany2020word} Let
$\eta\colon\Gamma\to\Delta$ be an algebraic immersion of core graphs.
Then 
\begin{equation}
C_{\eta}^{\mathrm{alg}}\left(N\right)=\begin{cases}
N^{\chi\left(\Gamma\right)}, & \mathrm{if}~\eta=\id~\mathrm{is~an~isomorphism,}\\
O\left(N^{\min\left(\chi\left(\Gamma\right),\chi\left(\Delta\right)\right)-1}\right) & \mathrm{otherwise.}
\end{cases}\label{eq:bounds on C^alg}
\end{equation}
\end{prop}

\begin{rem}
The bound \eqref{eq:bounds on C^alg} is enough to get the main result
of \cite{PP15}. Indeed, recall that \cite{PP15} analyses $\mathbb{E}_{w}\left[\zeta_{1}\right]=\Phi_{\eta_{w}}\left(N\right)$.
By definition, $\pi\left(w\right)$ is the smallest rank of a proper
algebraic extension of $\left\langle w\right\rangle $, namely, $1-\pi\left(w\right)$
is the largest Euler characteristic of (the Stallings core graph corresponding
to) a proper algebraic extension of $\left\langle w\right\rangle $.
In every decomposition $\Gamma_{w}\stackrel{\eta_{1}}{\to}\Sigma\stackrel{\eta_{2}}{\to}\Sigma'\stackrel{\eta_{3}}{\to}\Omega$
in $\algdecompt\left(\eta_{w}\right)$, the morphism $\eta_{2}\colon\Sigma\to\Sigma'$
is between two (core graphs representing) algebraic extensions of
$\left\langle w\right\rangle $. By \eqref{eq:bounds on C^alg}, $C_{\eta_{2}}^{\mathrm{alg}}\left(N\right)=O\left(N^{-\pi\left(w\right)}\right)$,
except for in the following two cases:
\begin{enumerate}
\item $\left(\eta_{1},\eta_{2},\eta_{3}\right)=\left(\id,\id,\eta_{w}\right)$,
in which case $C_{\eta_{2}}^{\mathrm{alg}}\left(N\right)=1$, and
\item $\Sigma=\Sigma'$ represents a proper algebraic extension of $w$
of rank $\pi\left(w\right)$, in which case $C_{\eta_{2}}^{\mathrm{alg}}\left(N\right)=N^{1-\pi\left(w\right)}$.
\end{enumerate}
Overall, we get \cite[Thm.~1.8]{PP15}:
\[
\mathbb{E}_{w}\left[\zeta_{1}\right]=\sum_{\left(\eta_{1},\eta_{2},\eta_{3}\right)\in\algdecompt\left(\eta_{w}\right)}C_{\eta_{2}}^{\mathrm{alg}}\left(N\right)=1+c_{w}N^{1-\pi\left(w\right)}+O\left(N^{-\pi\left(w\right)}\right),
\]
where $c_{w}\in\mathbb{Z}_{\ge1}$ is the number of proper algebraic
extensions of $\left\langle w\right\rangle $ of rank $\pi\left(w\right)$.
\end{rem}

The following Lemma is part of \cite[Lem.~6.17]{hanany2020word} and
strengthens Proposition \ref{prop:bound on C^alg} when $\chi(\Delta)=0$.
If $\eta\colon\Gamma\to\Delta$ is a surjective immersion of core
graphs such that $\Delta$ is a cycle, then $\eta$ is a covering
map with degree $\deg\left(\eta\right)$. Note that a core graph $\Gamma$
is a union of cycles if and only if $\chi\left(\Gamma\right)=0$. 
\begin{lem}
\label{lem:C^alg vanishes on non-isomorphisms of cycles}Let $\eta\colon\Gamma\to\Delta$
be an algebraic immersion of core graphs which are both unions of
cycles, and let $N_{0}$ be the maximal degree of $\eta$ restricted
to a component of $\Delta$. Then for every $N\ge N_{0}$,
\begin{equation}
C_{\eta}^{\mathrm{alg}}\left(N\right)=\begin{cases}
1 & \mathrm{if}~\eta~\mathrm{is~an~isomorphism},\\
0 & \mathrm{otherwise.}
\end{cases}\label{eq:C^alg on cycle graphs}
\end{equation}
\end{lem}

It turns out that $\LB$ has, too, a clean geometric interpretation:
whereas $\Phi$ counts lifts, $\LB$ counts \emph{injective} lifts.
\begin{lem}
\label{lem:L^B =00003D average number of inj lifts}\cite[Prop.~6.8]{hanany2020word}
Given an immersion $\eta\colon\Gamma\to\Delta$ of core graphs, $\LB_{\eta}(N)$
is equal to the average number of \emph{\uline{injective}} lifts
of $\eta$ to a random $N$-cover of $\Delta$. It also satisfies
\[
\LB_{\eta}(N)=N^{\chi(\Gamma)}\left(1+O\left(N^{-1}\right)\right).
\]
\end{lem}

We also need the following simple lemmas.
\begin{lem}
\label{lem:All Mobius inversions are  multiplicative on connected components of codomain}
All seven maps $\Phi,\LB,\CB,\RB,L^{\alg},$$C^{\mathrm{alg}}$ and
$R^{\alg}$ are multiplicative on connected components of the codomain.
\end{lem}

Namely, if, for example, $\eta\colon\Gamma\to\Delta$ is algebraic,
$\Delta_{1},\ldots,\Delta_{k}$ are the connected components of $\Delta$
and $\eta_{i}\colon\eta^{-1}\left(\Delta_{i}\right)\to\Delta_{i}$
are the associated restrictions of $\eta$, then $C_{\eta}^{\mathrm{alg}}=\prod_{i=1}^{k}C_{\eta_{i}}^{\mathrm{alg}}$.
\begin{proof}
First, $\Phi$ has this property: a random $N$-cover of $\Delta$
is an array of \emph{independent} random $N$-covers of each of $\Delta_{1},\ldots,\Delta_{k}$,
so the number of lifts of $\eta_{i}$ is independent of that of $\eta_{j}$
if $i\ne j$. This property of $\Phi$ is inherited by all its Möbius
inversions. We illustrate the proof with $C^{\alg}$, the other cases
being very similar. 

It is enough to prove that whenever $\eta\colon\Gamma\to\Delta$ is
an algebraic immersion of core graphs that satisfies $\eta=\eta_{1}\sqcup\eta_{2}$
(both non-empty morphisms), then $C_{\eta}^{\alg}=C_{\eta_{1}}^{\alg}\cdot C_{\eta_{2}}^{\alg}$.
We proceed by induction on $|\algdecompt(\eta)|$ and assume the property
holds for algebraic $\eta'$ whenever $|\algdecompt(\eta')|<|\algdecompt(\eta)|$.
Note that 
\[
\algdecompt\left(\eta\right)=\algdecompt\left(\eta_{1}\right)\times\algdecompt\left(\eta_{2}\right).
\]
Thus,
\begin{eqnarray}
\Phi_{\eta} & = & \Phi_{\eta_{1}}\Phi_{\eta_{2}}=\left[\sum_{\substack{\left(\gamma_{1},\gamma_{2},\gamma_{3}\right)\in\algdecompt\left(\eta_{1}\right)\colon}
}C_{\gamma_{2}}^{\alg}\right]\cdot\left[\sum_{\left(\delta_{1},\delta_{2},\delta_{3}\right)\in\algdecompt\left(\eta_{2}\right)}C_{\delta_{2}}^{\alg}\right]\nonumber \\
 & = & \left[\sum_{\substack{\left(\varepsilon_{1},\varepsilon_{2},\varepsilon_{3}\right)\in\algdecompt\left(\eta\right)\\
\left(\varepsilon_{1},\varepsilon_{2},\varepsilon_{3}\right)\ne\left(\id,\eta,\id\right)
}
}C_{\varepsilon_{2}}^{\alg}\right]+C_{\eta_{1}}^{\alg}C_{\eta_{2}}^{\alg},\label{eq:proof that Mobius inversion are multiplicative on conn comps}
\end{eqnarray}
where the last equality is by the induction hypothesis: whenever $(\gamma_{2},\delta_{2})\ne(\eta_{1},\eta_{2})$
we have 
\[
\left|\algdecompt\left(\gamma_{2}\sqcup\delta_{2}\right)\right|<\left|\algdecompt\left(\eta\right)\right|,
\]
so $C_{\gamma_{2}\sqcup\delta_{2}}^{\alg}=C_{\gamma_{2}}^{\alg}C_{\delta_{2}}^{\alg}$.
By the definition of $C^{\alg}$ in \eqref{eq:C-alg}, we conclude
that $C_{\eta}^{\alg}=C_{\eta_{1}}^{\alg}\cdot C_{\eta_{2}}^{\alg}$. 
\end{proof}
\begin{lem}
\label{lem:L^alg in terms of L}Let $\eta$ be an immersion of core
graphs. Then 
\begin{equation}
L_{\eta}^{\mathrm{alg}}=\sum_{\substack{\substack{\left(\eta_{1},\eta_{2}\right)\in\decomp(\eta)\colon}
\\
\eta_{1}\,\mathrm{is\,free}
}
}\LB_{\eta_{2}}.\label{eq:L^alg in term of L^B}
\end{equation}
\end{lem}

\begin{proof}
The proof is almost identical to that of \cite[Prop.~6.15]{hanany2020word},
which gives a similar statement for $C^{\mathrm{alg}}$. We repeat
it here because there is a minor error in that proof in \cite{hanany2020word}
(although not in the statement). Denote the right hand side of \eqref{eq:L^alg in term of L^B}
by $F_{\eta}\in\mathbb{Q}^{\mathbb{Z}_{\ge1}}$. By the unique decomposition
$\eta=\eta_{\mathrm{free}}\circ\eta_{\mathrm{alg}}$ from \eqref{eq:algebraic-free decomposition}
we have
\begin{eqnarray*}
\Phi_{\eta} & = & \sum_{\substack{\left(\eta_{1},\eta_{2}\right)\in\decomp(\eta)}
}\LB_{\eta_{2}}=\sum_{\substack{\left(\alpha,\beta,\eta_{2}\right)\colon\eta_{2}\circ\beta\circ\alpha=\eta,\\
\alpha~\mathrm{algebraic},~\beta~\mathrm{surjective~and~free}
}
}\LB_{\eta_{2}}\\
 & = & \sum_{\substack{\left(\alpha,\gamma\right)\in\algdecomp(\eta)}
}F_{\gamma}.
\end{eqnarray*}
By the definition \eqref{eq:L-alg and R-alg} and uniqueness of $L^{\mathrm{alg}}$,
we must have that $L_{\gamma}^{\mathrm{alg}}=F_{\gamma}$.
\end{proof}
\begin{cor}
\label{cor:L-alg order of chi}Given an immersion $\eta\colon\Gamma\to\Delta$
of core graphs, $L_{\eta}^{\alg}(N)=N^{\chi(\Gamma)}\left(1+O\left(N^{-1}\right)\right).$
\end{cor}

\begin{proof}
By Lemmas \ref{lem:L^B =00003D average number of inj lifts} and \ref{lem:L^alg in terms of L},
\begin{equation}
L_{\eta}^{\alg}(N)=\sum_{\substack{\Gamma\stackrel{\eta_{1}}{\twoheadrightarrow}\Sigma\stackrel{\eta_{2}}{\longrightarrow}\Delta\colon\,\eta=\eta_{2}\circ\eta_{1}\\
\eta_{1}~\mathrm{free~and~surjective}
}
}\LB_{\eta_{2}}\left(N\right)=\sum_{\substack{\Gamma\stackrel{\eta_{1}}{\twoheadrightarrow}\Sigma\stackrel{\eta_{2}}{\longrightarrow}\Delta\colon\,\eta=\eta_{2}\circ\eta_{1}\\
\eta_{1}~\mathrm{free~and~surjective}
}
}N^{\chi\left(\Sigma\right)}\left(1+O\left(N^{-1}\right)\right).\label{eq:estimating L-alg}
\end{equation}
By \cite[Prop.~4.3(2)]{hanany2020word}, whenever $\eta_{1}\colon\Gamma\to\Sigma$
is a surjective free morphism, $\chi(\Sigma)\le\chi(\Gamma)$, with
equality if and only if $\eta_{1}$ is an isomorphism. Thus, there
is a single summand of order $N^{\chi(\Gamma)}$ in \eqref{eq:estimating L-alg}:
the one corresponding to $(\eta_{1},\eta_{2})=(\id,\eta)$, and the
remaining ones are of strictly smaller order. 
\end{proof}

\subsection{The groups $G\wr S_{N}$ and their stable irreducible representations\label{subsec:reps of G wr S_N}}

Fix a compact group $G$. Throughout this section, we denote $G_{N}\defi G\wr S_{N}$.
The elements of $G_{N}$ are $\left(v,\sigma\right)\in G^{N}\times S_{N}$,
$v=(v_{1},\ldots,v_{N})$ with $\left(v,\sigma\right)\cdot\left(u,\tau\right)=((v_{1}u_{\sigma^{-1}\left(1\right)},\ldots,v_{N}u_{\sigma^{-1}\left(N\right)}),\sigma\tau)$.
Equivalently, the elements of $G_{N}$ are $N\times N$ monomial (also
known as generalized permutation) matrices with entries in $G\sqcup\left\{ 0\right\} $.
The element $\left(v,\sigma\right)$ corresponds to the monomial matrix
$A$ with $A_{i,\sigma^{-1}\left(i\right)}=v_{i}$, and the multiplication
is an ordinary matrix multiplication. 

Denote by ${\cal C}\left(G\right)$ the set of conjugacy classes of
$G$. The conjugacy class of $\left(v,\sigma\right)\in G_{N}$ is
determined by the numbers $(a_{\ell,C})_{\ell\in\left[N\right],C\in{\cal C}\left(G\right)}$,
where $a_{\ell,C}$ is the number of $\ell$-cycles in $\sigma$ with
product of the corresponding $v_{i}$'s belonging to the conjugacy
class $C$. Here, if the $\ell$-cycle in $\sigma$ is $\left(x_{1}~\ldots~x_{\ell}\right)$,
the product is $v_{x_{\ell}}v_{x_{\ell-1}}\cdots v_{x_{1}}$ \cite[\S2.3]{ceccherini2014representation}. 

The irreducible characters of $G_{N}$ are parameterized by functions
$\arrm\colon\irr(G)\to{\cal P}$ from the set of irreducible characters
of $G$ to the set of all partitions, with $\left|\arrm\right|\defi\sum_{\phi\in\irr(G)}\left|\arrm\left(\phi\right)\right|=N$
(see the detailed \cite[\S2.6]{ceccherini2014representation} when
$G$ is finite, or \cite[\S1.2]{hora2008limits} in the general case).
The support of $\arrm$ is $\supp(\arrm)\defi\left\{ \phi\in\irr(G)\,|\,\arrm(\phi)\ne\emptyset\right\} $.
We denote the irreducible character corresponding to $\arrm$ by $\chi^{\arrm}$.
In the following proposition, an \emph{ordered set partition} is a
set partition (of some set) with a prescribed order on the blocks
of the partition. For a partition $\mu\in{\cal P}$, we denote by
$\chi^{\mu}\in\irr(S_{|\mu|})$ the irreducible character of $S_{|\mu|}$
corresponding to the partition $\mu$.
\begin{prop}
\label{prop:irred chars of G wr S_N}
\begin{enumerate}
\item \label{enu:formula when support of size 1}Let $\mu\vdash N\ge1$
be a partition and $\phi\in\irr\left(G\right)$ an irreducible character
of $G$. Let $\arrm\colon\irr(G)\to{\cal P}$ be supported on $\phi$
with $\arrm(\phi)=\mu$. The corresponding irreducible character $\chi^{\phi,\mu}\defi\chi^{\arrm}$
of $G_{N}$ is given by
\[
\chi^{\phi,\mu}\left(v,\sigma\right)=\chi^{\mu}\left(\sigma\right)\prod_{\left(x_{1}~\ldots~x_{\ell}\right):~\mathrm{cycle~of~\sigma}}\phi\left(v_{x_{\ell}}\cdots v_{x_{1}}\right),
\]
where $\chi^{\mu}\in\irr\left(S_{N}\right)$, and the product goes
over all the cycles of $\sigma$.
\item \label{enu:formula for character as induction}Let $\arrm\colon\irr(G)\to{\cal P}$
with $|\arrm|=N$ and $\supp(\arrm)=\{\phi_{1},\ldots,\phi_{k}\}$.
Denote $d_{i}=|\arrm\left(\phi_{i}\right)|$, so that $d_{i}\ge1$
and $\sum d_{i}=N$. The irreducible character corresponding to $\arrm$
is given by 
\[
\chi^{\arrm}=\mathrm{Ind_{G_{d_{1}}\times\ldots\times G_{d_{k}}}^{G_{N}}\left(\chi^{\phi_{1},\arrm\left(\phi_{1}\right)}\boxtimes\ldots\boxtimes\chi^{\phi_{k},\arrm\left(\phi_{k}\right)}\right).}
\]
\item \label{enu:formula for value of induced character}Keeping the notation
from Item \ref{enu:formula for character as induction}, and fixing
$\left(v,\sigma\right)\in G_{N}$, let ${\cal B}_{\sigma}$ denote
the set of ordered set-partitions $\overline{B}=(B_{1},\ldots,B_{k})$
of $[N]$ into $k$ subsets $B_{1}\sqcup\ldots\sqcup B_{k}=[N]$ with
$\left|B_{i}\right|=d_{i}$ which are invariant under $\sigma\in S_{N}$,
namely, such that $\sigma(B_{i})=B_{i}$ for all $i$. Then
\begin{equation}
\chi^{\arrm}\left(v,\sigma\right)=\sum_{\overline{B}\in{\cal B}_{\sigma}}\prod_{i=1}^{k}\chi^{\arrm(\phi_{i})}\left(\sigma|_{B_{i}}\right)\prod_{\left(x_{1}~\ldots~x_{\ell}\right):~\mathrm{cycle~of~\sigma|_{B_{i}}}}\phi_{i}\left(v_{x_{\ell}}\cdots v_{x_{1}}\right).\label{eq:formula for induced char}
\end{equation}
\end{enumerate}
\end{prop}

\begin{proof}
The first two items are a combination of Theorems 2.4.4, 2.4.5 and
2.6.1 in \cite{ceccherini2014representation}. We now prove the third
one. Let $\overline{Q}=(Q_{1},\ldots,Q_{k})$ denote the ordered set-partition
of $[N]$ where $Q_{1}=[d_{1}]$, $Q_{2}=\left\{ d_{1}+1,\ldots,d_{1}+d_{2}\right\} $
and so on. For $y=\left(y_{1},\ldots,y_{k}\right)\in G_{d_{1}}\times\ldots\times G_{d_{k}}$,
denote
\[
\xi^{\arrm}\left(y\right)\defi\chi^{\phi_{1},\arrm\left(\phi_{1}\right)}\left(y_{1}\right)\cdots\chi^{\phi_{k},\arrm\left(\phi_{k}\right)}\left(y_{k}\right).
\]
By Item \ref{enu:formula for character as induction} and the definition
of induction, it $T\subseteq G_{N}$ is a set of representatives of
the right cosets of $G_{d_{1}}\times\ldots\times G_{d_{k}}$ in $G_{N}$,
then
\[
\chi^{\arrm}\left(v,\sigma\right)=\sum_{t\in T\colon t\left(v,\sigma\right)t^{-1}\in G_{d_{1}}\times\ldots\times G_{d_{k}}}\xi^{\arrm}\left(t\left(v,\sigma\right)t^{-1}\right).
\]
It is clear that one can have all representatives of the form $t=\left(1,\tau\right)$,
belonging to the subgroup isomorphic to $S_{N}$. In this case, $(1,\tau)(v,\sigma)(1,\tau)^{-1}\in G_{d_{1}}\times\ldots\times G_{d_{k}}$
if and only if $\tau\sigma\tau^{-1}\in S_{d_{1}}\times\ldots\times S_{d_{k}}$,
if and only if $\tau^{-1}.\overline{Q}\defi(\tau^{-1}.Q_{1},\ldots,\tau^{-1}.Q_{k})\in{\cal B}_{\sigma}$.
Moreover, the set $\{\tau^{-1}.\overline{Q}\}_{(1,\tau)\in T}$ is
precisely the set of all ordered set-partitions of $[N]$ to $k$
subsets of sizes $d_{1},\ldots,d_{k}$. Finally, whenever $\tau^{-1}.\overline{Q}\in{\cal B}_{\sigma}$,
we have $(1,\tau)(v,\sigma)(1,\tau)^{-1}=((v_{\tau^{-1}(1)},\ldots,v_{\tau^{-1}(N)}),\tau\sigma\tau^{-1})$,
the permutation $\sigma|_{\tau^{-1}.Q_{i}}$ is conjugate to $\tau\sigma\tau^{-1}|_{Q_{i}}$
and 
\[
\prod_{\left(x_{1}~\ldots~x_{\ell}\right):~\mathrm{cycle~of~\sigma|_{\tau^{-1}.Q_{i}}}}\phi_{i}\left(v_{x_{\ell}}\cdots v_{x_{1}}\right)=\prod_{\left(x_{1}~\ldots~x_{\ell}\right):~\mathrm{cycle~of~\tau\sigma\tau^{-1}|_{Q_{i}}}}\phi_{i}\left(v_{\tau^{-1}\left(x_{\ell}\right)}\cdots v_{\tau^{-1}\left(x_{1}\right)}\right).
\]
\end{proof}
For any partition $\mu$ and $N\ge|\mu|+\mu_{1}$ (if $\mu=\emptyset$
then $\mu_{1}=0$), denote 
\begin{equation}
\mu[N]\defi\left(N-|\mu|,\mu_{1},\ldots,\mu_{\ell(\mu)}\right)\label{eq:mu=00005BN=00005D}
\end{equation}
the partition of $N$ where outside the first row, its Young diagram
is precisely $\mu$.

\emph{Stable} irreducible characters of $G\wr S_{\bullet}$ are parameterized
by partition-valued functions $\arrm\colon\irr(G)\to{\cal P}$ with
$\left|\arrm\right|<\infty$ \cite{sam2019representations}. For every
$N\ge\left|\arrm\right|+\arrm\left(\mathrm{triv}\right)_{1}$ (again,
if $\arrm(\triv)=\emptyset$ define $\arrm(\triv)_{1}=0$), we define
$\arrm\left[N\right]$ as the function which agrees with $\arrm$
on all non-trivial characters and is $\arrm\left(\mathrm{triv}\right)[N']$
on the trivial representation, where $N'=N-|\arrm|+|\arrm(\triv)|$.
The stable irreducible character corresponding to $\arrm$ is 
\[
\chi^{\arrm[\bullet]}\defi\left\{ \chi^{\arrm\left[N\right]}\right\} {}_{N\ge\left|\arrm\right|+\arrm\left(\mathrm{triv}\right)_{1}}.
\]

\begin{prop}
\label{prop:dim of chi^arrm}Let $\arrm\colon\irr(G)\to{\cal P}$
with $\left|\arrm\right|=d<\infty$. Then
\[
\dim\left(\chi^{\arrm\left[N\right]}\right)=\Theta\left(N^{d}\right).
\]
\end{prop}

\begin{proof}
Use the notation from Proposition \ref{prop:irred chars of G wr S_N}\eqref{enu:formula for character as induction}
with $\arrm[N]$, assuming that $\phi_{1}=\triv$. By Item \ref{enu:formula when support of size 1}
of that proposition, for $i=2,\ldots,k$, the dimension of $\chi^{\phi_{i},\arrm[N](\phi_{i})}$
is 
\[
\chi^{\phi_{i},\arrm[N](\phi_{i})}\left(1_{G_{d_{i}}}\right)=\chi^{\arrm\left(\phi_{i}\right)}\left(1_{S_{d_{i}}}\right)\left(\dim\phi_{i}\right)^{d_{i}},
\]
which is a constant. For $i=1$ we have 
\[
\chi^{\arrm[N](\triv)}\left(1_{S_{d_{1}}}\right)=\chi^{\arrm(\triv)[N-d_{2}-\ldots-d_{k}]}\left(1_{S_{d_{1}}}\right)=\Theta\left(\left(N-d_{2}-\ldots-d_{k}\right)^{\left|\arrm(\triv)\right|}\right)=\Theta\left(N^{\left|\arrm(\triv)\right|}\right),
\]
where the middle equality is by $\dim\chi^{\mu[N]}=\Theta(N^{|\mu|})$,
which follows immediately from the hook length formula. Hence the
dimension of the exterior tensor product of these representations
on $G_{d_{1}}\times\ldots\times G_{d_{k}}$ is $\Theta(N^{|\arrm(\triv)|})$,
and by Proposition \ref{prop:irred chars of G wr S_N}\eqref{enu:formula for character as induction}
, the dimension of $\chi^{\arrm[N]}$ is of order
\[
\left[G_{N}\colon G_{d_{1}}\times\ldots\times G_{d_{k}}\right]\cdot N^{\left|\arrm(\triv)\right|}=\frac{N!}{(N-d_{2}-\ldots-d_{k})!d_{2}!\cdots d_{k}!}N^{\left|\arrm(\triv)\right|}=\Theta\left(N^{d}\right).
\]
\end{proof}
When $G$ is the trivial group and $G_{N}=S_{N}$, the entire theory
reduces to the much better-known theory of representation of the symmetric
group. In particular, stable representations are given by a single
partition $\mu=\arrm(\triv)$, and we denote
\[
\chi^{\mu[\bullet]}\defi\chi^{\arrm[\bullet]}=\left\{ \chi^{\mu[N]}\right\} _{N\ge|\mu|+\mu_{1}}.
\]

\subsection{The expectations $\mathbb{E}_{w}[\chi]$ and the parameters $\beta(w,\chi)$\label{subsec:The-expectations-and-beta}}

In the case of $S_{\bullet}$, for every stable irreducible character
$\chi^{\mu[\bullet]}=(\chi^{\mu[N]})_{N\ge|\mu|+\mu_{1}}$, it follows
already from \cite{hanany2020word} that $\mathbb{E}_{w}[\chi^{\mu[N]}]$,
the expected value of $\chi^{\mu[N]}$ under the $w$-measure in $S_{N}$,
coincides with a rational function $f_{w,\mu}(N)\in\mathbb{Q}(N)$
for every large enough $N$. As $\dim(\chi^{\mu[N]})$ is given by
a polynomial $g_{\mu}\in\mathbb{Q}[N]$ in $N$ of degree $|\mu|$,
we deduce that $\beta(w,\chi^{\mu[\bullet]})=-\frac{\deg(f_{w,\mu})}{|\mu|}$
is a rational number or $\infty$ (see \cite[\S\S4.2]{PSh23} for
a more detailed explanation).

Theorem \ref{thm:formula for stable irreps of G wr S_bullet} gives
a formula for every stable irreducible character $\chi=(\chi_{N})_{N\ge N_{0}}$
of $G\wr S_{\bullet}$ for an arbitrary compact group $G$. In Corollary
\ref{cor:E_w=00005Bchi=00005D  rational in N and beta rational} we
deduce from this formula similar conclusions, namely, that $\mathbb{E}_{w}[\chi_{N}]$
coincides with a rational function in $\mathbb{C}(N)$ for every large
enough $N$, and that $\beta(w,\chi)$ is either rational or $\infty$.

However, we do not need these results in order to \emph{define }$\beta(w,\chi)$.
Indeed, we may simply define

\begin{equation}
\beta\left(w,\chi\right)\defi\liminf_{N\to\infty}\frac{-\log\left|\mathbb{E}_{w}\left[\chi_{N}\right]\right|}{\log\left(\dim\chi_{N}\right)}.\label{eq:formal-def-of-beta}
\end{equation}
For the sake of this paper, we do not need the fact that $\beta(w,\chi)\in\mathbb{Q}\cup\{\infty\}$,
nor the fact that $\mathbb{E}_{w}[\chi_{N}]$ coincides with a rational
function.

\section{The average value of $(\protect\fix(w_{1})-1)\cdot\cdot\cdot(\protect\fix(w_{k})-1)$\label{sec:generalized PP15}}

The current section proves Theorem \ref{thm:generalized PP15}, which
generalizes the main result of \cite{PP15}. Fix a finite multiset
of non-trivial words $\left\{ w_{1},\ldots,w_{k}\right\} \subset\F\setminus\{1_{\F}\}$.
Denote by $\Gamma_{w_{1},\ldots,w_{k}}\defi\Gamma_{\{w_{1},\ldots,w_{k}\}}\defi\Gamma_{w_{1}}\sqcup\ldots\sqcup\Gamma_{w_{k}}$
the core graph obtained as a disjoint union of the cycles $\Gamma_{w_{1}},\ldots,\Gamma_{w_{k}}$
and equipped with an immersion to the bouquet 
\[
\eta_{w_{1},\ldots,w_{k}}\defi\eta_{\{w_{1},\ldots,w_{k}\}}\defi\eta_{w_{1}}\sqcup\ldots\sqcup\eta_{w_{k}}\colon\Gamma_{w_{1},\ldots,w_{k}}\to\Omega.
\]
 
\begin{defn}
\label{def:proper algebraic extension}A morphism of core graph $\gamma\colon\Gamma\to\Delta$
is called a \textbf{proper algebraic morphism} if it is algebraic
and when restricted to any connected component $\Delta$, it is \emph{not}
an isomorphism. Define\marginpar{$\protect\chialg$} 
\begin{equation}
\chialg\left(w_{1},\ldots,w_{k}\right)\defi\max\left\{ \chi\left(\Delta\right)\,\mid\,\substack{\Gamma_{w_{1},\ldots,w_{k}}\stackrel{\gamma}{\longrightarrow}\Delta\stackrel{\delta}{\longrightarrow}\Omega\,\,\mathrm{is\,a\,decomposition\,of}\,\eta_{w_{1},\ldots,w_{k}},\\
\mathrm{and}\,\gamma\,\mathrm{is\,a\,proper\,algebraic\,morphism}
}
\right\} ,\label{eq:def of chi-prop-alg}
\end{equation}
where the maximum over the empty set is $-\infty$. Finally, denote
by $\crit(w_{1},\ldots,w_{k})$\marginpar{${\scriptscriptstyle \protect\crit(w_{1},\ldots,w_{k})}$}
the set of ``critical'' proper algebraic extensions of $\Gamma_{w_{1},\ldots,w_{k}}$,
namely, the set consisting of all proper algebraic morphisms $\gamma\colon\Gamma_{w_{1},\ldots,w_{k}}\to\Delta$
with $\chi(\Delta)=\chialg(w_{1,}\ldots,w_{k})$ such that $\eta_{w_{1},\ldots,w_{k}}$
factors through $\gamma$.
\end{defn}

(As usual, the restriction of $\gamma$ to a connected component $\Delta'$
of $\Delta$ is $\gamma|_{\gamma{}^{-1}(\Delta')}\colon\gamma^{-1}(\Delta')\to\Delta'$.)
We remark that although not reflected in the notation, both $\chialg(w_{1},\ldots,w_{k})$
and $\crit(w_{1},\ldots,w_{k})$ depend only on the multiset $\left\{ w_{1},\ldots,w_{k}\right\} $:
there is no significance to the order of the words. Note that $\crit(w_{1},\ldots,w_{k})$
can be thought of as a subset of the finite set $\algdecomp(\eta_{w_{1},\ldots,w_{k}})$.
In particular, the set $\crit(w_{1},\ldots,w_{k})$ is, too, finite.

Recall the primitivity rank $\pi(w)$ from \eqref{eq:def of pi}.
For a single non-trivial word $1\ne w\in\F$, $\chialg(w)=1-\pi(w)$,
and \cite[Thm.~1.8]{PP15}, stated in \eqref{eq:PP15} above, translates
to that 
\[
\mathbb{E}_{S_{N}}\left[\mathrm{std}_{N}(w)\right]=\mathbb{E}_{S_{N}}\left[\fix(w)-1\right]=\left|\crit(w)\right|\cdot N^{\chialg(w)}\left(1+O\left(N^{-1}\right)\right).
\]
Theorem \ref{thm:generalized PP15} generalizes this result to a multiset
of words, and states that for $w_{1},\ldots,w_{k}\in\F\setminus\{1_{\F}\}$,
\begin{equation}
\mathbb{E}_{S_{N}}\left[\prod_{i=1}^{k}\left(\fix\left(w_{k}\right)-1\right)\right]=N^{\chialg\left(w_{1},\ldots,w_{k}\right)}\left(\left|\crit\left(w_{1},\ldots,w_{k}\right)\right|+O\left(\frac{1}{N}\right)\right).\label{eq:PP15 generalized}
\end{equation}

\begin{proof}[Proof of Theorem \ref{thm:generalized PP15}]
Denote the left hand side of \eqref{eq:PP15 generalized} by $E$.
Note that for any multiset of non-trivial words $U=\{u_{1},\ldots,u_{t}\}\in\F$
we have $\mathbb{E}_{S_{N}}[\prod_{j=1}^{t}\fix(u_{j})]=\Phi_{\eta_{U}}(N)$.
Hence,
\begin{eqnarray}
E & = & \sum_{S\subseteq\left[k\right]}\left(-1\right)^{k-|S|}\mathbb{E}_{S_{N}}\left[\prod_{s\in S}\fix\left(w_{s}\right)\right]=\sum_{S\subseteq\left[k\right]}\left(-1\right)^{k-|S|}\Phi_{\eta_{S}}\left(N\right)\nonumber \\
 & = & \sum_{S\subseteq\left[k\right]}\left(-1\right)^{k-|S|}\sum_{\left(\gamma_{1},\gamma_{2}\right)\in\algdecomp\left(\eta_{S}\right)}R_{\gamma_{1}}^{\mathrm{alg}}(N).\label{eq:equivalent exp}
\end{eqnarray}
If $\eta$ is an isomorphism of some core graph $\Sigma$, then $R_{\eta}(N)=\Phi_{\eta}(N)$
by definition and $\Phi_{\eta}(N)=N^{\chi(\Sigma)}$ (see \cite[Example 3.5]{hanany2020word}).
By Lemma \ref{lem:All Mobius inversions are  multiplicative on connected components of codomain},
$R^{\alg}$ is multiplicative on connected components of the codomain.
In particular, if $\gamma_{1}\colon\Gamma_{S}\to\Sigma$ is an isomorphism
on one of the components $\Sigma_{o}$ of $\Sigma$, then $\Sigma_{o}$
must be a cycle, hence $\chi(\Sigma_{o})=0$. So the value of $R_{\gamma_{1}}$
is the same as $R_{\gamma_{1}'}$ where $\gamma_{1}'$ is obtained
by removing $\Sigma_{o}$ from $\Sigma$ and its preimage from the
domain. We can also extend every $\gamma_{1}'$ defined on a sub-multiset
of the words to the entire multiset by isomorphisms on the remaining
words. 

Consider the equivalent classes of the different morphisms in \eqref{eq:equivalent exp}
up to extending and removing components with isomorphisms. We get
that every algebraic $\gamma_{1}\colon\Gamma_{w_{1},\ldots,w_{k}}\to\Sigma$
which is an isomorphism on at least one component satisfies that the
overall contribution of its equivalence class in \eqref{eq:equivalent exp}
vanishes. We are left precisely with morphisms that are proper algebraic
extensions:
\begin{equation}
E=\sum_{\substack{\left(\gamma_{1},\gamma_{2}\right)\in\algdecomp\left(\eta_{w_{1},\ldots,w_{k}}\right):\\
\gamma_{1}~\mathrm{is~proper~algebraic}
}
}R_{\gamma_{1}}^{\alg}(N)=\sum_{\substack{\left(\delta_{1},\delta_{2},\gamma_{2}\right)\in\algdecompt\left(\eta_{w_{1},\ldots,w_{k}}\right):\\
\delta_{2}\circ\delta_{1}~\mathrm{is~proper~algebraic}
}
}C_{\delta_{2}}^{\alg}(N),\label{eq:E as a sum over proper algebraic morphisms}
\end{equation}
where the last equality is by the definitions of $R^{\alg}$ and $C^{\alg}$.

Now consider a decomposition $\Gamma_{w_{1},\ldots,w_{k}}\stackrel{\delta_{1}}{\longrightarrow}\Sigma_{1}\stackrel{\delta_{2}}{\longrightarrow}\Sigma_{2}\stackrel{\gamma_{2}}{\longrightarrow}\Omega$
from $\algdecompt(\eta_{w_{1},\ldots,w_{k}})$ with $\delta_{2}\circ\delta_{1}$
proper algebraic. By definition, $\chi(\Sigma_{2})\le\chialg(w_{1},\ldots,w_{k})$.
By Proposition \ref{prop:bound on C^alg}, $C_{\delta_{2}}^{\alg}(N)=O(N^{\chi(\Sigma_{2})-1})=O(N^{\chialg(w_{1},\ldots,w_{k})-1})$,
unless $\delta_{2}$ is an isomorphism and $\delta_{1}\in\crit(w_{1},\ldots,w_{k})$.
In the latter case, $C_{\delta_{2}}^{\alg}(N)=N^{\chialg(w_{1},\ldots,w_{k})}.$
This completes the proof of Theorem \ref{thm:generalized PP15}. 
\end{proof}
We end this section with some observations regarding the possible
values of $\chialg(w_{1},\ldots,w_{k})$. Recall that the Euler characteristic
of a f.g.~free group $\F$ is defined to be the Euler characteristic
of its classifying space, namely, $\chi(\F)=\chi(\Omega)=1-\rk(\F)$.
\begin{prop}
\label{prop:values of chi-prop-alg}For any $k\ge1$ and non-trivial
words $w_{1},\ldots,w_{k}\in\F$, the possible values of $\chialg(w_{1},\ldots,w_{k})$
are $0,-1,\ldots,\chi(\F)$ or $-\infty$. Moreover,
\begin{enumerate}
\item $\chialg(w_{1},\ldots,w_{k})=-\infty$ if and only if the algebraic
part in the algebraic-free decomposition of $\eta_{w_{1},\ldots,w_{k}}$
is \emph{not} proper algebraic, if and only if for some $i$, $w_{i}$
is primitive in $\F$ with $\F=H*\langle w_{i}\rangle$, and for all
$j\ne i$ the word $w_{j}$ is conjugate into $H$.
\item \label{enu:constant term in E=00005Bstd std=00005D}$\chialg(w_{1},\ldots,w_{k})=0$
if and only if for every $i$, the word $w_{i}$ is either a proper
power or conjugate into the same cyclic subgroup as some $w_{j}$
with $j\ne i$.
\end{enumerate}
\end{prop}

The condition for $\chialg(w_{1},\ldots,w_{k})=-\infty$ is called
'strongly reducible' in \cite[Def.~2.9]{wilton2018essential}.
\begin{proof}
Consider the unique decomposition of $\eta_{w_{1},\ldots,w_{k}}$
as 
\[
\Gamma_{w_{1},\ldots,w_{k}}\stackrel{\eta_{\alg}}{\longrightarrow}\Sigma\stackrel{\eta_{\mathrm{free}}}{\longrightarrow}\Omega
\]
to an algebraic morphism and a free morphism, guaranteed by \eqref{eq:algebraic-free decomposition}.
If $(\gamma_{1},\gamma_{2})\in\algdecomp(\eta_{w_{1},\ldots,w_{k}})$
then there exists $\gamma_{2}'$ with $(\gamma_{1},\gamma_{2}')\in\algdecomp(\eta_{\alg})$
\cite[Thm.~4.9]{hanany2020word}. In this case, if $\eta_{\alg}$
is not a proper algebraic extension, then neither is $\gamma_{1}$.
Hence, if $\eta_{\alg}$ is \emph{not} a proper algebraic extension
if and only if $w_{1},\ldots,w_{k}$ admit no proper algebraic extensions,
namely, if and only if $\chialg(w_{1},\ldots,w_{k})=-\infty$. 

If $\eta_{\alg}$ is proper algebraic, the maximal Euler characteristic
in the definition \eqref{eq:def of chi-prop-alg} of\linebreak{}
$\chialg(w_{1},\ldots,w_{k})$ is at least 
\[
\chi\left(\Sigma\right)\ge\chi\left(\Omega\right)=\chi(\F),
\]
where the inequality is due to the fact that $\eta_{\mathrm{free}}$
is free and $\chi(\Gamma)\ge\chi(\Delta)$ in every free immersion
$\delta\colon\Gamma\to\Delta$ \cite[Prop.~4.3(3)]{hanany2020word}.
As any core graph has non-positive Euler characteristic, this proves
that $\chialg(w_{1},\ldots,w_{k})\in\left\{ 0,-1,\ldots,\chi(\F)\right\} \cup\left\{ -\infty\right\} $.
(All these values can be obtained, and even for single words -- see
\cite[Lem.~6.8]{puder2014primitive}.)

By Grushko's theorem, $\F$ splits canonically as $H_{1}*\ldots*H_{t}*J$
where for each $j=1,\ldots,t$, there is a non-empty subset $\emptyset\ne S_{j}\subseteq[k]$
so that $w_{i}$ is conjugate into $H_{j}$ for every $i\in S_{j}$,
each $H_{j}$ does not split freely relative to the set $\{w_{i}\,\mid\,i\in S_{j}\}$,
and $S_{1}\sqcup\ldots\sqcup S_{t}=[k]$. We claim that the associated
immersion $\gamma\colon\Gamma_{w_{1},\ldots,w_{k}}\to\Sigma$, where
$\Sigma$ is the multi core graph associated with the $\{H_{1}^{\F},\ldots,H_{t}^{\F}\}$,
is precisely $\eta_{\alg}$. Indeed, $\gamma$ is algebraic as it
is algebraic on any connected component of $\Sigma$ by definition,
the immersion $\Sigma\to\Omega$ is free by definition, and by uniqueness
this means that this decomposition coincides with the algebraic-free
decomposition of $\eta_{w_{1},\ldots,w_{k}}$. By the previous paragraph,
$\chialg(w_{1},\ldots,w_{k})=-\infty$ if and only if $\gamma$ is
\emph{not} a proper algebraic extension, which is precisely the case
if and only if some $S_{j}=\{w_{i}\}$ and $H_{j}=\langle w_{i}\rangle$,
meaning that $w_{i}$ is primitive in $\F$.

Finally, if $\Gamma_{w_{1},\ldots,w_{k}}\stackrel{\gamma}{\longrightarrow}\Delta\stackrel{\delta}{\longrightarrow}\Omega$
is a decomposition of $\eta_{w_{1},\ldots,w_{k}}$ with $\gamma$
algebraic and $\chi(\Delta)=0$, then each component of $\Delta$
corresponds to a cyclic subgroup of $\F$ and each $w_{i}$ is conjugate
into the subgroup corresponding to the cycle $\gamma(\Gamma_{w_{i}})$.
There are such proper algebraic extensions if and only if $w_{1},\ldots,w_{k}$
satisfy the conditions stated in the proposition for $\chialg(w_{1},\ldots,w_{k})=0$.
\end{proof}
\begin{rem}
Proposition \ref{prop:values of chi-prop-alg}\eqref{enu:constant term in E=00005Bstd std=00005D}
is analogous to similar results on the constant term of $\mathbb{E}[\chi(w_{1})\cdots\chi(w_{k})]$
where $\chi$ is the character of the defining representation in the
unitary group (see \cite[Thm.~2]{mingo2007second}, \cite[Thm.~4.1]{Radulescu06}
or \cite[Cor.~1.13]{MPunitary}), or in the orthogonal or compact
symplectic group \cite[Cor.~1.17]{MPorthsymp}. 
\end{rem}

\begin{cor}
\label{cor:prod of fix-1 strictly positive for large N and non-stronglly-reducible}
Let $w_{1},\ldots,w_{k}\in\F$ be non-trivial words, and assume that
$\chialg(w_{1},\ldots,w_{k})\ne-\infty$, namely, the algebraic part
in the algebraic-free decomposition of $\eta_{w_{1},\ldots,w_{k}}$
is proper algebraic. Then for every large enough $N$,
\[
\mathbb{E}_{S_{N}}\left[\prod_{i=1}^{k}\left(\fix\left(w_{i}\right)-1\right)\right]>0.
\]
\end{cor}

\section{The stable invariant $\protect\spf$\label{sec:spf}}

Let $G$ be a compact group and $\phi\in\irr(G)$ an irreducible character
of $G$. In this section we define the stable invariant $\spf$ and
prove some of its basic properties. This is a stable counterpart of
the invariant $\pi^{\phi}$ defined in \cite[Def.~1.5]{shomroni2023wreathI}.
As we explain below, this invariant generalizes $\spm$, the stable
mod-$m$ primitivity rank introduced in \cite[Def.~6.1]{PSh23}, and
to some extent also $\sp$, the stable primitivity rank, which was
originally defined in \cite[Def.~10.6]{wilton2024rational} and discussed
in \cite[\S\S4.1]{PSh23}. We first recall the definitions of the
latter two. The definitions use the following notion:
\begin{defn}
\label{def:efficient}Consider a diagram of core graphs 
\[
\xymatrix{P\ar@{->>}[d]^{\rho~~}\ar[r]^{b} & \Gamma\ar[d]^{f}\\
\Gamma_{w}\ar[r]^{\eta_{w}} & \Omega
}
\]
with $\rho$ a finite-degree covering map and $f$ an immersion. We
call the diagram \textbf{efficient} if $b$ is injective on some (and
therefore every) fiber of $\rho$ in $P$. This term is relevant whenever
we have a decomposition $b_{2}\circ b_{1}=\eta_{w^{\nu}}$ for some
partition $\nu\in{\cal P}$: call $b_{1}$ efficient if it is injective
on the fibers of $\Gamma_{w^{\nu}}\to\Gamma_{w}$.
\end{defn}

The notion of efficiency is slightly more involved when the map $f$
is not necessarily an immersion -- see \cite[Def.~3.1 and Rem.~4.2(4)]{PSh23}.
\begin{defn}
\label{def:spi}\cite[Def.~10.6]{wilton2024rational} Let $1\ne w\in\F$.
The \textbf{stable primitivity rank} of $w$ is 
\begin{equation}
\sp\left(w\right)\defi\inf\left\{ \frac{-\chi\left(\Gamma\right)}{\deg\left(\rho\right)}\,\middle|\,\begin{gathered}\xymatrix{P\ar@{->>}[d]^{\rho~~}\ar[r]^{b} & \Gamma\ar[d]^{f}\\
\Gamma_{w}\ar[r]^{\eta_{w}} & \Omega
}
\end{gathered}
~\begin{gathered}\mathrm{s.t.}~\Gamma~\mathrm{is~a~core~graph~and}~f~\mathrm{an~immersion,}\\
\rho~\mathrm{is~a~finite\text{-}degree~covering~map,}\\
b~\mathrm{is~proper~algebraic,\,and}\\
\,\mathrm{the~diagram~is~efficient}
\end{gathered}
\right\} .\label{eq:def of spi with core graphs}
\end{equation}
We also define $\sp\left(1_{\F}\right)\defi-1$. A diagram as in \eqref{eq:def of spi with core graphs}
which obtains the infimum, namely, such that $\sp(w)=\frac{-\chi(\Gamma)}{\deg(\rho)}$,
is called \textbf{$\sp$-extremal} for $w$.
\end{defn}

Of course, the diagram in \eqref{eq:def of spi with core graphs}
must commute. Let $w\ne1$. It is easy to see that $\sp(w)=\infty$
(namely, $w$ admits no diagrams as in \eqref{eq:def of spi with core graphs})
if and only if $w$ is primitive in $\F$. Recently, Wilton proved
that $\sp$-extremal diagrams exist for every non-primitive $1\ne w\in\F$
\cite{wilton2022rationality} (and therefore $\sp(w)\in\mathbb{Q}\cup\{\infty\}$).
Next, we define $\spm$:
\begin{defn}
\label{def:spm}\cite[Def.~6.1]{PSh23} Let $1\ne w\in\F$ and $1\ne m\in\mathbb{Z}_{\ge0}$.
The \textbf{mod-$m$ stable primitivity rank} of $w$ is 
\begin{equation}
\spm\left(w\right)\defi\inf\left\{ \frac{-\chi\left(\Gamma\right)}{\deg\left(\rho\right)}\,\middle|\,\begin{gathered}\xymatrix{P\ar@{->>}[d]^{\rho~~}\ar[r]^{b} & \Gamma\ar[d]^{f}\\
\Gamma_{w}\ar[r]^{\eta_{w}} & \Omega
}
\end{gathered}
~\begin{gathered}\mathrm{s.t.}~\Gamma~\mathrm{is~a~core~graph~and}~f~\mathrm{an~immersion,}\\
\rho~\mathrm{is~a~finite\text{-}degree~covering~map,}\\
\,\mathrm{the~diagram~is~efficient},\\
\text{and}\,\,\forall e\in E\left(\Gamma\right),~n_{b}\left(e\right)\in m\mathbb{Z}
\end{gathered}
\right\} ,\label{eq:spm with core graphs}
\end{equation}
where $P$ inherits an orientation from $\Gamma_{w}$, and $n_{b}\left(e\right)$
is the signed sum of times $e\in E\left(\Gamma\right)$ is covered
by $b$: a positive sign if the orientation of the edge in $P$ agrees
with that of $e$, and a negative sign otherwise.\footnote{Here we think of $\Gamma$ as a Serre graph -- see $\S\S$\ref{subsec:Core-graphs}.}
We also define $\spm\left(1_{\F}\right)\defi-1$. 
\end{defn}

To define $\spf$, we first discuss group-labeling of a graph $\Gamma$,
which may encode a homomorphism $\pi_{1}(\Gamma)\to G$. 
\begin{defn}
\label{def:G-labeling of a graph}Let $G$ be a compact group and
$\Gamma$ a finite Serre graph. A\textbf{ $G$-labeling} of $\Gamma$
is a map $\gamma\colon E\left(\Gamma\right)\to G$ satisfying $\gamma\left(\overline{e}\right)=\gamma\left(e\right)^{-1}$
for all $e\in E\left[\Gamma\right]$. A \textbf{Haar-random $G$-labeling}
of $\Gamma$ is a random $G$-labeling such that $\gamma\left(e\right)$
is Haar-random for every $e\in E\left(\Gamma\right)$, and the $G$-labels
on all edges are independent except for the constraint that $\gamma\left(\overline{e}\right)=\gamma\left(e\right)^{-1}$. 
\end{defn}

If $\Gamma$ is connected, we may identify $\pi_{1}(\Gamma)\cong\F_{r}$
for some $r$, where $\F_{r}$ is a rank-$r$ free group. The set
$\Hom(\pi_{1}(\Gamma),G)$ has a natural probability measure coming
from the Haar measure on $G^{r}$ (the direct product of $r$ copies
of $G$) through the isomorphism of sets $\Hom(\pi_{1}(\Gamma),G)\cong\Hom(\F_{r},G)\cong G^{r}$.
A Haar-random $G$-labeling can be used to encode a Haar-random homomorphism
$\pi_{1}(\Gamma)\to G$:
\begin{lem}
\label{lem:labeling to homomorphism to G}Let $\Gamma$ be a finite,
connected Serre graph with a fixed vertex $v\in V(\Gamma)$ and let
$G$ be a compact group. Every $G$-labeling of $\Gamma$ gives rise
to a well-defined homomorphism $\pi_{1}(\Gamma,v)\to G$ where for
any closed path $p$ at $v$ we map $[p]$ to the product in $G$
of the labels of the edges $p$ traverses. A Haar-random $G$-labeling
of $\Gamma$ gives rise to a Haar-random homomorphism $\pi_{1}(\Gamma,v)\to G$. 
\end{lem}

\begin{proof}
This is standard. See, for example, \cite[\S\S2.3]{hall2018ramanujan}.
\end{proof}
The $w$-measure on a compact group $G$ was defined in $\S$\ref{sec:Introduction}.
Equivalently, a $w$-random element of $G$ is obtained as the image
of $w$ through a Haar-random homomorphism $\varphi\colon\F\to G$,
or, equivalently, through a combination of $\eta_{w}\colon\Gamma_{w}\to\Omega$
and a Haar-random $G$-labeling of $\Omega$, as hinted to in the
following definition.
\begin{defn}
\label{def:phi-expectation of a morphism from w^sigma}Let $C$ be
an oriented cycle-graph, $\eta\colon C\to\Delta$ a graph morphism,
and $\gamma\colon E\left(\Delta\right)\to G$ a $G$-labeling of $\Delta$.
We denote by $\gamma\left(\eta\left(C\right)\right)$ the conjugacy
class in $G$ of the product of $G$-labels along the path of $C$
in $\Delta$. Namely, let $e_{1},\ldots,e_{m}$ be the edges of $C$
with orientation that agrees with that of $C$, ordered according
to their order in $C$ (so that $\iota\left(e_{i+1}\right)=\iota\left(\overline{e_{i}}\right)$
for $i=1,\ldots,m$ with $e_{m+1}\defi e_{1}$). Then 
\begin{equation}
\gamma\left(\eta\left(C\right)\right)\defi\left[\gamma\left(\eta\left(e_{1}\right)\right)\cdots\gamma\left(\eta\left(e_{m}\right)\right)\right]^{G}.\label{eq:conj class of a cycle under G-labeling}
\end{equation}
Now let $G$ be compact and $\phi\in\irr\left(G\right)$, let $P$
be a finite graph consisting of a disjoint union of oriented cycles
$C_{1},\ldots,C_{k}$, and let $\eta\colon P\to\Delta$ be a graph
morphism. The \textbf{$\eta$-expectation of $\phi$} is 
\[
\mathbb{E}_{\eta}\left[\phi\right]\defi\mathbb{E}_{\gamma}\left[\prod_{i=1}^{k}\phi\left(\gamma\left(\eta\left(C_{i}\right)\right)\right)\right].
\]
where $\gamma\colon E\left(\Delta\right)\to G$ is a Haar-random $G$-labeling
of $\Delta$. 
\end{defn}

We can now define $\spf$.
\begin{defn}
\label{def:spf}Let $G$ be a compact group and $\mathrm{triv}\ne\phi\in\irr\left(G\right)$
a non-trivial irreducible character. The stable $\phi$-primitivity
rank is
\begin{equation}
\spf\left(w\right)\defi\inf\left\{ \frac{-\chi\left(\Gamma\right)}{\deg\left(\rho\right)}\,\middle|\,\begin{gathered}\xymatrix{P\ar@{->>}[d]^{\rho~~}\ar[r]^{b} & \Gamma\ar[d]^{f}\\
\Gamma_{w}\ar[r]^{\eta_{w}} & \Omega
}
\end{gathered}
~\begin{gathered}\mathrm{s.t.}~\Gamma~\mathrm{is~a~core~graph~and}~f~\mathrm{an~immersion,}\\
\rho~\mathrm{is~a~finite\text{-}degree~covering~map,}\\
\,\mathrm{the~diagram~is~efficient},\\
\mathrm{and}~\mathbb{E}_{b}\left[\phi\right]\ne0
\end{gathered}
\right\} ,\label{eq:spf definition}
\end{equation}
where $P$ inherits the orientation from $\Gamma_{w}$.
\end{defn}

As mentioned above, $\spf$ generalizes $\spm$. Recall that we denote
$C_{m}\defi\{z\in S^{1}\,\mid\,z^{m}=1\}$ for $1\ne m\in\mathbb{Z}_{\ge0}$,
so $C_{0}=S^{1}$ and $C_{m}$ is finite non-trivial cyclic if $m\ge2$.
\begin{lem}
\label{lem:spm special case of spf}Let $1\ne m\in\mathbb{Z}_{\ge0}$,
$G=C_{m}$ and $\phi_{m}\in\irr\left(C_{m}\right)$ the standard character
given by the embedding $C_{m}=\left\{ z\in\mathbb{S}^{1}\,\middle|\,z^{m}=1\right\} \hookrightarrow\mathbb{S}^{1}$.
Then $\spf=\spm$. 

Moreover, in this case, in every diagram as in \eqref{eq:spf definition}
satisfying $\mathbb{E}_{b}\left[\phi\right]\ne0$, we have $\mathbb{E}_{b}\left[\phi\right]=1$.
\end{lem}

\begin{proof}
We need to show that in every diagram as in \eqref{eq:spf definition},
the condition $\mathbb{E}_{b}\left[\phi_{m}\right]\ne0$ is equivalent
to the condition that for every $e\in E(\Gamma)$, $n_{b}(e)\in m\mathbb{Z}$,
and that in this case $\mathbb{E}_{b}\left[\phi_{m}\right]=1$. Let
$\gamma\colon E(\Gamma)\to C_{m}$ be a Haar-random $C_{m}$-labeling
of $\Gamma$. Denote by $C_{1},\ldots,C_{k}$ the oriented cycles
in $P$. As $\phi_{m}$ is a linear (one-dimensional) character,
\begin{eqnarray*}
\mathbb{E}_{b}\left[\phi_{m}\right] & = & \mathbb{E}_{\gamma}\left[\prod_{i=1}^{k}\phi_{m}\left(\gamma\left(b\left(C_{i}\right)\right)\right)\right]=\mathbb{E}_{\gamma}\left[\prod_{e\in E(P)\,\mathrm{oriented}}\phi_{m}\left(\gamma\left(b\left(e\right)\right)\right)\right]=\\
 & = & \mathbb{E}_{\gamma}\left[\prod_{\left\{ e,\overline{e}\right\} \subseteq E(\Gamma)}\phi_{m}\left(\gamma\left(e\right)^{n_{b}(e)}\right)\right]=\prod_{\left\{ e,\overline{e}\right\} \subseteq E(\Gamma)}\mathbb{E}_{x\in C_{m}}\left[\phi_{m}\left(x^{n_{b}(e)}\right)\right].
\end{eqnarray*}
Each of the terms in the last product is equal to one if $n_{b}(e)\in m\mathbb{Z}$
and vanishes otherwise.
\end{proof}
The stable mod-$m$ primitivity rank $\spm$ for $m\ne1$ admits two
significant features not-necessarily possessed by $\spf$ for a non-linear
character $\phi$. First, we have the property from Lemma \ref{lem:spm special case of spf}
that $\mathbb{E}_{b}\left[\phi\right]\in\left\{ 0,1\right\} $, and,
in particular, $\mathbb{E}_{b}\left[\phi\right]$ is either zero or
a positive number. This means that its contributions in subsequent
formulas for the expected value of stable irreducible characters under
$w$-measure cannot cancel out. This is crucial in the proof of \cite[Thm.~1.5]{PSh23},
which is an analog of Theorem \ref{thm:spf bounds beta from below when supported on phi}
with $\spm$ in the stead of $\spf$ but with equality instead of
weak inequality.

Second, in the case of $\spm$, the condition that $\mathbb{E}_{b}\left[\phi\right]\ne0$
from \eqref{eq:spf definition} can be given a \emph{local} combinatorial
interpretation, namely, the condition that $n_{b}\left(e\right)\in m\mathbb{Z}$
for all $e\in E\left(\Gamma\right)$. This allows us to prove in \cite[Thm.~6.4]{PSh23}
the existence of $\spm$-extremal solutions, rationality and computability.
We do not know of such a local combinatorial condition in the case
of an arbitrary $\phi$. 
\begin{lem}
\label{lem:properties of E_phi}Let $\mathrm{triv}\ne\phi\in\irr\left(G\right)$
be a non-trivial irreducible character, and let $b\colon P\to\Gamma$
be an immersion from a union of oriented cycles $P$ to a core graph
$\Gamma$. 
\begin{enumerate}
\item \label{enu:E_phi(b)=00003DE_phi(b_alg)}If $b=b''\circ b'$ is a decomposition
of $b$ with $b''$ free. Then $\mathbb{E}_{b}\left[\phi\right]=\mathbb{E}_{b'}\left[\phi\right]$.
\item \label{enu:E_phi(b)=00003D0 if isomorphism component}If $b$ is an
isomorphism on some connected component of $\Gamma$, then $\mathbb{E}_{b}\left[\phi\right]=0$.
\end{enumerate}
\end{lem}

As usual, the assumption in Item \ref{enu:E_phi(b)=00003D0 if isomorphism component}
means that there is some connected component $\Gamma_{0}$ of $\Gamma$
with $b|_{b^{-1}\left(\Gamma_{0}\right)}\colon b^{-1}\left(\Gamma_{0}\right)\stackrel{\cong}{\to}\Gamma_{0}$
an isomorphism.
\begin{proof}
First, we may reduce both items to the case where $\Gamma$ is connected.
Indeed, assume that $\Gamma=\Gamma_{1}\sqcup\ldots\sqcup\Gamma_{k}$
are the connected components of $\Gamma$ and $b_{i}\colon b^{-1}\left(\Gamma_{i}\right)\to\Gamma_{i}$
is the restriction of $b$ to $b^{-1}\left(\Gamma_{i}\right)$. A
Haar random $G$-labeling of $\Gamma$ consists of independent, Haar
random $G$-labelings of each of the $\Gamma_{i}$ separately, so
$\mathbb{E}_{b}\left[\phi\right]=\mathbb{E}_{b_{1}}\left[\phi\right]\cdots\mathbb{E}_{b_{k}}\left[\phi\right]$.
By definition, $b''$ is free if and only if it is free when restricted
to every component of $\Gamma$. Thus, in both items it suffices to
assume that $\Gamma$ is connected, which we do for the rest of the
proof.

By Lemma \ref{lem:labeling to homomorphism to G}, a Haar-random $G$-labeling
of $\Gamma$ is equivalent to a Haar-random homomorphism $\alpha\colon\pi_{1}\left(\Gamma\right)\to G$.
If the (conjugacy classes of) elements of $\pi_{1}\left(\Gamma\right)$
represented by the image of $b$ are $y_{1},\ldots,y_{t}$, then $\mathbb{E}_{b}\left[\phi\right]=\mathbb{E}\left[\prod_{i=1}^{t}\phi\left(\alpha\left(y_{i}\right)\right)\right]$.
Now, assume that the assumed decomposition of $b$ is $P\stackrel{b'}{\longrightarrow}\Sigma\stackrel{b''}{\longrightarrow}\Gamma$
with $b''$ free. By definition, the connected components of $\Sigma$
correspond to subgroups $H_{1},\ldots,H_{m}\le\pi_{1}(\Gamma)$ with
$\pi_{1}\left(\Gamma\right)=H_{0}*H_{1}*\ldots*H_{m}$ for some $H_{0}\le\pi_{1}(\Gamma)$.
Every $y_{i}$ is conjugate into one of the $H_{j}$ ($j\ge1$). Clearly,
picking $\alpha$ is like picking independently Haar-random homomorphisms
$\alpha_{j}\colon H_{j}\to G$ for $j=0,1,\ldots,m$. Let $f\colon\left[t\right]\to\left[m\right]$
satisfy that $y_{i}$ is conjugate into $H_{f\left(i\right)}$. Then
\[
\mathbb{E}_{b}\left[\phi\right]=\mathbb{E}\left[\prod_{i=1}^{t}\phi\left(\alpha\left(y_{i}\right)\right)\right]=\prod_{j=1}^{m}\mathbb{E}\left[\prod_{i\in f^{-1}\left(j\right)}\phi\left(\alpha_{j}\left(y_{i}\right)\right)\right]=\mathbb{E}_{b'}\left[\phi\right],
\]
proving Item \ref{enu:E_phi(b)=00003DE_phi(b_alg)}.

For Item \ref{enu:E_phi(b)=00003D0 if isomorphism component}, assume
that $b\colon P\to\Gamma$ is an isomorphism with both $P$ and $\Gamma$
a single cycle. Then $\pi_{1}\left(\Gamma\right)\cong\mathbb{Z}=\left\langle x\right\rangle $
and the image of $x$ through a random $G$-labeling of $\Gamma$
is a Haar-random element of $G$. Here, $b\left(P\right)$ is (the
conjugacy class of) $x$, so it is (the conjugacy class of) a Haar-random
element of $G$. As $\phi\ne\mathrm{triv}$, we have $\mathbb{E}_{b}\left[\phi\right]=0$
by the orthogonality of characters.
\end{proof}
We conclude that in the definition of $\spf$ we may restrict to proper
algebraic morphisms.
\begin{cor}
\label{cor:spf with algebraic morphisms}Let $1\ne w\in\F$ and $\mathrm{triv}\ne\phi\in\irr\left(G\right)$.
Then 
\begin{equation}
\spf\left(w\right)=\inf\left\{ \frac{-\chi\left(\Gamma\right)}{\deg\left(\rho\right)}\,\middle|\,\begin{gathered}\xymatrix{P\ar@{->>}[d]^{\rho~~}\ar[r]^{b} & \Gamma\ar[d]^{f}\\
\Gamma_{w}\ar[r]^{\eta_{w}} & \Omega
}
\end{gathered}
~\begin{gathered}\mathrm{s.t.}~\Gamma~\mathrm{is~a~core~graph~and}~f~\mathrm{an~immersion,}\\
\rho~\mathrm{is~a~finite\text{-}degree~covering~map,}\\
b~\mathrm{is~proper~algebraic}\\
\,\mathrm{the~diagram~is~efficient},\\
\mathrm{and}\,\mathbb{E}_{b}\left[\phi\right]\ne0
\end{gathered}
\right\} .\label{eq:spf with algebraic}
\end{equation}
\end{cor}

\begin{proof}
Consider a diagram with $b\colon P\to\Gamma$ as in \eqref{eq:spf definition}.
It is enough to show that we can obtain a diagram as in \eqref{eq:spf with algebraic}
with $b'\colon P'\to\Gamma'$ and with the quotient $\frac{-\chi\left(\Gamma'\right)}{\deg\left(\rho'\right)}\le\frac{-\chi\left(\Gamma\right)}{\deg\left(\rho\right)}$.
Let let $P\stackrel{b_{\mathrm{alg}}}{\to}\Sigma\stackrel{b_{\mathrm{free}}}{\to}\Gamma$
be the algebraic-free decomposition of $b$ guaranteed by \eqref{eq:algebraic-free decomposition}.
By Lemma \ref{lem:properties of E_phi}\eqref{enu:E_phi(b)=00003DE_phi(b_alg)},
$\mathbb{E}_{b_{\mathrm{alg}}}\left[\phi\right]=\mathbb{E}_{b}\left[\phi\right]\ne0$.
By Lemma \ref{lem:properties of E_phi}\eqref{enu:E_phi(b)=00003D0 if isomorphism component},
$\Sigma$ contains no component on which $b_{\alg}$ is an isomorphism,
meaning that $b_{\alg}$ is proper algebraic. Hence 
\[
\begin{gathered}\xymatrix{P\ar@{->>}[d]^{\rho~~}\ar[r]^{b_{\mathrm{alg}}} & \Sigma\ar[d]^{f\circ b_{\mathrm{free}}}\\
\Gamma_{w}\ar[r]^{\eta_{w}} & \Omega
}
\end{gathered}
\]
satisfies all the assumptions of \eqref{eq:spf with algebraic}. By
\cite[Prop.~4.3(3)]{hanany2020word}, $\chi\left(\Gamma\right)\le\chi\left(\Sigma\right)$,
so $\frac{-\chi\left(\Sigma\right)}{\deg\left(\rho\right)}\le\frac{-\chi\left(\Gamma\right)}{\deg\left(\rho\right)}$.
\end{proof}
\begin{rem}
\label{rem:unified definition}When applied to the trivial character
(of an arbitrary compact group), the definition \eqref{eq:spf with algebraic}
recovers the stable primitivity rank $\sp(w)$. Indeed, $\mathbb{E}_{b}[\triv]$
is identically one, regardless of $b$. This yields a unified definition
to all the stable invariants discussed in the current section. 
\end{rem}

\begin{cor}
\label{cor:properties of spf}Let $G$ be a compact group, $\phi\in\irr\left(G\right)$
and $w\in\F$. Then,
\begin{enumerate}
\item \label{enu:spf AutF-invariant}$\spf$ is $\mathrm{Aut}\F$-invariant.
\item \label{enu:sp <=00003D spf}$\sp\left(w\right)\le\spf\left(w\right)$.
\item \label{enu:spf>=00003D1 for non-powers}If $w\ne1$ is a non-power
then $1\le\spf(w)$.
\end{enumerate}
\end{cor}

\begin{proof}
Item \ref{enu:spf AutF-invariant} follows from the fact that \eqref{eq:spf with algebraic}
can be given in completely group-theoretic terms. Indeed, it is easy
to see (and consult \cite[Rem.~4.2(5)]{PSh23}) that we may restrict
to diagrams as in \eqref{eq:spf with algebraic} with $\Gamma$ connected.
For any such diagram, pick an arbitrary vertex $v\in V\left(\Gamma\right)$
and let $H=f_{*}\left(\pi_{1}\left(\Gamma,v\right)\right)\le\F$.
Assume that $P=P_{1}\sqcup\ldots\sqcup P_{m}$ are the connected components
of $P$ with $\rho|_{P_{i}}\colon P_{i}\to\Gamma_{w}$ a degree-$n_{i}$
connected cover. For every $i=1,\ldots,m$, the path $b\left(P_{i}\right)$
corresponds to the conjugacy class in $H$ of $y_{i}=u_{i}w^{n_{i}}u_{i}^{-1}$
for some $u_{i}\in\F$. Now, $b$ being algebraic is equivalent to
$H$ being algebraic over the multiset of subgroups $\left\{ \left\langle y_{1}\right\rangle ,\ldots,\left\langle y_{m}\right\rangle \right\} $
(see $\S\S$\ref{subsec:Core-graphs}), and $b$ being, in addition,
\emph{proper} algebraic simply means that $m\ge2$ or $H\ne\langle y_{1}\rangle$.
The assumption that $b$ be injective on every fiber of $\rho$, is
equivalent to that for any $i,j\in\left\{ 1,\ldots,m\right\} $ and
$k\in\mathbb{Z}$, if $u_{i}w^{k}u_{j}^{-1}\in H$, then $i=j$ and
$n_{i}\mid k$. Finally, the quantity $\mathbb{E}_{\phi}[b]$ is equivalently
given by $\mathbb{E}\left[\prod_{i=1}^{m}\phi\left(\alpha\left(y_{i}\right)\right)\right]$,
where $\alpha\colon\pi_{1}(H)\to G$ is a Haar-random homomorphism.

Item \ref{enu:sp <=00003D spf} is immediate from Corollary \ref{cor:spf with algebraic morphisms}
and Definition \ref{def:spi}, and Item \ref{enu:spf>=00003D1 for non-powers}
follows from Item \ref{enu:sp <=00003D spf} and the fact that $\sp(w)\ge1$
for non-powers (see \cite[Thm.~4.4]{PSh23}).
\end{proof}

\section{Upper bounds for stable characters\label{sec:Lower-bounds-for-stable-chars}}

In the current section we prove Theorems \ref{thm:spi bounds beta from below when triv outside support}
and \ref{thm:spf bounds beta from below when supported on phi}. We
begin with two general formulas for the induction of the character
$\chi\boxtimes\triv$ from $G_{d}\times G_{N-d}$ to $G_{N}$. These
formulas are the content of Proposition \ref{prop:basic formula for induction f boxtensor triv}
and Corollary \ref{cor:basic formula for induction, with algebraic morphisms},
and will be used in the proofs of Theorems \ref{thm:spi bounds beta from below when triv outside support}
and \ref{thm:spf bounds beta from below when supported on phi}, as
well as in the proof of Theorem \ref{thm:formula for stable irreps of S_N}
in $\S$\ref{sec:Formula-for-stable-chars-of-S_N}. 

\subsection{A formula for induction of characters in $G\wr S_{N}$\label{subsec:formulas for induction}}

The arguments of the current $\S\S$\ref{subsec:formulas for induction}
use ideas from the forthcoming \cite{shomroni2023probabilisticHNC}.
Recall from $\S\S$\ref{subsec:reps of G wr S_N} that the conjugacy
class of an element $\left(v,\sigma\right)\in G_{N}$ is determined
by the numbers $\left(a_{\ell,C}\right)_{\ell,C}$, where $\ell\in[N]$
and $C$ is a conjugacy class in $G$, and $a_{\ell,C}$ is the number
of $\ell$-cycles in $\sigma$ such that the product of corresponding
elements in $v$ belongs to $C$.
\begin{defn}
\label{def:E_eta(f) }Let $\alpha\in S_{N}$ be a permutation, $1\ne w\in\F$
and $\eta\colon\Gamma_{w^{\alpha}}\to\Sigma$ an immersion of finite
graphs. This data gives rise to an \textbf{$\eta$-random} element
$(v,\sigma)\in G_{N}$, determined up to conjugation, as follows.
Let $\gamma\colon E(\Sigma)\to G$ be a Haar-random $G$-labeling.
Up to conjugation, the permutation $\sigma$ is deterministic and
given by $\alpha$. If $(x_{1}\,\ldots\,x_{\ell})$ is a cycle of
$\alpha$, consider the corresponding oriented cycle $D$ of $\Gamma_{w^{\alpha}}$
(isomorphic to $\Gamma_{w^{\ell}}$). The conjugacy class in $G$
corresponding to this cycle is then given by $\gamma(\eta(D))$, defined
in \eqref{eq:conj class of a cycle under G-labeling}.

If $f\colon G_{N}\to\mathbb{C}$ is a class function, we let 
\begin{equation}
\mathbb{E}_{\eta}\left[f\right]\defi\mathbb{E}\left[f\left(y\right)\right],\label{eq:E_eta(f)}
\end{equation}
where $y\in G_{N}$ is an $\eta$-random element. 
\end{defn}

\begin{rem}
\label{rem:connection between the two E_=00005Ceta}Definition \ref{def:E_eta(f) }
is related, of course, to Definition \ref{def:phi-expectation of a morphism from w^sigma}
of $\mathbb{E}_{\eta}[\phi]$. Note that $\phi$ is a class function
(an irreducible character, in fact) of $G$, whereas $f$ in \eqref{eq:E_eta(f)}
is a class function on $G_{N}$. When $f$ is the character $\chi^{\phi,\mu}$
from Proposition \ref{prop:irred chars of G wr S_N}\eqref{enu:formula when support of size 1}
and $\eta\colon\Gamma_{w^{\alpha}}\to\Sigma$ as in Definition \ref{def:E_eta(f) },
then 
\[
\mathbb{E}_{\eta}\left[\chi^{\phi,\mu}\right]=\chi^{\mu}\left(\alpha\right)\cdot\mathbb{E}_{\eta}\left[\phi\right].
\]
\end{rem}

The following Proposition is the main technical result in our forthcoming
proofs. 
\begin{prop}
\label{prop:basic formula for induction f boxtensor triv}Let $\chi\colon G_{d}\to\mathbb{C}$
be a character of $G_{d}$. Then
\[
\mathbb{E}_{w}\left[\mathrm{Ind}_{G_{d}\times G_{N-d}}^{G_{N}}\left(\chi\boxtimes\triv_{N-d}\right)\right]=\frac{1}{d!}\sum_{\alpha\in S_{d}}\sum_{\substack{\left(b,c\right)\in\decomp\left(\eta_{w^{\alpha}}\right)\colon\\
b~\mathrm{efficient}
}
}\mathbb{E}_{b}\left[\chi\right]\cdot\LB_{c}\left(N\right).
\]
\end{prop}

\begin{proof}
Along the proof denote 
\[
f_{N}\defi\mathrm{Ind}_{G_{d}\times G_{N-d}}^{G_{N}}\left(\chi\boxtimes\triv_{N-d}\right)\colon G_{N}\to\mathbb{C}.
\]
Repeating the argument in the proof of Proposition \ref{prop:irred chars of G wr S_N},
we obtain
\begin{equation}
f_{N}\left(v,\sigma\right)=\sum_{B\in\binom{[N]}{d}\colon\sigma(B)=B}\chi\left(v|_{B},\sigma|_{B}\right)\label{eq:base formula, step 1}
\end{equation}
for every $(v,\sigma)\in G_{N}$. Let $\left[N\right]{}_{d}$ denote
the set of ordered $d$-tuples of \emph{distinct} elements from $[N]$
(so $|[N]_{d}|=\frac{N!}{(N-d)!}$). As every set of size $d$ has
$d!$ possible orderings, if we go over such $d$-tuples instead of
over subsets of size $d$, \eqref{eq:base formula, step 1} becomes
\begin{equation}
f_{N}\left(v,\sigma\right)=\frac{1}{d!}\sum_{\alpha\in S_{d}}\sum_{B\in[N]_{d}\colon\sigma|_{B}=\alpha}\chi\left(v|_{B},\sigma|_{B}\right),\label{eq:base formula, step 2}
\end{equation}
where we write $\sigma|_{B}=\alpha$ to mean that the set of elements
in the tuple $B$ is invariant under $\sigma$, and that on the ordered
tuple, $\sigma$ coincides with $\alpha$.

Recall from Lemma \ref{lem:labeling to homomorphism to G} that a
Haar-random homomorphism $\F\to G$ can be realized by a Haar-random
$G$-labeling of the bouquet $\Omega$. Similarly, a Haar-random homomorphism
$\F\to G_{N}$ can be realized by a Haar-random $G$-labeling $\gamma\colon E(\hat{\Omega})\to G$
of a uniformly random $N$-cover $\hat{\Omega}$ of $\Omega$ \cite[Lem.~6.14]{PSh23}.
Here, an $N$-cover of $\Omega$ has $N$ labeled vertices $o_{i}=(o,i)$
for $i\in[N]$ and the $N$ edges above any petal of $\Omega$ correspond
to some (arbitrary) permutation in $S_{N}$ (see Footnote \ref{fn:random N cover}
and \cite[P.~9244]{hanany2020word}). A $w$-random element $(v,\sigma)\in G_{N}$
is then given as follows. First, for every $i\in[N]$ there is a unique
lift of the path in $\Omega$ corresponding to $w$ to a path in $\hat{\Omega}$
which begins at the vertex $o_{i}$. The permutation $\sigma\in S_{N}$
satisfies $\sigma(i)=j$ if and only if this lifted path ends at $o_{j}$.
Given $\sigma$, label by $v_{1},\ldots,v_{N}$ the $N$ vertices
of $\Gamma_{w^{\sigma}}$ representing the starting points of $w$,
so that there is a path reading $w$ from $v_{i}$ to $v_{j}$ if
and only if $\sigma(i)=j$.\footnote{In the way we explain the proof here, there is an immaterial issue:
we need to assume without loss of generality that $w$ is cyclically
reduced, for otherwise the starting point of $w$ is not part of the
Stallings core graphs $\Gamma_{w}$ or $\Gamma_{w^{\alpha}}$.} By the manner in which $\sigma$ was extracted from $\hat{\Omega}$,
there is an immersion $\eta\colon\Gamma_{w^{\sigma}}\to\hat{\Omega}$
mapping $v_{i}\mapsto o_{i}$ for all $i\in[N]$. The $w$-random
element of $G_{N}$ is the $\eta$-random element from Definition
\ref{def:E_eta(f) }.

Now let $(v,\sigma)\in G_{N}$ be a $w$-random element realized as
in the previous paragraph. Let $B=(b_{1},\ldots,b_{d})\in[N]_{d}$
be a $d$-tuple of distinct elements from $[N]$ and $\alpha\in S_{d}$.
The condition $\sigma|_{B}=\alpha$ means that there is an immersion
$b\colon\Gamma_{w^{\alpha}}\to\hat{\Omega}$ mapping $v_{i}\mapsto o_{b_{i}}$
for all $i\in[d]$ (here, $d$ vertices in $\Gamma_{w^{\alpha}}$
are labeled by $v_{1},\ldots,v_{d}$ as in the previous paragraph).
In particular, $b$ is injective on $\left\{ v_{1},\ldots,v_{d}\right\} $,
or, equivalently, efficient. Every such efficient $b$ factors uniquely
as $\Gamma_{w^{\alpha}}\stackrel{b_{1}}{\twoheadrightarrow}\Sigma\stackrel{\overline{b_{2}}}{\hookrightarrow}\hat{\Omega}$
with $b_{1}$ efficient and surjective, and $\overline{b_{2}}$ injective.
By symmetry, the distribution of the resulting random element $(v,\sigma)|_{B}$
does not depend on the particular $B$ but only on the morphism $b_{1}$,
and is given by a $b_{1}$-random element of $G_{d}$. So instead
of summing over all $B$'s, we may simply count the number of embeddings
of $\Sigma$ into $\hat{\Omega}$. Let $b_{2}=p\circ\overline{b_{2}}\colon\Sigma\to\Omega$.
By Lemma \ref{lem:L^B =00003D average number of inj lifts}, fixing
$b_{1}\colon\Gamma_{w^{\alpha}}\to\Sigma$ and $b_{2}\colon\Sigma\to\Omega$,
the average number of corresponding embeddings $\overline{b_{2}}\colon\Sigma\hookrightarrow\hat{\Omega}$
is precisely $\LB_{b_{2}}(N)$. We conclude that 
\[
\mathbb{E}_{w}\left[f_{N}\right]=\frac{1}{d!}\sum_{\substack{\alpha\in S_{d}}
}\sum_{\substack{\left(b_{1},b_{2}\right)\in\decomp\left(\eta_{w^{\alpha}}\right)\colon\\
b_{1}~\mathrm{efficient}
}
}\mathbb{E}_{b_{1}}\left[\chi\right]\cdot\LB_{b_{2}}\left(N\right).
\]
\end{proof}
We also need a version of the last proposition with algebraic morphisms.
The proof will rely on the following lemma, whose second part is a
generalization of Lemma \ref{lem:properties of E_phi}\eqref{enu:E_phi(b)=00003DE_phi(b_alg)}.
\begin{lem}
\label{lem:efficiency and E_eta invariant under free extensions}Let
$\alpha\in S_{d}$ and $\eta\colon\Gamma_{w^{\alpha}}\to\Delta$ be
an immersion of Stallings core graphs\footnote{Namely, an immersion of core graphs, where there is, in addition,
some $\eta_{2}\colon\Delta\to\Omega$ with $\eta_{2}\circ\eta=\eta_{w^{\alpha}}$}. Let $\eta=\eta''\circ\eta'$ be a decomposition with $\eta''$ free.
Then, 
\begin{enumerate}
\item \label{enu:efficiency invariant under free extensions}$\eta$ is
efficient if and only if $\eta'$ is efficient, and
\item \label{enu:E_eta=00005Bf=00005D invariant under free extensions}a
$\eta$-random element of $G_{d}$ is the same as a $\eta'$-random
element of $G_{d}$. In particular, $\mathbb{E}_{\eta}[f]=\mathbb{E}_{\eta'}[f]$
for every class function $f\colon G_{d}\to\mathbb{C}$.
\end{enumerate}
\end{lem}

\begin{proof}
\eqref{enu:efficiency invariant under free extensions} It is clear
that if $\eta$ is efficient then $\eta'$ is efficient. Conversely,
assume that $\eta$ is \emph{not} efficient. This means that it factors
through $q\colon\Gamma_{w^{\alpha}}\to\Gamma_{w^{\beta}}$ for some
$\beta\in S_{d'}$ and $d'<d$. Let $\gamma_{\mathrm{free}}\circ\gamma_{\alg}=\eta'$
be the algebraic-free decomposition of $\eta'$ guaranteed by \eqref{eq:algebraic-free decomposition},
so $(\eta''\circ\gamma_{\mathrm{free}})\circ\gamma_{\alg}$ is the
algebraic-free decomposition of $\eta$. Denote $\Sigma_{1}=\mathrm{codom}(\gamma_{\alg})$
and $\Sigma_{2}=\cod\left(\eta'\right)$. 
\[
\begin{gathered}\xymatrix{\Gamma_{w^{\alpha}}\ar@{->>}[dd]^{\rho~~}\ar[rr]^{\gamma_{\mathrm{alg}}}\ar@{->>}[rd]^{q~~} &  & \Sigma_{1}\ar[r]^{\gamma_{\mathrm{free}}} & \Sigma_{2}\ar[r]^{\eta''} & \Delta\ar[dd]\\
 & \Gamma_{w^{\beta}}\ar[ru]^{\exists!\theta}\\
\Gamma_{w}\ar[rrrr]^{\eta_{w}} &  &  &  & \Omega
}
\end{gathered}
\]
Note that $q$ is algebraic: every component of $\Gamma_{w^{\alpha}}$
represents a finite index subgroup in the component of its image under
$q$. By \cite[Thm.~4.9]{hanany2020word}, there is a (unique) morphism
$\theta\colon\Gamma_{w^{\beta}}\to\Sigma_{1}$ such that $\theta\circ q=\gamma_{\alg}$.
Hence $\eta'$ factors, too, through $q$, meaning that it is not
efficient. 

\eqref{enu:E_eta=00005Bf=00005D invariant under free extensions}
Recall from Definition \ref{def:E_eta(f) } that the permutation part
of a $\eta$-random element $(v,\sigma)\in G_{d}$ is determined by
the given $\alpha\in S_{d}$, so it is the same in a $\eta'$-random
element. The remainder of the proof is the same argument as in Lemma
\ref{lem:properties of E_phi}\eqref{enu:E_phi(b)=00003DE_phi(b_alg)},
so we only sketch it: The conjugacy class of $G$ corresponding to
a given cycle $D\subseteq\Gamma_{w^{\alpha}}$ is determined by the
image of $D$ through $\gamma\circ\eta$, where $\gamma\colon E(\Delta)\to G$
is a Haar-random $G$-labeling. Such a Haar-random labeling is equivalent
to sampling a Haar-random homomorphism from $\pi_{1}(\Delta')\to G$
for every connected component $\Delta'$ of $\Delta$. Finally, picking
a random homomorphism from $J=H_{1}*H_{2}$ is the same as picking
a random homomorphism from each of $H_{1}$ and $H_{2}$ separately
and independently.
\end{proof}
\begin{cor}
\label{cor:basic formula for induction, with algebraic morphisms}Let
$\chi\colon G_{d}\to\mathbb{C}$ be a character of $G_{d}$. Then
\[
\mathbb{E}_{w}\left[\mathrm{Ind}_{G_{d}\times G_{N-d}}^{G_{N}}\left(\chi\boxtimes\triv_{N-d}\right)\right]=\frac{1}{d!}\sum_{\alpha\in S_{d}}\sum_{\substack{\left(b,c\right)\in\algdecomp\left(\eta_{w^{\alpha}}\right)\colon\\
b~\mathrm{efficient}
}
}\mathbb{E}_{b}\left[\chi\right]\cdot L_{c}^{\alg}\left(N\right).
\]
\end{cor}

\begin{proof}
Given Proposition \ref{prop:basic formula for induction f boxtensor triv},
it is enough to show that for every $d\in\mathbb{Z}_{\ge0}$ and $\alpha\in S_{d}$, 

\[
\sum_{\substack{\left(b,c\right)\in\decomp\left(\eta_{w^{\alpha}}\right):\colon\\
b~\mathrm{efficient}
}
}\mathbb{E}_{b}\left[\chi\right]\cdot\LB_{c}\left(N\right)=\sum_{\substack{\left(b,c\right)\in\algdecomp\left(\eta_{w^{\alpha}}\right)\colon\\
b~\mathrm{efficient}
}
}\mathbb{E}_{b}\left[\chi\right]\cdot L_{c}^{\alg}\left(N\right).
\]
Indeed, by Lemma \ref{lem:efficiency and E_eta invariant under free extensions},
if $b=b_{2}\circ b_{1}$ is the algebraic-free decomposition of $b$
guaranteed by \eqref{eq:algebraic-free decomposition}, then $b$
is efficient if and only if $b_{1}$ is, and $\mathbb{E}_{b}[\chi]=\mathbb{E}_{b_{1}}[\chi]$.
Of course, $b$ is surjective if and only if $b_{2}$ is surjective.
Hence, 
\begin{eqnarray*}
\sum_{\substack{\left(b,c\right)\in\decomp\left(\eta_{w^{\alpha}}\right):\colon\\
b~\mathrm{efficient}
}
}\mathbb{E}_{b}\left[\chi\right]\cdot\LB_{c}\left(N\right) & = & \sum_{\substack{\left(b_{1},b_{2},c\right)\in\decompt\left(\eta_{w^{\alpha}}\right)\colon\\
b_{1}~\mathrm{efficient~and~algebraic,}\\
b_{2}~\mathrm{free}
}
}\mathbb{E}_{b_{1}}\left[\chi\right]\cdot\LB_{c}\left(N\right)\\
 & \stackrel{\mathrm{Lemma}~\ref{lem:L^alg in terms of L}}{=} & \sum_{\substack{\left(b_{1},\bar{c}\right)\in\algdecomp\left(\eta_{w^{\alpha}}\right)\colon\\
b_{1}~\mathrm{efficient}
}
}\mathbb{E}_{b_{1}}\left[\chi\right]\cdot L_{\bar{c}}^{\alg}\left(N\right).
\end{eqnarray*}
\end{proof}

\subsection{Proof of Theorems \ref{thm:spi bounds beta from below when triv outside support}
and \ref{thm:spf bounds beta from below when supported on phi}\label{subsec:Proof-of-lower bounds of stable chars}}

In the setting of Theorems \ref{thm:spi bounds beta from below when triv outside support}
and \ref{thm:spf bounds beta from below when supported on phi}, we
consider a stable irreducible character $\chi^{\arrm[\bullet]}$ of
$G_{\bullet}=G\wr S_{\bullet}$ such that $\triv\notin\supp(\arrm)$
(in Theorem \ref{thm:spf bounds beta from below when supported on phi}
we assume, moreover, that $\arrm$ is supported on a single $\phi\in\irr(G)$).
The main advantage of this assumption is that the irreducible character
$\chi^{\arrm[N]}$ is then given by induction of characters as in
Proposition \ref{prop:basic formula for induction f boxtensor triv}
and Corollary \ref{cor:basic formula for induction, with algebraic morphisms}.
\begin{lem}
\label{lem:if no triv, then stable is induction}Let $\arrm\colon\irr(G)\to{\cal P}$
satisfy $|\arrm|=d$ and $\triv\notin\supp(\arrm)$. Then whenever
$\chi^{\arrm[N]}$ is defined (namely, for $N\ge d$), 
\[
\chi^{\arrm[N]}=\mathrm{Ind}_{G_{d}\times G_{N-d}}^{G_{N}}\left(\chi^{\arrm}\boxtimes\triv_{N-d}\right).
\]
\end{lem}

\begin{proof}
Denote $\supp(\arrm)=\{\phi_{1},\ldots,\phi_{k}\}$ and $d_{i}=|\arrm\left(\phi_{i}\right)|$,
so that $d_{i}\ge1$ and $\sum d_{i}=d$. By definition, in the notation
of Proposition \ref{prop:irred chars of G wr S_N}\eqref{enu:formula when support of size 1},
$\chi^{\triv,\arrm[N](\triv)}=\chi^{\triv,(N-d)}$ is the trivial
character, so by Proposition \ref{prop:irred chars of G wr S_N}\eqref{enu:formula for character as induction},
\[
\chi^{\arrm\left[N\right]}=\mathrm{Ind_{G_{d_{1}}\times\ldots\times G_{d_{k}}\times G_{N-d}}^{G_{N}}\left(\chi^{\phi_{1},\arrm\left(\phi_{1}\right)}\boxtimes\ldots\boxtimes\chi^{\phi_{k},\arrm\left(\phi_{k}\right)}\boxtimes\triv\right).}
\]
By the same proposition 
\[
\chi^{\arrm}=\mathrm{Ind_{G_{d_{1}}\times\ldots\times G_{d_{k}}}^{G_{d}}\left(\chi^{\phi_{1},\arrm\left(\phi_{1}\right)}\boxtimes\ldots\boxtimes\chi^{\phi_{k},\arrm\left(\phi_{k}\right)}\right)},
\]
so
\[
\chi^{\arrm}\boxtimes\triv=\mathrm{Ind_{G_{d_{1}}\times\ldots\times G_{d_{k}}\times G_{N-d}}^{G_{d}\times G_{N-d}}\left(\chi^{\phi_{1},\arrm\left(\phi_{1}\right)}\boxtimes\ldots\boxtimes\chi^{\phi_{k},\arrm\left(\phi_{k}\right)}\boxtimes\triv\right)}.
\]
The result now follows by the transitivity of character induction.
\end{proof}
Recall the notation $\mathbb{E}_{\eta}[f]$ from Definition \ref{def:E_eta(f) }.
The following lemma is based on the same idea as Lemma \ref{lem:properties of E_phi}\eqref{enu:E_phi(b)=00003D0 if isomorphism component}.
\begin{lem}
\label{lem:E_eta=00005Bchi=00005D=00003D0 if isom on comp and triv not is supp}Let
$\arrm\colon\irr(G)\to{\cal P}$ satisfy $|\arrm|=d$ and $\triv\notin\supp(\arrm)$,
let $\alpha\in S_{d}$, $1\ne w\in\F$ and let $\eta\colon\Gamma_{w^{\alpha}}\to\Sigma$
be an immersion of finite graphs. If $\eta$ is an isomorphism on
one of the connected components of $\Sigma$, then $\mathbb{E}_{\eta}[\chi^{\arrm}]=0$.
\end{lem}

As above, the last assumption in the statement means there is a connected
component $\Sigma'$ of $\Sigma$ such that $\eta|_{\eta^{-1}(\Sigma')}\colon\eta^{-1}(\Sigma')\to\Sigma'$
is an isomorphism. In particular, $\text{\ensuremath{\Sigma'}}$ is
a cycle.
\begin{proof}
Let $D$ be the cycle of $\Gamma_{w^{\alpha}}$ so that $\eta$ is
an isomorphism on $\eta(D)$. By Definition \ref{def:E_eta(f) },
an $\eta$-random element $(v,\sigma)$ of $G_{d}$ has $\sigma$
conjugate to the given $\alpha$. By Proposition \ref{prop:irred chars of G wr S_N}\eqref{enu:formula for value of induced character},
the value $\chi^{\arrm}$ is given as a sum over set-partitions $\overline{B}$
of $[d]$ which are invariant under $\sigma$. In particular, one
of these subsets must contain all the indices involved in $D$. This
subset corresponds to some character $\phi\in\supp(\arrm)$, and one
of the terms in the formula \eqref{eq:formula for induced char} is
$\phi$ evaluated on a product of the Haar-random $G$-labels of $\eta(D)$
(with a prescribed order). By definition, these labels are independent
of the remaining labels of $\Sigma$, and the fact that $\eta$ is
an isomorphism on $\eta(D)$ means that this term is the only one
where these labels play out. Hence, the expectation of this particular
term is a factor of the overall contribution of $\overline{B}$, and
as $\phi\ne\triv$, the expectation of this term vanishes (exactly
as in the proof of Lemma \ref{lem:properties of E_phi}\eqref{enu:E_phi(b)=00003D0 if isomorphism component}).
\end{proof}
\begin{proof}[Proof of Theorem \ref{thm:spi bounds beta from below when triv outside support}]
 Assume that $\arrm\colon\irr(G)\to{\cal P}$ satisfies $|\arrm|=d\ge1$
and $\triv\notin\supp(\arrm)$. We ought to show that $\sp(w)\le\beta(w,\chi^{\arrm[\bullet]})$.
By Lemma \ref{lem:if no triv, then stable is induction} and Corollaries
\ref{cor:basic formula for induction, with algebraic morphisms} and
\ref{cor:L-alg order of chi},

\begin{eqnarray}
\mathbb{E}_{w}\left[\chi^{\arrm[N]}\right] & = & \frac{1}{d!}\sum_{\alpha\in S_{d}}\sum_{\substack{\left(b,c\right)\in\algdecomp\left(\eta_{w^{\alpha}}\right)\colon\\
b~\mathrm{efficient}
}
}\mathbb{E}_{b}\left[\chi^{\arrm}\right]\cdot L_{c}^{\alg}\left(N\right).\nonumber \\
 & = & \frac{1}{d!}\sum_{\alpha\in S_{d}}\sum_{\substack{\Gamma_{w^{\alpha}}\stackrel{b}{\longrightarrow}\Sigma\stackrel{c}{\longrightarrow}\Omega:~c\circ b=\eta_{w^{\alpha}},~\mathrm{and}\\
b~\mathrm{efficient~and~algebraic}
}
}\mathbb{E}_{b}\left[\chi^{\arrm}\right]\cdot N^{\chi(\Sigma)}\left(1+O\left(N^{-1}\right)\right)\nonumber \\
 & = & \frac{1}{d!}\sum_{\alpha\in S_{d}}\sum_{\substack{\Gamma_{w^{\alpha}}\stackrel{b}{\longrightarrow}\Sigma\stackrel{c}{\longrightarrow}\Omega:~c\circ b=\eta_{w^{\alpha}},~\mathrm{and}\\
b~\mathrm{efficient~and~proper~algebraic}
}
}\mathbb{E}_{b}\left[\chi^{\arrm}\right]\cdot N^{\chi(\Sigma)}\left(1+O\left(N^{-1}\right)\right),\label{eq:eq in proof of spi as lower bound}
\end{eqnarray}
where the last equality is by Lemma \ref{lem:E_eta=00005Bchi=00005D=00003D0 if isom on comp and triv not is supp}:
$\mathbb{E}_{b}[\chi^{\arrm}]=0$ unless $b$ is \emph{proper} algebraic.
The value of $\mathbb{E}_{b}[\chi^{\arrm}]$ does not depend on $N$.
Hence, every non-vanishing summand in \eqref{eq:eq in proof of spi as lower bound}
is of order $N^{\chi(\Sigma)}$ where $b\colon\Gamma_{w^{\alpha}}\to\Sigma$
is proper algebraic and efficient, and $\alpha\in S_{d}$. Every such
morphism gives rise to a diagram satisfying all the constraints in
Definition \ref{def:spi} of $\sp(w)$, and which is of degree $d$.
Hence, by definition, $\sp(w)\le\frac{-\chi(\Sigma)}{d}$. Finally,
by Proposition \ref{prop:dim of chi^arrm}, $\dim(\chi^{\arrm[N]})=\Theta(N^{d})$,
so
\[
\mathbb{E}_{w}\left[\chi^{\arrm[N]}\right]=O\left(N^{-d\cdot\sp(w)}\right)=O\left(\dim\left(\chi^{\arrm[N]}\right)^{-\sp(w)}\right).
\]
Namely, $\beta(w,\chi^{\arrm})\ge\sp(w)$.
\end{proof}
\begin{proof}[Proof of Theorem \ref{thm:spf bounds beta from below when supported on phi}]
 Now assume that $\arrm\colon\irr(G)\to{\cal P}$ satisfies $\supp(\arrm)=\{\phi\}$
with $\triv\ne\phi\in\irr(G)$ and denote $\mu=\arrm(\phi)$ and $d=|\mu|$.
We ought to show that $\spf(w)\le\beta(w,\chi^{\arrm[\bullet]})$.
Recall from Proposition \ref{prop:irred chars of G wr S_N}\eqref{enu:formula when support of size 1}
that $\chi^{\arrm}=\chi^{\phi,\mu}$, and from Remark \ref{rem:connection between the two E_=00005Ceta}
that if $b\colon\Gamma_{w^{\alpha}}\to\Sigma$ for some $\alpha\in S_{d}$
then $\mathbb{E}_{b}\left[\chi^{\phi,\mu}\right]=\chi^{\mu}\left(\alpha\right)\cdot\mathbb{E}_{b}\left[\phi\right]$.
Hence, under the current assumptions, \eqref{eq:eq in proof of spi as lower bound}
becomes\footnote{For this proof we could have also used Proposition \ref{prop:basic formula for induction f boxtensor triv}
instead of (indirectly) Corollary \ref{cor:basic formula for induction, with algebraic morphisms},
as $\spf$ can be defined with more general morphisms which need not
be (proper) algebraic.}
\begin{eqnarray}
\mathbb{E}_{w}\left[\chi^{\arrm[N]}\right] & = & \frac{1}{d!}\sum_{\alpha\in S_{d}}\chi^{\mu}\left(\alpha\right)\sum_{\substack{\Gamma_{w^{\alpha}}\stackrel{b}{\longrightarrow}\Sigma\stackrel{c}{\longrightarrow}\Omega:~c\circ b=\eta_{w^{\alpha}},~\mathrm{and}\\
b~\mathrm{efficient~and~proper~algebraic}
}
}\mathbb{E}_{b}\left[\phi\right]\cdot N^{\chi(\Sigma)}\left(1+O\left(N^{-1}\right)\right).\label{eq:in proof of spf as lower bound}
\end{eqnarray}
In every non-vanishing term in the last sum, we have $\mathbb{E}_{b}[\phi]\ne0$,
so the corresponding morphism $b\colon\Gamma_{w^{\alpha}}\to\Sigma$
gives rise to a diagram satisfying all the constraints in Definition
\ref{def:spf} of $\spf(w)$ (and even of \eqref{eq:spf with algebraic}),
and which is of degree $d$. Hence, by definition, $\spf(w)\le\frac{-\chi(\Sigma)}{d}$.
We conclude the proof exactly as we concluded the previous one.
\end{proof}

\section{\label{sec:sp is profinite}The stable primitivity rank is profinite}

We now prove Theorem \ref{thm:spi profinite} stating that $\sp$
is a profinite invariant, namely, that it coincides on every two words
which induce the same distributions on any finite group. More specifically,
we show that if $w_{1}$ and $w_{2}$ induce the same distribution
on every group of the form $S_{m}\wr S_{N}$, then $\sp(w_{1})=\sp(w_{2})$.
We do it by showing that the infimum of $\beta(w,\chi)$ over the
set ${\cal J}$ defined in the statement of the theorem coincides
with $\sp(w)$.

Given $1\ne w\in\F,$ we would like to show that every diagram as
in Definition \ref{def:spi} can be traced in a contribution to some
Fourier coefficient of the $w$-measure in some compact group. The
formula \eqref{eq:eq in proof of spi as lower bound} gives something
very close to that: for every compact group $G$, when restricted
to some $\arrm\colon\irr(G)\to{\cal P}$ with $|\arrm|=d$ and $\triv\notin\supp(\arrm)$,
we get a summand for every diagram of degree $d$ which is valid in
Definition \ref{def:spi}. Moreover, this contribution, if not vanishing,
is of the 'right' order of magnitude, $N^{\chi(\Sigma)}$, so if the
diagram is $\sp$-extremal, it is of order $N^{-d\cdot\sp(w)}=O(\dim(\chi^{\arrm[N]})^{-\sp(w)})$.
Recall also that by \cite{wilton2022rationality}, $\sp$-extremal
diagrams exist for every $w$. 

The problem, of course, is with the 'coefficient' $\mathbb{E}_{b}[\chi^{\arrm}]$
in every summand: this coefficient might vanish, and even if it does
not vanish, different summands of the same extremal order $N^{-d\cdot\sp(w)}$
might cancel out. The point of the following proof is to find a specific
group $G$ (which will be $S_{m}$) and a specific $\arrm$ so that
the coefficients $\mathbb{E}_{b}[\chi^{\arrm}]$ are always strictly
positive.
\begin{proof}[Proof of Theorem \ref{thm:spi profinite}]
 Fix $1\ne w\in\F$. Recall that ${\cal J}=\bigsqcup_{m\ge2}{\cal J}_{m}$
and ${\cal J}_{m}$ denotes the set of all stable irreducible characters
of $S_{m}\wr S_{\bullet}$ corresponding to some $\arrm\colon\irr(S_{m})\to{\cal P}$
with $\triv\notin\supp(\arrm)\ne\emptyset$. We ought to show that
\[
\sp(w)=\inf_{\chi\in{\cal J}}\beta\left(w,\chi\right)
\]
and that the infimum is attained. Theorem \ref{thm:spi bounds beta from below when triv outside support}
guarantees that $\sp(w)\le\beta(w,\chi)$ for all $\chi\in{\cal J}$.
It remains to show there exist some $\chi\in{\cal J}$ with $\sp(w)=\beta(w,\chi)$.

Consider an $\sp$-extremal diagram for $w$, guaranteed to exist
by \cite{wilton2022rationality}.
\[
\xymatrix{P\ar@{->>}[d]^{\rho~~}\ar[r]^{b_{0}} & \Sigma_{0}\ar[d]^{f}\\
\Gamma_{w}\ar[r]^{\eta_{w}} & \Omega
}
\]
Denote by $d$ the degree of $\rho$, so $P=\Gamma_{w^{\alpha_{0}}}$
for some $\alpha_{0}\in S_{d}$, the morphism $b_{0}\colon\Gamma_{w^{\alpha}}\to\Sigma_{0}$
is efficient and proper algebraic, and $\sp(w)=\frac{-\chi(\Sigma_{0})}{d}$.

Let $m\in\mathbb{Z}_{\ge2}$. Recall that the irreducible complex
representations of $S_{m}$ are given by partitions of $m$. Let $\std_{m}\in\irr(S_{m})$
be the standard character, corresponding to the partition $(m-1,1$),
let $\mu=(d)\vdash d$, and define $\arrm\colon\irr(S_{m})\to{\cal P}$
by $\arrm(\std_{m})=\mu$ and $\arrm(\phi)=\emptyset$ for every $\phi\ne\std_{m}$.

Now consider the formula \eqref{eq:in proof of spf as lower bound}
for $\mathbb{E}_{w}[\chi^{\arrm[N]}]$. The character $\chi^{\mu}$
of $S_{d}$ is the trivial character, so the formula becomes
\begin{eqnarray}
\mathbb{E}_{w}\left[\chi^{\arrm[N]}\right] & = & \frac{1}{d!}\sum_{\alpha\in S_{d}}\sum_{\substack{\Gamma_{w^{\alpha}}\stackrel{b}{\longrightarrow}\Sigma\stackrel{c}{\longrightarrow\Omega}:~c\circ b=\eta_{w^{\alpha}},~\mathrm{and}\\
b~\mathrm{efficient~and~proper~algebraic}
}
}\mathbb{E}_{b}\left[\std_{m}\right]\cdot N^{\chi(\Sigma)}\left(1+O\left(N^{-1}\right)\right).\label{eq:in proof of profinite step 1}
\end{eqnarray}
Note that at least one of the summands in \eqref{eq:in proof of profinite step 1}
corresponds to the $\sp-$extremal diagram from above, namely $(b,c)=(b_{0},f)$
and $\Sigma=\Sigma_{0}$. Now consider any $\left(b,c\right)\in\algdecomp\left(\eta_{w^{\alpha}}\right)$
with $b\colon\Gamma_{w^{\alpha}}\to\Sigma$ efficient and proper algebraic
as in \eqref{eq:in proof of profinite step 1}. We claim that $\mathbb{E}_{b}[\std_{m}]>0$
for every large enough $m$. First, if $\Sigma_{1},\ldots,\Sigma_{t}$
are the connected components $\Sigma$ and $b_{1},\ldots,b_{t}$ the
corresponding restrictions of $b$, then $\mathbb{E}_{b}[\std_{m}]=\prod_{i=1}^{t}\mathbb{E}_{b_{i}}[\std_{m}]$.
So it is enough to prove the claim under the assumption that $\Sigma$
is connected. In this case, the different cycles of $\Gamma_{w^{\alpha}}$
correspond to a multiset of elements $u_{1},\ldots,u_{t}$ in the
free group $\pi_{1}\left(\Sigma\right)$ and 
\[
\mathbb{E}_{b}\left[\std_{m}\right]=\mathbb{E}_{S_{m}}\left[\prod_{j=1}^{t}\std_{m}\left(u_{j}\right)\right].
\]
The fact that $b$ is proper algebraic means that the assumptions
of Corollary \ref{cor:prod of fix-1 strictly positive for large N and non-stronglly-reducible}
are satisfied, so $\mathbb{E}_{b}[\std_{m}]>0$ for every large enough
$m$.

Finally, let $m$ be large enough so that $\mathbb{E}_{b}[\std_{m}]>0$
for each of the (finitely many) summands in \eqref{eq:in proof of profinite step 1}.
For such values of $m$ we get that the corresponding $\arrm\colon\irr(S_{m})\to{\cal P}$
satisfies 
\[
\mathbb{E}_{w}\left[\chi^{\arrm[N]}\right]=\Theta\left(N^{\chi\left(\Sigma_{0}\right)}\right)=\Theta\left(N^{-d\cdot\sp(w)}\right)=\Theta\left(\dim\left(\chi^{\arrm[N]}\right){}^{-\sp(w)}\right),
\]
namely, $\beta(w,\chi^{\arrm[\bullet]})=\sp(w)$.
\end{proof}
The last proof showed that for every $1\ne w\in\F$ and every large
enough $m$, there is some stable irreducible character $\chi^{\arrm[\bullet]}$
of $S_{m}\wr S_{\bullet}$, with $\arrm$ supported on $\std_{m}$,
such that $\beta(w,\chi^{\arrm[\bullet]})=\sp(w)$. Note that such
$\arrm$ satisfies the assumptions of Theorem \ref{thm:spf bounds beta from below when supported on phi},
so $\sp^{\std_{m}}(w)\le\beta(w,\chi^{\arrm[\bullet]})$. By Corollary
\ref{cor:properties of spf}\eqref{enu:sp <=00003D spf}, $\sp(w)\le\sp^{\std_{m}}(w)$.
We conclude the following.
\begin{cor}
For every $w\in\F$,
\[
\sp(w)=\inf_{m\in\mathbb{Z}_{\ge2}}\sp^{\std_{m}}\left(w\right)=\inf_{\phi}\spf\left(w\right),
\]
where $\std_{m}\in\irr(S_{m})$ is the standard $(m-1)$-dimensional
representation.\footnote{The proof shows, moreover, that for every fixed $w$, $\sp^{\std_{m}}(w)$
stabilizes for every large enough $m$ and becomes equal to $\sp(w)$.}
\end{cor}

\section{Formulas for stable characters of symmetric groups and its generalizations\label{sec:Formula-for-stable-chars-of-S_N}}

As shown in Lemma \ref{lem:if no triv, then stable is induction},
as long as $\arrm(\triv)=\emptyset$, the stable character $\chi^{\arrm}$
can be defined using induction of characters, and one can estimate
and bound the corresponding Fourier coefficients of $w$-measures
as done in Theorems \ref{thm:spi bounds beta from below when triv outside support},
\ref{thm:spf bounds beta from below when supported on phi} and \ref{thm:spi profinite}.
However, the restriction to such $\arrm$ leaves out many stable characters
and, in particular, \emph{all} non-trivial stable characters of the
symmetric groups $S_{\bullet}$. In the current section, relying on
our analysis above of the induces characters, we derive a formula
for all stable Fourier coefficients of $w$-measures on $S_{\bullet}$
and, more generally, on $G\wr S_{\bullet}$. We begin with the formula
in the special case of $S_{\bullet}$, which was stated as Theorem
\ref{thm:formula for stable irreps of S_N} above. In $\S$\ref{subsec:general formula for G wr S}
we state and prove the formula for the general case.
\begin{thm*}[Theorem \ref{thm:formula for stable irreps of S_N}]
 Let $d\in\mathbb{Z}_{\ge0}$, let $\mu\vdash d$ be a partition
and $1\ne w\in\F$ a non-power. For every $N\ge d+\mu_{1}$ we have
\begin{equation}
\mathbb{E}_{w}\left[\chi^{\mu\left[N\right]}\right]=\frac{1}{d!}\sum_{\sigma\in S_{d}}\chi^{\mu}\left(\sigma\right)\sum_{\substack{\left(\eta_{1},\eta_{2},\eta_{3}\right)\in\algdecompt\left(w^{\sigma}\to\F\right)\colon\\
\eta_{1}~\mathrm{is~efficient,~and}\\
\cod\left(\eta_{2}\right)~\mathrm{has~no~cycles}
}
}C_{\eta_{2}}^{\mathrm{alg}}\left(N\right).\label{eq:formula for E_w of irreps of S_n}
\end{equation}
\end{thm*}
In the statement, $w^{\sigma}\to\F$ is the morphism $\eta_{w^{\sigma}}\colon\Gamma_{w^{\sigma}}\to\Omega$
from Definition \ref{def:Gamma_w} and $C_{\eta_{2}}^{\alg}$ is the
Möbius inversion defined in $\S\S\S$\ref{subsec:M=0000F6bius-inversions}.
Let $\Gamma_{w^{\sigma}}\stackrel{\eta_{1}}{\to}\Sigma\stackrel{\eta_{2}}{\to}\Sigma'\stackrel{\eta_{3}}{\to}\Omega$
be an element of $\algdecompt\left(\eta_{w^{\sigma}}\right)$. By
``$\eta_{1}$ is efficient'' we mean that the diagram
\[
\xymatrix{\Gamma_{w^{\sigma}}\ar[d]_{\rho}\ar[r]^{\eta_{1}} & \Sigma\ar[d]^{\eta_{3}\circ\eta_{2}}\\
\Gamma_{w}\ar[r]_{\eta_{w}} & \Omega
}
\]
is efficient, namely, that $\eta_{1}$ is injective on every fiber
of $\rho$. By ``$\cod\left(\eta_{2}\right)$ has no cycles'', we
mean that $\Sigma'$ has no connected component which is a cycle-graph.
When $d=0$, the only permutation in $S_{0}$ is the empty permutation
$\sigma$, and then $\Gamma_{w^{\sigma}}$ is the empty graph\footnote{\label{fn:d=00003D0}Indeed, one can make the formula \eqref{eq:formula for E_w of irreps of S_n}
work also for $d=0$. Recall that $\chi^{\emptyset\left[N\right]}$
is the trivial character, so $\mathbb{E}_{w}[\chi^{\emptyset\left[N\right]}]=1$.
On the other hand, the only algebraic morphism from the empty graph
$\emptyset$ is the identity $\id\colon\emptyset\to\emptyset$. So
$\algdecompt\left(\emptyset\to\Omega\right)=\left\{ \left(\id,\id,\emptyset\to\Omega\right)\right\} $,
and by \eqref{eq:bounds on C^alg}, $C_{\emptyset\to\emptyset}^{\mathrm{alg}}\left(N\right)=N^{\chi\left(\emptyset\right)}=1$.}.

Throughout the proof of Theorem \ref{thm:formula for stable irreps of S_N},
we will use the notation
\begin{equation}
f_{w}\left(\nu\right)=f_{w}\left(\sigma\right)=\sum_{\substack{\left(\eta_{1},\eta_{2},\eta_{3}\right)\in\algdecompt\left(w^{\sigma}\to\F\right)\colon\\
\eta_{1}~\mathrm{is~efficient,~and}\\
\cod\left(\eta_{2}\right)~\mathrm{has~no~cycles}
}
}C_{\eta_{2}}^{\mathrm{alg}}\left(N\right),\label{eq:f_w(nu)}
\end{equation}
where $\nu\vdash d$ is the partition corresponding to the cycle structure
of $\sigma\in S_{d}$. Our goal is, therefore, to prove that $\mathbb{E}_{w}[\chi^{\mu\left[N\right]}]=\frac{1}{d!}\sum_{\sigma\in S_{d}}\chi^{\mu}\left(\sigma\right)f_{w}\left(\sigma\right)$.
\begin{rem}
\label{rem:about proper powers}For $w\in\F$ a proper power, $\sp\left(w\right)=0$
and Conjecture \ref{conj:spi as statistical invariant} holds (see
\cite[\S\S4.2]{PSh23}). Yet, the formula \eqref{eq:formula for E_w of irreps of S_n}
does not hold as is (Corollary \ref{prop:step 1 goal} below fails).
However, the formula and proof hold for every $w\ne1$, including
proper powers, if we replace the condition that $\mathrm{codom}(\eta_{2})$
has no cycles, with the condition that $\eta_{2}\circ\eta_{1}$ is
proper algebraic (the only change needed in the proof is in the definition
of $\Sigma_{0}$ in the proof of Corollary \ref{prop:step 1 goal}).
\end{rem}

\subsection{Sketch of proof and an illustration\label{subsec:Sketch-of-proof}}

We start with a sketch and an illustration of the proof of Theorem
\ref{thm:formula for stable irreps of S_N}. The proof consists of
the following three steps.
\begin{enumerate}
\item \label{enu:formula for Ind without cycles} A simple adjustment to
Corollary \ref{cor:basic formula for induction, with algebraic morphisms},
stated as Corollary \ref{cor:formula for ind with eta_1 efficient},
gives that $\mathbb{E}_{w}[\mathrm{ind}_{S_{d}\times S_{N-d}}^{S_{N}}\left(\chi^{\mu}\boxtimes\mathrm{triv}\right)]$
is equal to a sum as in the right hand side of \eqref{eq:formula for E_w of irreps of S_n},
except that we omit the restriction that $\cod\left(\eta_{2}\right)$
has no cycles (while keeping the condition that $\eta_{1}$ is efficient).
Using properties of $C^{\alg}$ from $\S\S\S$\ref{subsec:Useful-properties-of-Mobius},
we conclude that the summation from Corollary \ref{cor:basic formula for induction, with algebraic morphisms}
can be rewritten as a summation over decompositions where $\eta_{1}$
is efficient \emph{and }$\cod\left(\eta_{2}\right)$ has no cycles,
although now the sum involves also permutations in $S_{m}$ with $m<d$.
\item \label{enu:Pieri}(Proposition \ref{prop:pieri}) Using Pieri's formula
we decompose $\mathrm{ind}_{S_{d}\times S_{N-d}}^{S_{N}}\left(\chi^{\mu}\boxtimes\mathrm{triv}\right)$
into stable irreducible characters.
\item \label{enu:endgame}(Section \ref{subsec:Finishing-the-proof of S_n formula})
We show that the summation from Step \ref{enu:formula for Ind without cycles}
can be rewritten as a sum of terms that look like the right hand side
of \eqref{eq:formula for E_w of irreps of S_n}, taken for all stable
irreducible characters of $S_{\bullet}$ appearing in the formula
from Step \ref{enu:Pieri}. We can then conclude by induction on $d$.
\end{enumerate}
Let us illustrate these steps with the stable character $\chi^{\mu[\bullet]}$
when $\mu=\left(2,1\right)\vdash3$. Corollary \ref{cor:formula for ind with eta_1 efficient}
gives
\begin{align}
\mathbb{E}_{w}\left[\mathrm{Ind}_{S_{3}\times S_{N-3}}^{S_{N}}\left(\chi^{\left(2,1\right)}\boxtimes\mathrm{triv}\right)\right]=\frac{1}{3!}\sum_{\sigma\in S_{3}}\chi^{\left(2,1\right)}\left(\sigma\right)\sum_{\substack{\left(\eta_{1},\eta_{2},\eta_{3}\right)\in\algdecompt\left(w^{\sigma}\to\F\right)\colon\\
\eta_{1}~\mathrm{is~efficient}
}
}C_{\eta_{2}}^{\mathrm{alg}}\left(N\right)\nonumber \\
=\frac{1}{3}\left[\begin{gathered}\sum_{\substack{\left(\eta_{1},\eta_{2},\eta_{3}\right)\in\algdecompt\left(w^{\left(1,1,1\right)}\to\F\right)\colon\\
\eta_{1}~\mathrm{is~efficient}
}
}C_{\eta_{2}}^{\mathrm{alg}}\left(N\right)\end{gathered}
\right]-\frac{1}{3}\left[\begin{gathered}\sum_{\substack{\left(\eta_{1},\eta_{2},\eta_{3}\right)\in\algdecompt\left(w^{\left(3\right)}\to\F\right)\colon\\
\eta_{1}~\mathrm{is~efficient}
}
}C_{\eta_{2}}^{\mathrm{alg}}\left(N\right)\end{gathered}
\right].\label{eq:two sums}
\end{align}
Now consider the first sum in \eqref{eq:two sums}. The domain, $\Gamma_{w^{\left(1,1,1\right)}}$,
has three connected components, each of which isomorphic to $\Gamma_{w}$.
As $w$ is a non-power, the efficiency of $\eta_{1}$ guarantees that
$\eta_{1}$ is always an isomorphism when restricted to (the preimage
of) cycles in its codomain. By Step \ref{enu:formula for Ind without cycles},
we may ignore the decompositions where $\cod\left(\eta_{2}\right)$
contains a cycle and $\eta_{2}\circ\eta_{1}$ is not an isomorphism
on its preimage, and if $\eta_{2}\circ\eta_{1}$ is an isomorphism
on some component, we may ignore the component. The components on
which $\eta_{2}\circ\eta_{1}$ is an isomorphism consist of $0$,
$1$, $2$ or $3$ copies of $\Gamma_{w}$. There are three possibilities
for choosing 1 or 2 components, and one possibility for choosing $0$
or 3. Thus the first sum of \eqref{eq:two sums} is equal to

\[
\frac{1}{3}\left[f_{w}\left(\left(1,1,1\right)\right)+3\cdot f_{w}\left(\left(1,1\right)\right)+3\cdot f_{w}\left(\left(1\right)\right)+f_{w}\left(\emptyset\right)\right].
\]
In the second sum of \eqref{eq:two sums} the codomain $\Gamma_{w^{\left(3\right)}}$
has only one component, so it can be rewritten as
\[
-\frac{1}{3}\left[f_{w}\left(\left(3\right)\right)+f_{w}\left(\emptyset\right)\right].
\]
Combining the two we get
\begin{eqnarray}
\mathbb{E}_{w}\left[\mathrm{Ind}_{S_{3}\times S_{N-3}}^{S_{N}}\left(\chi^{\left(2,1\right)}\boxtimes\mathrm{triv}\right)\right] & = & \frac{1}{3}f_{w}\left(\left(1,1,1\right)\right)-\frac{1}{3}f_{w}\left(\left(3\right)\right)+f_{w}\left(\left(1,1\right)\right)+f_{w}\left(\left(1\right)\right).\label{eq:formula for Ind as a sum of f_w}
\end{eqnarray}
Using Pieri's rule, Step \ref{enu:Pieri} gives the following formula
for every large enough $N$:
\begin{equation}
\mathrm{Ind}_{S_{3}\times S_{N-3}}^{S_{N}}\left(\chi^{\left(2,1\right)}\boxtimes\mathrm{triv}\right)=\chi^{\left(2,1\right)\left[N\right]}+\chi^{\left(2\right)\left[N\right]}+\chi^{\left(1,1\right)\left[N\right]}+\chi^{\left(1\right)\left[N\right]}.\label{eq:Pieri - example}
\end{equation}
Let us write the formula \eqref{eq:formula for E_w of irreps of S_n}
that we want to prove for the characters showing up in \eqref{eq:Pieri - example}:
\begin{eqnarray}
\mathbb{E}_{w}\left[\chi^{\left(2,1\right)\left[N\right]}\right] & \stackrel{?}{=} & \frac{1}{3!}\sum_{\sigma\in S_{3}}\chi^{\left(2,1\right)}\left(\sigma\right)f_{w}\left(\sigma\right)=\frac{1}{3}f_{w}\left(\left(1,1,1\right)\right)-\frac{1}{3}f_{w}\left(\left(3\right)\right)\nonumber \\
\mathbb{E}_{w}\left[\chi^{\left(2\right)\left[N\right]}\right] & \stackrel{?}{=} & \frac{1}{2!}\sum_{\sigma\in S_{2}}\chi^{\left(2\right)}\left(\sigma\right)f_{w}\left(\sigma\right)=\frac{1}{2}f_{w}\left(\left(1,1\right)\right)+\frac{1}{2}f_{w}\left(\left(2\right)\right)\label{eq:4 equations}\\
\mathbb{E}_{w}\left[\chi^{\left(1,1\right)\left[N\right]}\right] & \stackrel{?}{=} & \frac{1}{2!}\sum_{\sigma\in S_{2}}\chi^{\left(1,1\right)}\left(\sigma\right)f_{w}\left(\sigma\right)=\frac{1}{2}f_{w}\left(\left(1,1\right)\right)-\frac{1}{2}f_{w}\left(\left(2\right)\right)\nonumber \\
\mathbb{E}_{w}\left[\chi^{\left(1\right)\left[N\right]}\right] & \stackrel{?}{=} & \frac{1}{1!}\sum_{\sigma\in S_{1}}\chi^{\left(1\right)}\left(\sigma\right)f_{w}\left(\sigma\right)=f_{w}\left(\left(1\right)\right).\nonumber 
\end{eqnarray}
Now, from \eqref{eq:formula for Ind as a sum of f_w} and \eqref{eq:Pieri - example}
we conclude that the \emph{sum} of the four equations in \eqref{eq:4 equations}
holds. The content of Step \ref{enu:endgame} is that this works for
all $\mu$. We are done by induction on $d$ as all the stable characters
in the formula \eqref{eq:Pieri - example} correspond to partitions
of $d'<3$, except for $\chi^{\mu\left[N\right]}$ itself.

\subsection{Step 1: A formula without cycles in the codomain of $\eta_{2}$ \label{subsec:A-formula-without cycles}}

We start with the following version of Corollary \ref{cor:basic formula for induction, with algebraic morphisms}.
\begin{cor}
\label{cor:formula for ind with eta_1 efficient}Let $d\in\mathbb{Z}_{\ge0}$,
$\mu\vdash d$ and $1\ne w\in\F$. Then for every $N\ge d$,
\[
\mathbb{E}_{w}\left[\mathrm{Ind}_{S_{d}\times S_{N-d}}^{S_{N}}\left(\chi^{\mu}\boxtimes\mathrm{\mathrm{triv}}\right)\right]=\frac{1}{d!}\sum_{\sigma\in S_{d}}\chi^{\mu}\left(\sigma\right)\sum_{\substack{\left(\eta_{1},\eta_{2},\eta_{3}\right)\in\algdecompt\left(w^{\sigma}\to\F\right)\colon\\
\eta_{1}~\mathrm{is~efficient}
}
}C_{\eta_{2}}^{\mathrm{alg}}\left(N\right).
\]
\end{cor}

\begin{proof}
Using the fact that $L_{\eta}^{\mathrm{alg}}=\sum_{(\eta_{1},\eta_{2})\algdecomp(\eta)}C_{\eta_{1}}^{\mathrm{alg}}$,
Corollary \ref{cor:basic formula for induction, with algebraic morphisms}
gives
\begin{equation}
\mathbb{E}_{w}\left[\mathrm{Ind}_{S_{d}\times S_{N-d}}^{S_{N}}\left(\chi^{\mu}\boxtimes\mathrm{\mathrm{triv}}\right)\right]=\frac{1}{d!}\sum_{\alpha\in S_{d}}\sum_{\substack{\left(\eta_{1},\eta_{2},\eta_{3}\right)\in\algdecompt\left(w^{\alpha}\to\F\right)\colon\\
\eta_{1}~\mathrm{is~efficient}
}
}\mathbb{E}_{\eta_{1}}\left[\chi^{\mu}\right]\cdot C_{\eta_{2}}^{\mathrm{alg}}\left(N\right).\label{eq:formula for ind with C and E_eta}
\end{equation}
Recall that the term $\mathbb{E}_{\eta_{1}}[\chi^{\mu}]$, introduced
in Definition \ref{def:E_eta(f) }, gives the expected value of $\chi^{\mu}$
on a random element $\left(v,\sigma\right)$ of $G\wr S_{d}$ with
the permutation $\sigma$ conjugate to $\alpha$ and $v$ sampled
according to a random $G$-labeling of $\mathrm{codom}(\eta_{1})$.
But in the current case, $G$ is the trivial group, so this random
element is deterministically the permutation $\alpha$ (up to conjugation)
and $\mathbb{E}_{\eta_{1}}[\chi^{\mu}]=\chi^{\mu}(\alpha)$, independently
of $\eta_{1}$.
\end{proof}
Recall the notation of $f_{w}$ from \eqref{eq:f_w(nu)}. For $\sigma\in S_{d}$
we denote by $\mathrm{cyc}\left(\sigma\right)$ the set of cycles
of $\sigma$.
\begin{prop}
\label{prop:step 1 goal}Let $d\in\mathbb{Z}_{\ge0}$, $\mu\vdash d$
and $1\ne w\in\F$ a non-power. Then for every $N\ge d$,
\[
\mathbb{E}_{w}\left[\mathrm{Ind}_{S_{d}\times S_{N-d}}^{S_{N}}\left(\chi^{\mu}\boxtimes\mathrm{\mathrm{triv}}\right)\right]=\frac{1}{d!}\sum_{\sigma\in S_{d}}\chi^{\mu}\left(\sigma\right)\sum_{A\subseteq\mathrm{cyc}\left(\sigma\right)}f_{w}\left(\sigma|_{A}\right).
\]
\end{prop}

\begin{proof}
Corollary \ref{cor:formula for ind with eta_1 efficient} gives the
same formula with $\sum_{A\subseteq\mathrm{cyc}\left(\sigma\right)}f_{w}\left(\sigma|_{A}\right)$
replaced by 
\[
\sum_{\substack{\left(\eta_{1},\eta_{2},\eta_{3}\right)\in\algdecompt\left(w^{\sigma}\to\F\right)\colon\\
\eta_{1}~\mathrm{is~efficient}
}
}C_{\eta_{2}}^{\mathrm{alg}}\left(N\right).
\]
Fix such a decomposition $(\eta_{1},\eta_{2},\eta_{3})$. Denote $\Sigma_{0}\sqcup\Sigma_{1}=\mathrm{codom}(\eta_{2})$
where $\Sigma_{0}$ denotes the union of components of $\mathrm{codom}(\eta_{2})$
which are\emph{ }cycles, and let $\Gamma_{i}=(\eta_{2}\circ\eta_{1})^{-1}(\Sigma_{i})\subseteq\Gamma_{w}$
for $i=0,1$. Let $A\subseteq\mathrm{cyc}\left(\sigma\right)$ denote
the cycles of $\sigma$ corresponding to the connected components
of $\Gamma_{1}$, and let $\overline{\eta_{2}}$ denote the restriction
of $\eta_{2}$ to $\eta_{1}\left(\Gamma_{1}\right)\to\Sigma_{1}$.
Each component of $\Gamma_{0}$ is mapped by $\eta_{2}\circ\eta_{1}$
and by $\eta_{1}$ to a cycle. As $\eta_{1}$ is efficient and $w$
a non-power, $\eta_{1}$ must be an isomorphism on these components.
If $\eta_{2}$ is not an isomorphism on these cycles, then $C_{\eta_{2}}^{\mathrm{alg}}\left(N\right)=0$
for $N\ge d$ by Lemmas \ref{lem:C^alg vanishes on non-isomorphisms of cycles}
and \ref{lem:All Mobius inversions are  multiplicative on connected components of codomain}.
If $\eta_{2}$ is, too, an isomorphism on these cycles, then by the
same lemmas, $C_{\eta_{2}}^{\mathrm{alg}}\left(N\right)=C_{\overline{\eta_{2}}}^{\mathrm{alg}}\left(N\right)$
for $N\ge d$. Of course, when restricted to $\Gamma_{1}$, $\eta_{1}$
is still efficient and algebraic and $\overline{\eta_{2}}$ is algebraic
with no cycles in its codomain. This completes the proof.
\end{proof}

\subsection{Step 2: A variant of Pieri's formula\label{subsec:pieri}}

Recall that every partition $\mu=\left(\mu_{1},\ldots,\mu_{\ell}\right)$
is identified with a Young diagram with rows of lengths $\mu_{1},\ldots,\mu_{\ell}$. 
\begin{defn}
\label{def:P^- of mu}For a partition $\mu$ denote by $P^{-}\left(\mu\right)$
the set of partitions obtained from $\mu$ by removing at most one
cell from every column.
\end{defn}

For example, $P^{-}\left(\left(3,2\right)\right)=\left\{ \left(3,2\right),\left(2,2\right),\left(3,1\right),\left(2,1\right),\left(3\right),\left(2\right)\right\} $.
\begin{prop}
\label{prop:pieri}Let $\mu\vdash d$ be a partition. For $N\ge\left|\mu\right|+\mu_{1}$,
\begin{equation}
\mathrm{Ind}_{S_{d}\times S_{N-d}}^{S_{N}}\left(\chi^{\mu}\boxtimes\mathrm{triv}\right)=\sum_{\nu\in P^{-}\left(\mu\right)}\chi^{\nu\left[N\right]}.\label{eq:Pieri with P-}
\end{equation}
\end{prop}

\begin{proof}
The following formula is sometimes called Pieri's formula, or a special
case of the Littlewood-Richardson rule (see \cite[\S~I.9]{macdonald1998symmetric}):
\begin{equation}
\mathrm{Ind}_{S_{d}\times S_{N-d}}^{S_{N}}\left(\chi^{\mu}\boxtimes\mathrm{triv}\right)=\sum_{\tau\in P_{N}^{+}\left(\mu\right)}\chi^{\tau},\label{eq:Pieri formula with P+}
\end{equation}
where $P_{N}^{+}\left(\mu\right)$ are all the partitions of $N$
obtained from $\mu$ by \emph{adding} at most one cell to every column.
But adding at most one cell to every column of $\mu$ is the same
as removing at most one cell from every column (so taking some $\nu\in P^{-}\left(\mu\right)$)
and then adding exactly one cell to every column, including all the
original columns of $\mu$ that may be empty in $\nu$, and the latter
is guaranteed when expanding $\nu$ to $\nu\left[N\right]$ when $N\ge\left|\mu\right|+\mu_{1}$.
(For smaller values of $N$, the right hand side of \eqref{eq:Pieri with P-}
is not even defined.)
\end{proof}
We will also need the following ``inverse'' of Pieri's formula.
\begin{lem}
\label{lem:Inverse Pieri}Let $0\le k\le d$ be integers, $\mu\vdash d$
and $\tau\in S_{k}$. For $\alpha\in S_{d-k}$ denote by $\tau\sqcup\alpha$
the permutation of $S_{d}$ obtained from $\left(\tau,\alpha\right)$
via the natural embedding $S_{k}\times S_{d-k}\hookrightarrow S_{d}$.
Then 
\begin{equation}
\frac{1}{\left(d-k\right)!}\sum_{\alpha\in S_{d-k}}\chi^{\mu}\left(\tau\sqcup\alpha\right)=\sum_{\nu\in P^{-}\left(\mu\right)\colon\left|\nu\right|=k}\chi^{\nu}\left(\tau\right).\label{eq:inverse pieri}
\end{equation}
\end{lem}

\begin{proof}
The irreducible representations of $S_{k}\times S_{d-k}$ are precisely
$\{\chi^{\nu}\boxtimes\chi^{\vartheta}\}_{\nu\vdash k,\vartheta\vdash d-k}$.
Denote by $m_{\nu,\vartheta}\in\mathbb{Z}_{\ge0}$ the multiplicity
of $\chi^{\nu}\boxtimes\chi^{\vartheta}$ in the restriction of $\chi^{\mu}$
to $S_{k}\times S_{d-k}$, namely, 
\[
\mathrm{Res}_{S_{k}\times S_{d-k}}^{S_{d}}\left(\chi^{\mu}\right)=\bigoplus_{\nu,\vartheta}m_{\nu,\vartheta}\cdot\left(\chi^{\nu}\boxtimes\chi^{\vartheta}\right).
\]
So 
\[
\frac{1}{\left(d-k\right)!}\sum_{a\in S_{d-k}}\chi^{\mu}\left(\tau\sqcup\alpha\right)=\sum_{\nu\vdash k,\vartheta\vdash d-k}m_{\nu,\vartheta}\chi^{\nu}\left(\tau\right)\cdot\frac{1}{\left(d-k\right)!}\sum_{\alpha\in S_{d-k}}\chi^{\vartheta}\left(\alpha\right).
\]
Unless $\chi^{\vartheta}$ is the trivial character, the average $\frac{1}{\left(d-k\right)!}\sum_{\alpha\in S_{d-k}}\chi^{\vartheta}\left(\alpha\right)$
vanishes. By Frobenius reciprocity,
\begin{eqnarray*}
\frac{1}{\left(d-k\right)!}\sum_{a\in S_{d-k}}\chi^{\mu}\left(\tau\sqcup\alpha\right) & = & \sum_{\nu\vdash k}m_{\nu,\mathrm{triv}_{d-k}}\chi^{\nu}\left(\tau\right)=\sum_{\nu\vdash k}\chi^{\nu}\left(\tau\right)\left\langle \mathrm{Res}_{S_{k}\times S_{d-k}}^{S_{d}}\left(\chi^{\mu}\right),\chi^{\nu}\boxtimes\mathrm{triv}\right\rangle _{S_{k}\times S_{d-k}}\\
 & = & \sum_{\nu\vdash k}\chi^{\nu}\left(\tau\right)\left\langle \chi^{\mu},\mathrm{Ind}_{S_{k}\times S_{d-k}}^{S_{d}}\left(\chi^{\nu}\boxtimes\mathrm{triv}\right)\right\rangle _{S_{d}}=\sum_{\nu\in P^{-}\left(\mu\right)\colon\left|\nu\right|=k}\chi^{\nu}\left(\tau\right),
\end{eqnarray*}
where the last equality is by Pieri's formula \eqref{eq:Pieri with P-}.
\end{proof}

\subsection{\label{subsec:Finishing-the-proof of S_n formula}Finishing the proof
of Theorem \ref{thm:formula for stable irreps of S_N}}

Recall the formula for $\mathbb{E}_{w}[\mathrm{Ind}_{S_{d}\times S_{N-d}}^{S_{N}}\left(\chi^{\mu}\boxtimes\mathrm{triv}\right)]$
from Corollary \ref{prop:step 1 goal}. The (not yet proven) formula
$\mathbb{E}_{w}[\chi^{\nu\left[N\right]}]=\frac{1}{d!}\sum_{\sigma\in S_{d}}\chi^{\nu}\left(\sigma\right)f_{w}\left(\sigma\right)$
in Theorem \ref{thm:formula for stable irreps of S_N} together with
Proposition \ref{prop:pieri} gives another (not yet proven) formula
for $\mathbb{E}_{w}[\mathrm{Ind}_{S_{d}\times S_{N-d}}^{S_{N}}\left(\chi^{\mu}\boxtimes\mathrm{triv}\right)]$.
Lemma \ref{lem:sum=00003Dsum} shows that, at least when summed over
all $\nu\in P^{-}\left(\mu\right)$ for every fixed $\mu$, these
two formulas are equal, so the corresponding sum of the formulas from
Theorem \ref{thm:formula for stable irreps of S_N} is true. 
\begin{lem}
\label{lem:sum=00003Dsum}Let $\mu\vdash d$. For every $N\ge d+\mu_{1}$,
$\mathbb{E}_{w}[\mathrm{Ind}_{S_{d}\times S_{N-d}}^{S_{N}}\left(\chi^{\mu}\boxtimes\mathrm{triv}\right)]$
is equal to the sum predicted by Proposition \ref{prop:pieri} and
Theorem \ref{thm:formula for stable irreps of S_N}, namely, 
\[
\mathbb{E}_{w}\left[\mathrm{Ind}_{S_{d}\times S_{N-d}}^{S_{N}}\left(\chi^{\mu}\boxtimes\mathrm{triv}\right)\right]=\sum_{\nu\in P^{-}\left(\mu\right)}\frac{1}{\left|\nu\right|!}\sum_{\tau\in S_{\left|\nu\right|}}\chi^{\nu}\left(\tau\right)f_{w}\left(\tau\right).
\]
\end{lem}

\begin{proof}
We use Corollary \ref{prop:step 1 goal} and aggregate the terms according
to $\tau=\sigma|_{A}$: 
\begin{eqnarray}
\mathbb{E}_{w}\left[\mathrm{Ind}_{S_{d}\times S_{N-d}}^{S_{N}}\left(\chi^{\mu}\boxtimes\mathrm{triv}\right)\right] & = & \frac{1}{d!}\sum_{\sigma\in S_{d}}\chi^{\mu}\left(\sigma\right)\sum_{A\subseteq\mathrm{cyc}\left(\sigma\right)}f_{w}\left(\sigma|_{A}\right)\nonumber \\
 & = & \frac{1}{d!}\sum_{k=0}^{d}\sum_{\tau\in S_{k}}f_{w}\left(\tau\right)\sum_{\sigma\in S_{d},~A\subseteq\mathrm{cyc}\left(\sigma\right)\colon\sigma|_{A}=\tau}\chi^{\mu}\left(\sigma\right).\label{eq:tmp1}
\end{eqnarray}
For a fixed $\tau\in S_{k}$, there are $\binom{d}{k}$ options for
the subset of $\left[d\right]$ covered by $A$. By symmetry, up to
a factor of $\binom{d}{k}$, we may assume this subset is always $\left[k\right]$,
so \eqref{eq:tmp1} is equal to
\[
\frac{1}{d!}\sum_{k=0}^{d}\sum_{\tau\in S_{k}}f_{w}\left(\tau\right)\binom{d}{k}\sum_{\alpha\in S_{d-k}}\chi^{\mu}\left(\tau\sqcup\alpha\right)\stackrel{\mathrm{Lemma}~\ref{lem:Inverse Pieri}}{=}\sum_{k=0}^{d}\frac{1}{k!}\sum_{\tau\in S_{k}}f_{w}\left(\tau\right)\sum_{\nu\in P^{-}\left(\mu\right)\colon\left|\nu\right|=k}\chi^{\nu}\left(\tau\right),
\]
which is what we need to prove.
\end{proof}

\begin{proof}[Proof of Theorem \ref{thm:formula for stable irreps of S_N}]
 The proof now proceeds by induction on $d$. We already explained
in Footnote \ref{fn:d=00003D0} above why the formula works for $d=0$.
Now assume that $d\ge1$ and that the theorem is true for all partitions
of $d'$ for all $d'<d$. Let $\mu\vdash d$. Then the formula \eqref{eq:formula for E_w of irreps of S_n}
holds for $\mu$ by Proposition \ref{prop:pieri}, Lemma \ref{lem:sum=00003Dsum},
and the fact that $\mu$ is the only partition in $P^{-}\left(\mu\right)$
with $\ge d$ cells.
\end{proof}

\subsection{A formula for general stable irreducible characters of $G\wr S_{\bullet}$\label{subsec:general formula for G wr S}}

Let $G$ be again any compact group. We end this section with a formula
for arbitrary stable irreducible characters of $G\wr S_{\bullet}$,
a formula which applies to $\chi^{\arrm[\bullet]}$ even when $\triv\in\supp(\arrm)$.
Recall the notation $\mathbb{E}_{\eta}[\chi]$ from Definition \ref{def:E_eta(f) }.
\begin{thm}
\label{thm:formula for stable irreps of G wr S_bullet}Let $d\in\mathbb{Z}_{\ge0}$
and let $\arrm\colon\irr(G)\to{\cal P}$ with $|\arrm|=d$, and let
$1\ne w\in\F$ be a non-power. For every $N\ge d+\arrm(\triv)_{1}$
we have 
\begin{equation}
\mathbb{E}_{w}\left[\chi^{\arrm\left[N\right]}\right]=\frac{1}{d!}\sum_{\sigma\in S_{d}}\sum_{\substack{\left(\eta_{1},\eta_{2},\eta_{3}\right)\in\algdecompt\left(w^{\sigma}\to\F\right)\colon\\
\eta_{1}~\mathrm{is~efficient,~and}\\
\cod\left(\eta_{2}\right)~\mathrm{has~no~cycles}
}
}\mathbb{E}_{\eta_{1}}\left[\chi^{\arrm}\right]\cdot C_{\eta_{2}}^{\mathrm{alg}}\left(N\right).\label{eq:formula of w-expectation of stable characters of G wr S}
\end{equation}
\end{thm}

The proof follows roughly the same lines of the proof of Theorem \ref{thm:formula for stable irreps of S_N}
in $\S\S$\ref{subsec:Sketch-of-proof}-\ref{subsec:Finishing-the-proof of S_n formula}.
We go over the steps of the proof and stress the adaptations required
for establishing the more general Theorem \ref{thm:formula for stable irreps of G wr S_bullet}.
\begin{itemize}
\item The general strategy of the proof is the same as sketched in $\S\S$\ref{subsec:Sketch-of-proof}.
For the sake of the proof we now denote for every finitely supported
$\arrm\colon\irr(G)\to{\cal P}$ and $\sigma\in S_{|\arrm|}$, 
\[
g_{w}\left(\arrm,\sigma\right)=\sum_{\substack{\left(\eta_{1},\eta_{2},\eta_{3}\right)\in\algdecompt\left(w^{\sigma}\to\F\right)\colon\\
\eta_{1}~\mathrm{is~efficient,~and}\\
\cod\left(\eta_{2}\right)~\mathrm{has~no~cycles}
}
}\mathbb{E}_{\eta_{1}}\left[\chi^{\arrm}\right]\cdot C_{\eta_{2}}^{\mathrm{alg}}\left(N\right),
\]
so we ought to show that $\mathbb{E}_{w}[\chi^{\arrm[N]}]=\frac{1}{d!}\sum_{\sigma\in S_{d}}g_{w}\left(\arrm,\sigma\right).$
\item Instead of Corollary \ref{cor:formula for ind with eta_1 efficient}
we need the formula \eqref{eq:formula for ind with C and E_eta}:
\begin{equation}
\mathbb{E}_{w}\left[\mathrm{Ind}_{G_{d}\times G_{N-d}}^{G_{N}}\left(\chi^{\arrm}\boxtimes\mathrm{\mathrm{triv}}\right)\right]=\frac{1}{d!}\sum_{\sigma\in S_{d}}\sum_{\substack{\left(\eta_{1},\eta_{2},\eta_{3}\right)\in\algdecompt\left(w^{\sigma}\to\F\right)\colon\\
\eta_{1}~\mathrm{is~efficient}
}
}\mathbb{E}_{\eta_{1}}\left[\chi^{\arrm}\right]\cdot C_{\eta_{2}}^{\mathrm{alg}}\left(N\right).\label{eq:step 0 for G wr S}
\end{equation}
We skip, for now, the analog of Proposition \ref{prop:step 1 goal}
-- it will be proven together with the analog of Lemma \ref{lem:sum=00003Dsum}
in Lemma \ref{lem:sum=00003Dsum for G wr S} below.
\item For the results around Pieri's rule, we first generalize our notation.
For $\arrm\colon\irr(G)\to{\cal P}$, denote by $P^{-}\left(\arrm\right)$
the set of partition-valued functions on $\irr(G)$ obtained from
$\arrm$ by removing at most one cell from every column of $\arrm(\triv)$
(and leaving the other values of $\arrm$ unchanged). Pieri's rule
for the groups $G\wr S_{N}$ (e.g., \cite[Cor.~4.8]{ingram2009wreath})
translates, as in the proof of Proposition \ref{prop:pieri}, to the
following. Let $\arrm\colon\irr(G)\to{\cal P}$ satisfy $|\arrm|=d$.
For $N\ge d+\arrm(\triv)_{1}$,
\begin{equation}
\mathrm{Ind}_{G_{d}\times G_{N-d}}^{G_{N}}\left(\chi^{\arrm}\boxtimes\mathrm{triv}\right)=\sum_{\arrn\in P^{-}\left(\arrm\right)}\chi^{\arrn\left[N\right]}.\label{eq:Pieri for G wr S}
\end{equation}
\item Lemma \ref{lem:Inverse Pieri} generalizes to the following statement
(again, with the same proof). Let $0\le k\le d$ be integers, $\arrm\colon\irr(G)\to{\cal P}$
satisfy $|\arrm|=d$ and $\tau\in G_{k}$. For $\alpha\in G_{d-k}$
denote by $\tau\sqcup\alpha$ the element of $G_{d}$ obtained from
$\left(\tau,\alpha\right)$ via the natural embedding $G_{k}\times G_{d-k}\hookrightarrow G_{d}$.
Then 
\begin{equation}
\int_{\alpha\in G_{d-k}}\chi^{\arrm}\left(\tau\sqcup\alpha\right)=\sum_{\arrn\in P^{-}\left(\arrm\right)\colon\left|\arrn\right|=k}\chi^{\arrn}\left(\tau\right),\label{eq:inverse Pieri for G wr S}
\end{equation}
where the integral is with respect to the Haar measure on $G_{d-k}$.
\end{itemize}
We now present a generalization of Lemma \ref{lem:sum=00003Dsum},
encapsulating also a generalization of Proposition \ref{prop:step 1 goal}.
\begin{lem}
\label{lem:sum=00003Dsum for G wr S}Let $\arrm\colon\irr(G)\to{\cal P}$
with $|\arrm|=d$. For every $N\ge d+\arrm(\triv)_{1}$, 
\[
\mathbb{E}_{w}\left[\mathrm{Ind}_{G_{d}\times G_{N-d}}^{G_{N}}\left(\chi^{\arrm}\boxtimes\mathrm{triv}\right)\right]=\sum_{\arrn\in P^{-}\left(\arrm\right)}\frac{1}{\left|\arrn\right|!}\,\,\sum_{\beta\in S_{\left|\arrn\right|}}g\left(\arrn,\beta\right).
\]
\end{lem}

\begin{proof}
Start with \eqref{eq:step 0 for G wr S}, and consider a specific
$(\eta_{1},\eta_{2},\eta_{3})\in\algdecompt(w^{\sigma}\to\F)$ with
$\eta_{1}$ efficient. We use the same notation as in the proof of
Proposition \ref{prop:step 1 goal}, and denote also by $\overline{\eta_{1}}$
the restriction of $\eta_{1}$ to $\Gamma_{1}\to\eta_{1}\left(\Gamma_{1}\right)$
and by $k$ the size of the union of cycles in $A$. Again, we notice
that $\eta_{1}$ must restrict to an isomorphism $\Gamma_{0}\to\eta_{1}(\Gamma_{0})$,
and that unless $\eta_{2}$ is also an isomorphism on $\eta_{1}(\Gamma_{0})\to\Sigma_{0}$,
then $C_{\eta_{2}}^{\alg}(N)=0$ for $N\ge d$ and we may ignore this
summand. 

We are left with the cases where $\eta_{2}\circ\eta_{1}$ is an isomorphism
on $\Gamma_{0}$. In such cases, as mentioned above, $C_{\eta_{2}}^{\mathrm{alg}}\left(N\right)=C_{\overline{\eta_{2}}}^{\mathrm{alg}}\left(N\right)$.
As for $\mathbb{E}_{\eta_{1}}[\chi^{\arrm}]$, note that each cycle
of $\sigma$ in $\Gamma_{0}$ is mapped to a Haar-random element of
$G$ by the Haar random $G$-labeling of $\eta_{1}(\Gamma_{0})$. 

Fix an integer $0\le k\le d$, a permutation $\beta\in S_{k}$ and
$(\overline{\eta_{1}},\overline{\eta_{2}},\overline{\eta_{3}})\in\algdecompt(w^{\beta}\to\F)$
with $\overline{\eta_{1}}$ efficient and $\mathrm{codom}(\overline{\eta_{2}})$
without cycles, and consider the set $E(\beta,\overline{\eta_{1}},\overline{\eta_{2}},\overline{\eta_{3}})$
of all the summands in \eqref{eq:step 0 for G wr S} such that $\sigma|_{A}=\beta$,
$\eta_{2}\circ\eta_{1}$ is an isomorphism on $\Gamma_{0}$, $\eta_{1}|_{\Gamma_{I}}=\overline{\eta_{1}}$
and $\eta_{2}|_{\eta_{1}(\Gamma_{1})}=\overline{\eta_{2}}$. By symmetry,
up to a factor of $\binom{d}{k}$, we may assume that the cycles in
$A$ permute the numbers in $[k]$. There are $(d-k)!$ possible permutations
$\sigma\in S_{d}$ with $\sigma|_{[k]}=\beta$, and each of these
corresponds to a unique summand in $E(\beta,\overline{\eta_{1}},\overline{\eta_{2}},\overline{\eta_{3}})$.
Moreover, averaging over all these $(d-k)!$ possibilities for $\sigma$,
we ``inflate'' the $\overline{\eta_{1}}$-random element $\tau$
of $G_{k}$ to the random element $\tau\sqcup\alpha_{d-k}\in G_{d}$,
where $\alpha_{d-k}\in G_{d-k}$ is a Haar-random element. 

Overall, by grouping together all summands in $E(\beta,\overline{\eta_{1}},\overline{\eta_{2}},\overline{\eta_{3}})$,
we obtain that \linebreak{}
$\mathbb{E}_{w}[\mathrm{Ind}_{G_{d}\times G_{N-d}}^{G_{N}}(\chi^{\arrm}\boxtimes\mathrm{triv})]$
is equal to 
\begin{eqnarray*}
 &  & \frac{1}{d!}\sum_{k=0}^{d}\binom{d}{k}\sum_{\beta\in S_{k}}\sum_{{\scriptscriptstyle \substack{\left(\overline{\eta_{1}},\overline{\eta_{2}},\overline{\eta_{3}}\right)\in\algdecompt\left(\eta_{w^{\beta}}\right)\colon\\
\overline{\eta_{1}}~\mathrm{is~efficient,~and}\\
\cod\left(\overline{\eta_{2}}\right)~\mathrm{has~no~cycles}
}
}}\left(d-k\right)!\cdot\mathbb{E}\left[\chi^{\arrm}\left(\tau\sqcup\alpha_{d-k}\right)\right]\cdot C_{\overline{\eta_{2}}}^{\mathrm{alg}}\left(N\right),\\
 & \stackrel{\eqref{eq:inverse Pieri for G wr S}}{=} & \sum_{k=0}^{d}\frac{1}{k!}\sum_{\beta\in S_{k}}\sum_{{\scriptscriptstyle \substack{\left(\overline{\eta_{1}},\overline{\eta_{2}},\overline{\eta_{3}}\right)\in\algdecompt\left(\eta_{w^{\beta}}\right)\colon\\
\overline{\eta_{1}}~\mathrm{is~efficient,~and}\\
\cod\left(\overline{\eta_{2}}\right)~\mathrm{has~no~cycles}
}
}}C_{\overline{\eta_{2}}}^{\mathrm{alg}}\left(N\right)\sum_{\arrn\in P^{-}\left(\arrm\right)\colon\,\left|\arrn\right|=k}\mathbb{E}_{\overline{\eta_{1}}}\left[\chi^{\arrn}\right]\\
 & = & \sum_{k=0}^{d}\frac{1}{k!}\sum_{\arrn\in P^{-}\left(\arrm\right)\colon\left|\arrn\right|=k}\,\,\sum_{\beta\in S_{k}}g\left(\arrn,\beta\right)=\sum_{\arrn\in P^{-}\left(\arrm\right)}\frac{1}{\left|\arrn\right|!}\,\,\sum_{\beta\in S_{\left|\arrn\right|}}g\left(\arrn,\beta\right).
\end{eqnarray*}
\end{proof}
The endgame is the same as in the proof of Theorem \ref{thm:formula for stable irreps of S_N}.
\begin{proof}[Proof of Theorem \ref{thm:formula for stable irreps of G wr S_bullet}]
 The proof now proceeds by induction on $d=|\arrm|$. The case $d=0$
corresponds to the trivial stable character, and the Theorem holds
as explained in Footnote \ref{fn:d=00003D0}. Now assume that $d\ge1$
and that the theorem is true for all partitions of $d'$ for all $d'<d$.
Consider some $\arrm$ with $|\arrm|=d$. Then the formula \eqref{eq:formula of w-expectation of stable characters of G wr S}
holds for $\arrm$ by \eqref{eq:Pieri for G wr S}, Lemma \ref{lem:sum=00003Dsum for G wr S},
and the fact that $\arrm$ is the only element in $P^{-}\left(\arrm\right)$
of size $\ge d$.
\end{proof}
\begin{cor}
\label{cor:E_w=00005Bchi=00005D  rational in N and beta rational}Fix
$w\in\F$, a compact group $G$ and a finitely supported $\arrm\colon\irr(G)\to{\cal P}$.
Then, 
\begin{enumerate}
\item there exist a rational function $f_{w,\arrm}\in\mathbb{C}(N)$ such
that $\mathbb{E}_{w}[\chi^{\arrm[N]}]=f_{w,\arrm}(N)$ for every large
enough $N$,\footnote{When $G$ is finite, this is also proven in \cite[Thm.~1.6]{shomroni2023wreathII}.}
and 
\item if $\left|\arrm\right|>0$ and $w\ne1_{\F}$, then $\beta(w,\chi^{\arrm[\bullet]})\in\mathbb{Q}_{\ge0}\cup\{\infty\}$.
\end{enumerate}
\end{cor}

(Of course, $\beta(1_{\F},\chi^{\arrm[\bullet]})=-1$.)
\begin{proof}
For the first item, note that in the formula in Theorem \ref{thm:formula for stable irreps of G wr S_bullet}
there are finitely many summands, and in each of them the term $\mathbb{E}_{\eta_{1}}[\chi^{\arrm}]$
is a constant in $\mathbb{C}$. The term $C_{\eta_{2}}^{\alg}(N)$
coincides with a rational function in $\mathbb{Q}(N)$ for every large
enough $N$ by \cite[Cor.~6.9 and Prop.~6.15]{hanany2020word} (here,
'large enough' depends on $\eta_{2}$).

For the second item, note that the degree of the rational function
corresponding to $C_{\eta_{2}}^{\alg}(N)$ is non-positive by Proposition
\ref{prop:bound on C^alg} (by 'degree' of a rational function $\frac{P}{Q}$
with $P,Q\in\mathbb{Q}[N]$ we mean here $\deg P-\deg Q$ or $-\infty$
for the zero function). Hence $\deg(f_{w,\arrm})\in\mathbb{Z}_{\le0}\cup\{-\infty$\}.
By the proof of Proposition \ref{prop:dim of chi^arrm}, the dimension
$\dim(\chi^{\arrm[N]})$ coincides with a polynomial $g_{\arrm}\in\mathbb{C}[N]$
of degree $\left|\arrm\right|$. Hence,
\[
\beta\left(w,\chi_{}^{\arrm[\bullet]}\right)=-\frac{\deg(f_{w,\arrm})}{\left|\arrm\right|}\in\mathbb{Q}_{\ge0}\cup\{\infty\}.
\]
 
\end{proof}
\bibliographystyle{alpha}
\bibliography{stable-invariants-bib}

\begin{thebibliography}{CSST14}

\bibitem[BS87]{broder1987second}
A.~Broder and E.~Shamir.
\newblock On the second eigenvalue of random regular graphs.
\newblock In {\em 28th Annual Symposium on Foundations of Computer Science
  ({SFCS} 1987)}, pages 286--294. IEEE, 1987.

\bibitem[Cal09]{calegari2009scl}
D.~Calegari.
\newblock {\em scl}, volume~20 of {\em MSJ Memoirs}.
\newblock Mathematical Society of Japan, Tokyo, 2009.

\bibitem[Cas23]{cassidy2023projection}
E.~Cassidy.
\newblock Projection formulas and a refinement of {S}chur-{W}eyl-{J}ones
  duality for symmetric groups.
\newblock preprint arXiv:2312.01839v2, 2023.

\bibitem[Cas24]{cassidy2024random}
E.~Cassidy.
\newblock Random permutations acting on $k$--tuples have near--optimal spectral
  gap for $k=$poly$(n)$.
\newblock preprint arXiv:2412.13941, 2024.

\bibitem[CEF15]{church2015fi}
T.~Church, J.~S. Ellenberg, and B.~Farb.
\newblock {FI}-modules and stability for representations of symmetric groups.
\newblock {\em Duke Math. J.}, 164(9):1833--1910, 2015.

\bibitem[CSST14]{ceccherini2014representation}
T.~Ceccherini-Silberstein, F.~Scarabotti, and F.~Tolli.
\newblock {\em Representation theory and harmonic analysis of wreath products
  of finite groups}, volume 410 of {\em London Mathematical Society Lecture
  Note Series}.
\newblock Cambridge University Press, 2014.

\bibitem[HHH08]{hora2008limits}
A.~Hora, T.~Hirai, and E.~Hirai.
\newblock Limits of characters of wreath products $\mathfrak{S}_n$({T}) of a
  compact group {T} with the symmetric groups and characters of
  $\mathfrak{S}_\infty$({T}) {II}: From a viewpoint of probability theory.
\newblock {\em J. Math. Soc. Japan}, 60(4):1187--1217, 2008.

\bibitem[HP23]{hanany2020word}
L.~Hanany and D.~Puder.
\newblock Word measures on symmetric groups.
\newblock {\em Int. Math. Res. Not. IMRN}, 2023(11):9221--9297, 2023.

\bibitem[HPS18]{hall2018ramanujan}
C.~Hall, D.~Puder, and W.~Sawin.
\newblock Ramanujan coverings of graphs.
\newblock {\em Adv. Math.}, 323:367--410, 2018.

\bibitem[HW16]{helfer2016counting}
J.~Helfer and D.~T. Wise.
\newblock Counting cycles in labeled graphs: the nonpositive immersion property
  for one-relator groups.
\newblock {\em Int. Math. Res. Not. IMRN}, 2016(9):2813--2827, 2016.

\bibitem[IJS09]{ingram2009wreath}
F.~Ingram, N.~Jing, and E.~Stitzinger.
\newblock Wreath product symmetric functions.
\newblock {\em Internat. J. Algebra}, 3(1):1--19, 2009.

\bibitem[KM02]{KM02}
I.~Kapovich and A.~Myasnikov.
\newblock Stallings foldings and subgroups of free groups.
\newblock {\em J. Algebra}, 248(2):608--668, 2002.

\bibitem[LW17]{louder2017stackings}
L.~Louder and H.~Wilton.
\newblock Stackings and the {W}-cycles conjecture.
\newblock {\em Canad. Math. Bull.}, 60(3):604--612, 2017.

\bibitem[Mac98]{macdonald1998symmetric}
I.~G. Macdonald.
\newblock {\em Symmetric functions and Hall polynomials}.
\newblock Oxford university press, 1998.

\bibitem[MP19]{MPunitary}
M.~Magee and D.~Puder.
\newblock Matrix group integrals, surfaces, and mapping class groups {I}:
  {$\mathrm{U}(n)$}.
\newblock {\em Invent. Math.}, 218(2):341--411, 2019.

\bibitem[MP24]{MPorthsymp}
M.~Magee and D.~Puder.
\newblock Matrix group integrals, surfaces, and mapping class groups {II}:
  $\mathrm{O}(n)$ and $\mathrm{Sp}(n)$.
\newblock {\em Math. Ann.}, 388(2):1437--1494, 2024.

\bibitem[M{\'S}S07]{mingo2007second}
J.~A. Mingo, P.~{\'S}niady, and R.~Speicher.
\newblock Second order freeness and fluctuations of random matrices: {II}.
  {U}nitary random matrices.
\newblock {\em Adv. Math.}, 209(1):212--240, 2007.

\bibitem[Nic94]{nica1994number}
A.~Nica.
\newblock On the number of cycles of given length of a free word in several
  random permutations.
\newblock {\em Random Structures Algorithms}, 5(5):703--730, 1994.

\bibitem[PP15]{PP15}
D.~Puder and O.~Parzanchevski.
\newblock Measure preserving words are primitive.
\newblock {\em J. Amer. Math. Soc.}, 28(1):63--97, 2015.

\bibitem[PS25]{PSh23}
D.~Puder and Y.~Shomroni.
\newblock Stable invariants of words from random matrices.
\newblock preprint, arXiv:2311.17733v3, 2025.

\bibitem[Pud14]{puder2014primitive}
D.~Puder.
\newblock Primitive words, free factors and measure preservation.
\newblock {\em Israel J. Math.}, 201(1):25--73, 2014.

\bibitem[R{\u a}d06]{Radulescu06}
F.~R{\u a}dulescu.
\newblock Combinatorial aspects of {C}onnes's embedding conjecture and
  asymptotic distribution of traces of products of unitaries.
\newblock In {\em Proceedings of the Operator Algebra Conference, Bucharest}.
  Theta Foundation, 2006.

\bibitem[Sho23]{shomroni2023wreathII}
Y.~Shomroni.
\newblock Word measures on wreath products {II}.
\newblock preprint, arXiv:2311.11316, 2023.

\bibitem[Sho25]{shomroni2023wreathI}
Y.~Shomroni.
\newblock Word measures on wreath products {I}.
\newblock {\em Israel J. Math.}, 2025.
\newblock appeared online, doi.org/10.1007/s11856-025-2776-4.

\bibitem[Sho26]{shomroni2023probabilisticHNC}
Y.~Shomroni.
\newblock Probabilistic {H}anna {N}eumann conjectures.
\newblock In preparation, 2026+.

\bibitem[SS19]{sam2019representations}
S.~V. Sam and A.~Snowden.
\newblock Representations of categories of {G}-maps.
\newblock {\em J. Reine Angew. Math.}, 2019(750):197--226, 2019.

\bibitem[Sta83]{stallings1983topology}
J.~R. Stallings.
\newblock Topology of finite graphs.
\newblock {\em Invent. Math.}, 71(3):551--565, 1983.

\bibitem[Tak51]{takahasi1951note}
M.~Takahasi.
\newblock Note on chain conditions in free groups.
\newblock {\em Osaka Math. J}, 3(2):221--225, 1951.

\bibitem[Wil18]{wilton2018essential}
H.~Wilton.
\newblock Essential surfaces in graph pairs.
\newblock {\em J. Amer. Math. Soc.}, 31(4):893--919, 2018.

\bibitem[Wil24]{wilton2024rational}
H.~Wilton.
\newblock Rational curvature invariants for 2-complexes.
\newblock {\em Proc. Roy. Soc. A}, 480(2296):20240025, 2024.

\bibitem[Wil25]{wilton2022rationality}
H.~Wilton.
\newblock Rationality theorems for curvature invariants of 2-complexes.
\newblock {\em Math. Ann.}, 392:3765--3796, 2025.

\end{thebibliography}

\noindent Doron Puder, School of Mathematical Sciences, Tel Aviv University,
Tel Aviv, 6997801, Israel\\
\texttt{doronpuder@gmail.com}~\\

\noindent Yotam Shomroni, School of Mathematical Sciences, Tel Aviv
University, Tel Aviv, 6997801, Israel\\
\texttt{yotam.shomroni@gmail.com}~\\

\end{document}